\documentclass[a4paper,12pt]{article}
\usepackage{amssymb,amsfonts,amsmath}
\usepackage{blindtext}
\usepackage[]{hyperref}
\usepackage[T1]{fontenc}
\usepackage{titlesec}
\usepackage{slashbox,booktabs}
\usepackage{amsthm}
\usepackage{empheq}
\usepackage[english]{babel}
\usepackage[utf8]{inputenc}
\usepackage{latexsym}
\usepackage{epsfig}
\usepackage[all]{xy}
\usepackage{tikz}
\usepackage{cancel}

\usepackage{amscd}
\usepackage{lipsum}

\newcounter{Chapcounter}

\newcommand{\chapter}[1]{%
    {\centering%
        \addtocounter{Chapcounter}{1}%
        \LARGE {\textbf{ Chapter \theChapcounter: ~\parbox{0.8\textwidth}{#1}}}%
    }%
    \addcontentsline{toc}{section}{ Chapter:~\theChapcounter~~ #1}%
}

\newcommand{\PP}{{\mathbb{P}}}
\newcommand{\ZZ}{{\mathbb{Z}}}
\newcommand{\CC}{{\mathbb{C}}}

\begin{document}

\theoremstyle{plain}
\newtheorem {theorem} {Theorem} [section] 
\newtheorem {lemma} {Lemma} [section]
\newtheorem {definition} {Definition} [section]
\newtheorem {corollary}  {Corollary}  [section]
\newtheorem {proposition} {Proposition} [section]
\newtheorem {remark} {Remark} [section]
\newtheorem*{question*}{Question}




\title{A version of Kapranov's Chow quotient and smooth 
moduli space of point configurations in $\PP^2$.}

\author{Lev Radzivilovsky}
 \renewcommand{\today}
 
\maketitle

\begin{abstract}
A new moduli space for configurations of $n$ ordered points in a projective plane, 
which is a version of Kapranov's "Chow quotient of Grassmanians"
is introduced. The new construction is a Chow quotient as well 
but with additional lines connecting pairs of marked points 
(inspired by the idea of blow up).
Both Kapranov's construction and the new construction provide an algebraic variety 
which is a compactification for the space of generic configurations of $n$ 
distinct points in projective plane.
The difference is, that in Kapranov's construction, the space for configurations of 6 
points in two-dimensional plane is not smooth; 
while with the new construction, the space for
configuration of 6 points in plane is smooth.
\end{abstract}
 
\newpage

\tableofcontents

\newpage

\medskip


\section{Overview.}

This work is inspired by Kapranov's construction [Kap], 
which he has called "Chow quotients of Grassmanians" (further investigated in [KT]).

The construction provides a moduli space of ordered configurations of $n$ points in $P^k$.
It is very elegant and has several nice properties. It is a compactification 
for the space of of ordered configurations of $n$ points in generic position, 
up to projective transformations. It is also an algebraic variety.

In the case of configurations of $n$ points in $P^1$, Kapranov's construction 
gives the previously known $\overline{M}_{0,n}$ - the space of stable rational nodal 
curves with $n$ marked points, which are trees of stable projective lines.
$\overline{M}_{0,n}$ was extensively investigated ([Kn], or the space $X_n$ in [Keel],
where its cohomology ring is investigated); and has
important applications (a good exposition is [FP]). 
It turns out that $\overline{M}_{0,n}$ is a smooth algebraic variety.

The next particular case of Kapranov's construction is the space of configurations
of $n$ marked points in $P^2$. We shall denote this particular case by $Y_n$.
Similarly to the case of projective line, a generic point of $Y_n$ is a generic 
configuration of $n$ ordered marked points on a projective plane up to projective transformations;
there are also non-generic points which are made of several non-generic but stable configuration
of $n$ marked points. The combinatorics of how precisely several non-generic but stable 
configurations are combined into one point of the space is non-obvious (or in Kapranov's terms,
how several orbits form one Chow orbit), but it follows from the simple rule: the homology class 
of each Chow orbit is the same, and non-generic Chow orbit is a limit point of generic orbits.

All Kapranov's "Chow quotients of Grassmanians" are projective algebraic varieties; 
however not all of them are smooth. $Y_5$ is still smooth, since by duality it is isomorphic to 
$\overline{M}_{0,5}$.
So the first particular case which might be non-smooth is $Y_6$, and as a matter of fact is non-smooth.

In this work a new moduli space $\widetilde{Y_n}$ is introduced, which is similar to Kapranov's 
Chow quotient of Grassmanians, and is based on the same technology of Chow quotients.
The new space is also a compactification for the moduli space of generic configurations
of $n$ points in the projective plane, up to projective transformations.
The idea of this construction is similar to blow-up; before we apply Chow quotient, 
we consider not just configurations of $n$ marked points in projective planes, 
but also the lines connecting the pairs of marked points.

The main result of this work is that $\widetilde{Y_6}$ - our construction of the moduli space 
for the case of configurations of 6-tuples of points in a projective plane is smooth.
Thus we have a moduli (at least for the case of six points) 
a moduli space which is a smooth projective variety.

The construction makes sense for any $n$, and for higher dimensions as well, 
and it always provides a compact projective variety (as well as original Kapranov's construction);
and one might hope it is smooth. However, even our proof for the smoothness of $\widetilde{Y_6}$ 
consists of heavy case-checking, and is not easily generalizable.

In section 2 we recall Kapranov's construction of Chow quotients.
Then we give a visual description for points of the moduli space (which is useful 
both for Kapranov's Chow quotient of Grassmannian and the new space which will be constructed). 
In $2.2$ we show singular points in $Y_6$, the moduli space of $6$-tuples of points in $\PP^2$ 
which gives the motivation for defining the new moduli space.
In the end of section 2, we recall some projective invariants which will later give us
coordinated functions on the new space.
The first type of invariant is the cross-ratio
(for quadruples of lines connecting one of the marked points to some other 4 marked points),
the second (introduced in the subsection 2.3).
The second type of invariants which will be later useful for some more degenerate case
is the triple ratio of six-tuples of points (introduced in the subsection 2.4).

In section 3 that, the definition of the new moduli space is introduced.
The idea for the new moduli space is to plug the blow-up idea into the Chow quotient.
Therefore, in addition to the marked points, we introduce the lines connecting the pairs
of marked points and only then apply the Chow quotient.
Subsection $3.2$ is dedicated to discussion of the algebraic structure of $\widetilde{Y_n}$,
then in $3.3$ we discuss $\widetilde{Y_5}$ which will further allow us to describe 
different 6-tuples, as 6-tuple is a 5-tuple with an extra point.

Sections 4 and 5 contain the proof of smoothness of $\widetilde{Y_6}$.
For the lack of elegant proof, 
we discuss how to classify all possible Chow orbits. In many types 
of configurations it is easy to show that several types of points of our moduli space 
the space is smooth. In some of the more degenerate types of configuration, 
we look at the local coordinates and show that the cotangent space at that point
is spanned by 4 differentials of local coordinates.

Section 4 is dedicated to listing of all the possible types of Chow orbits, while explaining 
why is it possible to omit classes of types of Chow orbits for which smoothness is easy to eaplain;
section 5 is dedicated to calculation of the dimension of cotangent spacee for the less obvious 
configuration using local coordinates.

\medskip
\textbf{Acknowledgements.} I want to thank my supervisor Evgenii Shustin,
for a lot of beautiful geometry, and Shachar Carmeli for a lot of moral support and
his constant willingness to discuss any kind of mathematics at any moment.
I want also to thank Yoav Krauz for useful mathematical discussions 
(especially for telling me about the triple ratio). 

\section{Background.}
\subsection{Chow quotient of Grassmanians.}
Kapranov has introduced the idea of "Chow quotients" in \cite{Kap}, 
and used them to build the space of Chow quotients of Grassmanians.
The space is very interesting, and it was further studied in \cite{KT}.
We recall the construction here for convenience.

An algebraic group $G$ acts on complex a projective variety $X$.
For each homology class $\delta$ in $X$, 
the set of algebraic cycles, representing $\delta$, is a projective 
algebraic variety, which will be denoted $h_\delta$. 

All generic orbits produce the same homology class $h_\delta$.
Some non-generic orbits produce a different homology class,
or even have lower dimension.
We consider the closure of the set of all orbits in $h_\delta$.

We shall focus on the case when $X$ is $(\PP^2)^n$,
and $G$ is $GL_3$. Geometrically, we think of configurations
of $n$ marked points up to projective transformations.
However, in \cite{Kap} the general case of $\PP^m$ is considered.

We would like to discuss, how to visualize a degenerate
point of Chow quotient on projective plane. 

Let us start with a simpler example, when $X$ is $(\PP^1)^n$,
and $G=GL_2$. We want to show, how one could recover 
the picture of "trees of projective lines", 
starting with Chow quotients. 
For simplicity, consider just $4$ points.
Assume that marked points $1$ and $2$ converge to the zero point 
$(0:1)$ at uniform rate. Until we reach the limit, one can 
rescale the picture, so that marked points $1$ and $2$ are 
fixed, and the other two points $3$ and $4$ converge to the infinite
point $(1:0)$. 	

Up to projective transformation, these two pictures are the same one,
equivalently, one could say that they belong to the same orbit.
Now, one could ask what is the limit point. On one of the pictures,
points $1$ and $2$ coincide; on another picture, points $3$ and $4$ 
coincide. These are precisely the "branches" of the tree. 

We could have taken more points and separate them into two classes,
but the idea is the same. Every limit point of moduli space, which
is degenerate, can be described as several orbits taken together,
which can be geometrically regarded as several different configurations
of projective lines with marked points. We shall call these configurations
\textbf{charts}. On different charts, different points coincide.
The description of which configurations can be combined with which
is well known (a tree). It can also be deduced from the definition of Chow quotient.
The other way to formulate it, that for any $4$ marked points,
on there is at least one chart where one could compute the cross-ratio,
and if there are two such charts, the result should be the same.

Let us generalize this idea for the case of $\PP^2$.
We shall start by an example of a degenerate point of the Chow quotient.
Let $P_1,P_2,...,P_k$ and $Q_1,...,Q_m$ be marked points.
If we perform a homothety with center at $0$ and coefficient $\epsilon$ to points
$P_1,P_2,...,P_k$ (to be specific, the points $P_i=(x_i:y_i:1)$ are 
replaced by $P_i=(\epsilon x_i : \epsilon y_i: 1)$, 
while points $Q_1,...,Q_m$ remain at their places.
Up to a projective transformation, it is the same as to perform a homothety
with coefficient $\frac{1}{\epsilon}$ to all points $Q_1,...,Q_m$, to be specific $Q_i=(u_i:v_i:1)$
are replaced by $(u_i : v_i : \epsilon )$, while
keeping the points $P_1,P_2,...,P_k$ at their places.
Take the limit when $\epsilon\to 0$.
We get two different configurations: in one, points $P_1,...,P_k$
collide at the origin, in the other, points $Q_1,...,Q_m$ arrange themselves 
along the infinite line.

From the point of view of Chow quotient, these two limit configurations are taken together
(they are not projectively equivalent, in other words they give different 
orbits, but they belong to the same Chow orbit).
Different orbits of $PGL_3$ might have different dimension, e. g. when all 
the marked points are on the same line. The orbits which 
have maximal dimension are the ones with stabilizer $\{Id\}$. For instance, 
if some $4$ marked points are in general position then the dimension 
of the orbit is maximal.
It turns out that the opposite is also true: If there are no $4$ points in general position 
then the configuration is unstable, in the sense that the orbit have lower dimension. 
The unstable configurations does not occur as components of Chow orbits.  

Let us apply directly the definition of Chow quotient in order to describe the Chow orbits.
By K{\"u}nneth formula, the homology of $(\PP^2)^n$ is a tensor product of $n$ copies of 
$H_*(\PP^2,\ZZ)$. Let $\alpha_i$ denote the homology class of a line in the $i$-th 
plane of the product and $\beta_i$ the fundamental class of this plane. 
Then $H_*((\PP^2)^n,\ZZ) \cong \ZZ\left\langle 1,\alpha_1,\beta_1\right\rangle \otimes ... \otimes \ZZ
\left\langle1,\alpha_n,\beta_n\right\rangle$. 
The homology class of the closure of a generic $PGL_3$-orbit is given by the 
formula 

\begin{align*}
\sum_{i_1<i_2<i_3<i_4	} \beta_{i_1} \otimes \beta_{i_2} \otimes\beta_{i_3} \otimes \beta_{i_4}
+\sum_{i_1<i_2<i_3 , i_4<i_5} \beta_{i_1} 
\otimes \beta_{i_2} \otimes\beta_{i_3} \otimes \alpha_{i_4}\otimes \alpha_{i_5} +\\ 
+\sum_{i_1<i_2 , i_3<i_4<i_5<i_6} \beta_{i_1} 
\otimes \beta_{i_2} \otimes\alpha_{i_3} \otimes \alpha_{i_4}\otimes \alpha_{i_5}\otimes \alpha_{i_6}+\\
+\sum_{i_1 , i_2<...<i_7} \beta_{i_1} 
\otimes \alpha_{i_2} \otimes ...\otimes \alpha_{i_6}\otimes \alpha_{i_7} 
+\sum_{i_1<i_2<...<i_7<i_8} \alpha_{i_1} 
\otimes \alpha_{i_2} \otimes...\otimes \alpha_{i_8}
\end{align*}	

(This is a special case of formula $(2.1.7)$ in \cite{Kap}.) In order to show that, one could work out the intersection number between
a generic orbit and each basic homology class of codimension $8$. 
For example, by stating that $4$ generic points are required 
to be taken by a matrix $A\in GL_3$ to some specific locations, 
one can compute a coefficient of 
$\beta_{i_1} \otimes \beta_{i_2} \otimes\beta_{i_3} \otimes \beta_{i_4}$,
and by stating that $3$ generic points are required to be taken to specific 
places and $2$ more generic points are to be sent to some specific lines 
one can compute the coefficient of $\beta_{i_1} 
\otimes \beta_{i_2} \otimes\beta_{i_3} \otimes \alpha_{i_4}\otimes \alpha_{i_5}$
and so on. Each $\alpha$ condition gives a linear equation on the matrix entries 
of the projective transformation; each $\beta$ condition gives two linear equations.
Given $8$ generic condition, the matrix is defined uniquely up to a factor,
and generically, we get a unique projective transformation.
Since the condition is transversal (because equations are linear and generically independent)
one sees that the coefficient in the formula for homology class is $1$.

In non-generic cases, we can get several different orbits of maximal 
dimension in one Chow orbit, but the same homological class, 
which means that on exactly one of the charts we have points $i_1, i_2, i_3, i_4$
forming a non-degenerate configuration - that correspond to the first summand 
in our formula, and more laws of that kind corresponding to different other
summands.

The relation between the Chow orbits is not arbitrary, but we were not able 
to formulate it in a simple way to express it for surfaces, 
similar to  "trees of projective lines" for curves.
One geometric way to see which configuration is connected to which is the following.
Consider an orbit - stable, but not generic. Perturb slightly it into 
a generic configuration. In that new configuration, each $3$ points are 
not on one line. Choose some $4$ points, which were not generic at the original 
configuration, and drag the configuration back to the original position, while
constantly performing projective transformation to keep those $4$ points fixed.
The limit configuration would be another configuration in the same Chow orbit. 
If there are some $4$ points which are degenerate in both configurations you had,
repeat this process once more for those $4$ points, etc. 

\subsection{Singular points of $Y_6$}
At first glance, Kapranov's construction looks so nice that one 
wouldn't expect it to have singular points. However, even with configurations of $6$ points
it has singularities. 

We shall show a simple geometric way to see the singular points of $Y_6$.
The points will be named $A,B,C,D,E,F$. First, consider the configuration of 
points $A,B,C$ forming a triangle, and points $D,E,F$ being generic points 
on the lines $BC,AC,AB$ respectively. Each configuration of that kind is
defined, up to projective transformation, by just one parameter: 
$c=\frac{AF\cdot BD\cdot CE}{FB\cdot DC\cdot EA}$, considered in oriented distances, which reminds of 
Ceva theorem, so we shall call it the \textbf{Ceva ratio}. It is invariant 
under projective transformations.	
Actually, in the Chow quotient, we see $4$ different charts in this configuration:

\begin{center}
            	\includegraphics[scale=0.5]{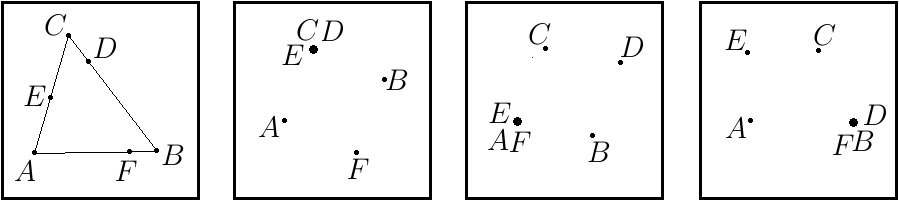}
\end{center}

because each time a triple of points are arranged along a line, on another chart
the other points collide, and in our case we have 3 such triples.
One could choose a projective transformation so that $A, B, C$ would 
become vertices of a specific regular triangle. We could also at the 
same time send the intersection point of $AD$ and $BE$ to the 
center of the triangle, and then $D$ and $E$ would become midpoints
of triangles sides, but the location of $F$ on $AB$ would be governed 
by the value of the Ceva parameter $c$.

Alternatively, we could make $D$ and $F$ be midpoints of the sides 
by some projective transformation, and the Ceva parameter would govern 
the location of $E$, or we could make $E$ and $F$ be the midpoints, 
and the Ceva parameter would govern a location of $D$, but all these
$3$ different charts are projectively equivalent.

Now, consider what happens at the limit when $c\to \infty$ or $c\to 0$.
Now these $3$ pictures become different, and we get the point of Chow 
quotient corresponding to $6$ charts:

\begin{center}
\includegraphics{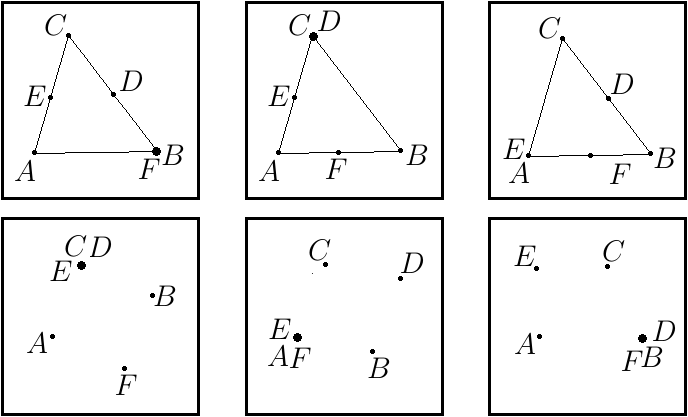}
\end{center}

If you keep the points $A,D,B,E$ in the left chart stable, 
and perturb points $C$ and $F$, you get some part of a neighborhood
of this Chow orbit which is similar to the neighborhood of 
$\CC^4$ (in the usual topology). 
In the configurations you get the cross ratio of lines 
$AB, AD, AF, AE$ is close to $0$, but the cross ratio of lines 
$BA, BD, BF, BE$ accepts all possible values (and is 
undefined at the original configuration).

However, if you perturb points $D$ and $A$ in the second chart
(from the left) the cross ratio of lines 
$BA, BD, BF, BE$ is close to $0$, and yet you get another part 
of the small neighborhood of the original Chow orbit, 
so that part of the neighborhood is different from 
the other part. 

So, we see directly that the neighborhood of this Chow orbit
does not look like a small part of $\CC^4$ but is more complicated.

\subsection{Cross-ratios.}

As we know from elementary projective geometry, for any $4$ distinct lines 
via a point in a projective plane the cross-ratio is well-defined.
Even if two of the lines coincide, it is still well-defined
(to be $0$, $1$ or $\infty$), but if $3$ or more lines coincide, they are 
no longer defined (however, if two disjoint pairs coincide, but not 
all lines are at the same place, the cross-ratio is still well-defined).

For a generic configuration of $n$ points, the cross-ratio can be defined 
for a set of $5$ points - the "central" point $O$ and $4$ additional points
$A,B,C,D$ as a cross-ratio of lines $OA,OB,OC,OD$. We shall denote it by 
$\left[A,B;C,D\right]_O$.
We use the convention 
\[\left[A,B;C,D\right]={\frac{AC}{BC}}:{\frac{AD}{BD}}\]
for the cross-ratio of $4$ points on one line, and similar ordering for 
$4$ lines through a point.
These functions are algebraic in coordinates of the points, 
and invariant under projective transformations. 

They can be extended to the entire $\widetilde{Y_n}$.
Indeed, on each generic Chow orbit the function is constant, 
and on all different configurations of a non-generic Chow orbit 
are being approximated by the same generic Chow orbit, so on each of 
different charts of the same Chow orbit the cross-ratio is either 
defined to be the same thing, or undefined. However, one can 
be sure that at least one configuration of that Chow orbit 
allows us to compute the cross-ratio. Indeed, take the configuration
on which the points $A,B,C,D$ are in general position. If point $O$
is distinct from all of them, then $\left[A,B;C,D\right]_O$ is obviously
well-defined. Indeed, no more than $2$ of the lines may coincide. 
In the case of $Y_n$, a problem might arise if point $O$ coincides with 
one of the points $A,B,C$ or $D$, since one of the lines might be not 
defined, but not in $\widetilde{Y_n}$, since
we may use the lines 
$\ell_{O,A},\ell_{O,B},\ell_{O,C},\ell_{O,D}$,
which are always well-defined, and in the chosen chart no more than $2$ of them coincide.

Those functions satisfy $2$ types of identities:
\begin{equation} \label{eq:1}
	\left[A,B;C,D\right]_E\cdot\left[A,B;D,E\right]_C\cdot\left[A,B;E,C\right]_D=1
\end{equation}
\begin{equation} \label{eq:2}
\left[A,B;C,D\right]_F\cdot\left[A,B;D,E\right]_F\cdot\left[A,B;E,C\right]_F=1
\end{equation}

Of course, care should be taken when dealing with zeroes and infinities; 
each cross-ratio should be considered as a fraction, in which the nominator or 
the denominator can be zero, but not simultaneously. The identity holds if product 
of nominators equals to product of denominators; so if one of cross-ratios is zero, 
another cross-ratio has to be infinity, and the third might be anything.

One might see this cross-ratio functions as projective coordinates.
For instance, if 
\[A=(1:0:0), B=(0:1:0), C=(0:0:1), D=(1:1:1), P=(x:y:z)\]
Then $\left[A,B;P,D\right]_C=\frac{y}{x}$.
Multiplying this identity with two similar identities for $\frac{z}{y}$ and $\frac{y}{x}$
we get \ref{eq:1}.

The identity \ref{eq:2}, we project all other points from $F$ 
to a projective line not passing through $F$, and we get a relation between
$5$ cross-ratios of $5$ points on a projective line and if we write them 
explicitly we get a product of things that cancel out.

\subsection{Triple ratios.}
Given 6 marked points $A,B,C,P,Q,R$ we may define another interesting projective invariant
which we shall call triple ratio. For that we have to assume points $A, B, C$ form a non-degenerate 
triangle, and lines $AP, BQ,$ and $CR$ are well-defined (which means $P$ is distinct from $A$ 
or that the line $AP$ is chosen anyway as in our construction). Let $P', Q',$ and $R'$
be the intersection points of $AP, BQ,$ and $CR$ with the respective sides 
of the triangle $BC, CA,$ and $AB$.
Then we define the triple ratio 
\[ \left\{A,B,C;P,Q,R\right\} = \frac{AR' \cdot BP' \cdot CQ' }{ R'B \cdot P'C \cdot Q'A} \] 
in oriented intervals (for example, in the classical Ceva theorem the triple ratio is 1, 
and in the classical Menelaus theorem the triple ratio is -1).

The triple ratio is a projective invariant in a generic case. 
If we choose any generic line which intersects the sides of the triangle at $X,Y$, and $Z$ respectively,
it is easy to see that 
\[ \left\{A,B,C;P,Q,R\right\} = -\left[A,B;R,Z\right]\cdot \left[B,C;P,X\right]\cdot \left[C,A;Q,Y\right] \] 
(the minus comes from Menelaus theorem).
Since the cross-ratio is a projective invariant, the triple ratio is also a projective invariant.

In the case when one of the lines $AP, BQ,$ and $CR$ coincide with a side of the triangle $ABC$,
the triple ratio becomes $0$ or $\infty$, because we have a zero either 
in the nominator or in the denominator. In some of the more degenerate configurations, we get zero 
in both nominator and the denominator, so in those configurations the triple ratio is not well-defined.

However, similarly to the cross-ratio, in our construction of Chow quotient with extra lines, 
each triple ratio of six marked points is well-defined. Firstly, in all charts in which 
the triple ratio is well-defined, it is the same (since so it is in generic configuration, 
and the Chow orbit is obtained by a limit of a generic object). Secondly, if for instance 
we need to find a chart in which $\left\{P_1,P_2,P_3; P_4,P_5,P_6\right\}$ is well defined.
By definition, there is a chart in which $\ell_{1,4},\ell_{2,5},\ell_{1,3},\ell_{2,3}$ are placed generically.
We may assume that $\ell_{1,3},\ell_{2,3}$ are two of the sides of a regular triangle, 
and $\ell_{1,4},\ell_{2,5}$ are its altitudes, which will force the points $P_1,P_2,P_3$ 
be the vertexes of that regular triangle, points $P_4'$ and $P_5'$ the midpoints of the sides opposite
to $P_1$ and $P_2$ respectively. The location of $P_6'$ defines the triple ratio; and on this chart 
the triple ratio is well-defined.

\section{The moduli space $\widetilde Y_n$.}

\subsection{The new construction.}
We consider Chow quotient of a different space. 

Our space is a subspace of
 $\left(\PP^2\right)^n\times \left(\PP^{2^*}\right)^{ {n \choose 2} }$ 
of points $P_1, P_2, P_3, ... , P_n$ and lines $\ell_{i,j}$ numbered by 
subsets of $2$ elements in $\left\{1,2,3,...,n\right\}$, with additional condition: 
line $\ell_{i,j}$ contains $P_i$ and $P_j$.

Our group is $GL_3$, acting by projective transformations on the plane (with induced action
on the dual plane). 

The Chow quotient of that space by this group is our new space $\widetilde{Y_n}$.

Our result is 

\textbf{Theorem.} $\widetilde{Y_6}$ is smooth.

Since it is a Chow quotient, $\widetilde{Y_n}$ is a projective algebraic variety. 

Generic configurations of points give just one point in $\widetilde{Y_n}$,
while non-generic may give additional points.
For instance, in the example of singular point of $Y_6$ which is shown in $2.2$ lines would be added 
at all colliding points, some lines will be added at the colliding points,
so that the direction in which the points collided will be taken into account.

\subsection{General remarks on algebraic structure.}
Given a chart, we may choose 5 points and try to compute some cross-ratio they define. 
The lines passing through the fifth point form a projective line. For any 4 points in 
a projective line, the cross-ratio is uniquely defined in $\PP^1$, unless 3 of the lines coincide
(because in that case, we get zero both in the nominator and in the denominator). 
Therefore, for 5 points in the plane, the cross-ratio will not be computable only if
4 of the points are on the same line. For the case of $Y_n$, a problem might also arise 
when two points coincide on a chart. However, for $\widetilde{Y_n}$ the line corresponding 
to a pair of points is always given. 

Cross-ratios ratios are well-defined on $\widetilde{Y_n}$ since 
there is just one chart on which given $4$ points are non-degenerate.
So, to compute $\left[A,B;C,D\right]_E$ we may choose the chart where 
the quadruple $A,B,C,E$ is generic. The lines $AE, BE, CE$ are distinct,
and the line $DE$ is uniquely defined (even if $D$ and $E$ coincides) 
precisely because it is $\widetilde{Y_n}$. So for each cross-ratio,
we get a chart where it is defined by 4 lines, and no more than 2 of those coincide,
so the cross-ratio is defined uniquely on one of the charts.

If the same cross-ratio is defined on 2 different charts, 
then it turns out to be the same on both of them. 
Indeed, it is obvious on a generic configuration 
(corresponding to a generic orbit of group action), 
and the Chow orbits are limiting points of generic orbits.

If for some point of $\widetilde{Y_n}$ we find a chart in which neither 4 points are colinear, 
we may compute all the cross-ratios from the same chart. Unfortunately, it is not always the case, 
but since it covers many cases, we shall discuss it first, specifically for the case of $\widetilde{Y_6}$. 
By choosing 4 marked points in generic position (we may assume those are $P_1,P_2,P_3,P_4$, and considering variation in 
a small neighborhood of the other two points ($P_5$ and $P_6$), we get local coordinates on the moduli space around these points.
Indeed, in the neighborhood of those points, 
all cross-ratios are uniquely defined, and hence any chart containing 
some four points in generic position is uniquely defined.

Those cases can be decomposed into 3 types: one is when neither $P_5$ nor $P_6$ coincide with any other points in the given chart, 
the second is when one of them coincides with one of the first 4 points, and the third case when $P_5$ coincides with $P_6$. 
When more collisions take place, then 4 points are collinear, so those cases have to be considered separately and more carefully.
In the cases when there are no more than 3 collinear points, we may place points $P_1,P_2,P_3,P_4$ as the standard frame 
\[(1,0,0),(0,1,0),(0,0,1),(1,1,1)\]
and in the first subcase, we may choose two pairs of Cartesian coordinates for both $P_5$ and $P_6$.
In the other two subcases, we must choose local coordinates so that in a small neighborhood the changes of the direction of the line 
connecting the two colliding points will be small. 
It is easily achieved (like usually with blow-up) if we take one of 
the coordinates of a point as a cross-ratio of four lines, 
meeting in one of the colliding points, 
and one of the lines connecting the two colliding points.
The second coordinate should be related to a quadruple 
of lines meeting away from the point of collision, with one of the lines passing 
through the colliding point, but in a direction transversal to the line of collision.

In these cases we also get a lot of non-zero cross ratios; 
and we can recover all other charts from this chart, 
including the charts which are not stabilized by points (which will be further discussed 
in the subsection 8.1).

In other, more complicated, cases some case-checking will be required.
We shall organize the case-checking in the following way.
Given a point in $\widetilde{Y_6}$, we may forget point number $6$, by erasing
it from all charts (and disregarding the charts that lose stability). 
Vice-versa, to get a configuration of 6 points, we may list all the types 
of the possible configurations for 5 points, and try inserting $P_6$ 
again into the charts we have for 5 points, and see if more charts are required.
In that way, we may even skip some families of cases.

\subsection{Geometric description of $Y_5 = \widetilde{Y_5}$ and its points.}

We have 3 types of Chow orbits:

(1) The generic configurations, 5-tuples of points, neither two of which coincide 
and neither three are collinear. In terms of cross-ratios, it is the same as to say 
that no cross ration is $0,1$ or $\infty$.

\begin{center}
	\includegraphics[scale=0.5]{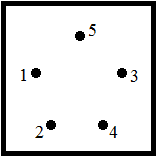}
\end{center}

All other types can be obtained as degeneration of those generic configurations.

(2) Simple degenerations: two points might collide in a generic direction; 
on the other hand, three points may become collinear. 

Indeed, if we deform a generic configuration is such 
a way that points 1 and 2 collide (but in generic direction),
assuming point 1 is the origin of the Cartesian plane and point 2 is at distance $\epsilon$ 
from it, while other points are fixed, and $\epsilon$ tends to zero. We may scale the picture,
so that points $1$ and $2$ are stable and other points approach infinite line, thus becoming 
collinear at the limit.

\begin{center}
	\includegraphics[scale=0.5]{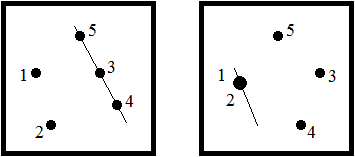}
\end{center}

We get two charts. On one (the right) points 1 and 2 collide; on another (the left) points 3, 4 and 5 are collinear. 
For the case of $\widetilde{Y_5}$ we have to add a line via the two colliding points, however, it is not new information, 
as $\left[P_2,P_3,P_4,P_5\right]_{P_1}$ should be the same on both charts. So, in this case there is no difference between 
the definitions of $\widetilde{Y_5}$ and of $Y_5$.

(3) Deep degeneration: 
in the previous picture, we may move $P_2$ in the left chart so that it comes to a line 
$P_1P_2$. In this case, we get 3 charts, as in the picture:

\begin{center}
	\includegraphics[scale=0.5]{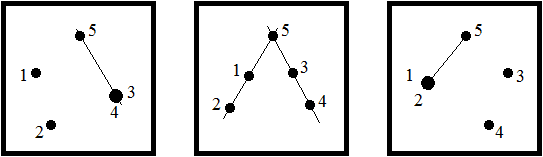}
\end{center}
In thus case also, we see no difference between $\widetilde{Y_5}$ and $Y_5$: indeed, 
in the middle chart all lines are defined by the points, while on the left and on the right,
there are colliding points, but the direction of the line connecting the two colliding points is determined 
uniquely by compatibility of the cross-ration with the middle chart.

It is not possible to degenerate any further, because all charts are stable. So one conclusion is that $Y_5$
is the same as $\widetilde{Y_5}$. 
Indeed, all lines $\ell_{i,j}$ are defined uniquely by the marked points; 
in the charts with coinciding points, the lines are fixed by the cross-ratios 
which can be obtained from another charts by consistency.
Another thing we see from this geometric description is that $Y_5$ is smooth.
Indeed, we see that for each point of the moduli space, there is a chart (the left in the case 2, the middle in the case 3)
which contains 4 points in general position. 
We may assume that those 4 points are $(1,0,0), (0,1,0), (0,0,1)$ and $(1,1,1)$.
The last point may be specified by two ratios of its projective coordinates, 
thus giving a set of two local coordinates in $Y_5$.
Coincidentally, these ratios are cross-ratios, and also the Chow forms used to define the Chow variety. 

\section{Preparations for proof of smoothness of $\widetilde{Y}_6$ - case listing.}

We shall list all the possible Chow orbits, and prove that each is a smooth point. 
In many cases it is obvious; in some an explanation will be required.
So in the end we shall list the hard cases and describe the local coordinates for each case. 

Before we start listing different configurations of points, we should separate two types of stable charts.
Many charts might be stabilized by points, however some of the charts may be stabilized by 
the additional lines (in cases where two many points coincide). 
It would be convenient to list the cases according to the configurations of points and add the lines 
later, but before we we do it we should describe all the possible configurations of 6 points and related 
15 lines which are not stabilized by points only, but also by lines.

We shall called a chart of this kind \textbf{not interesting}, if all information in this chart exists in 
a chart which is stabilized by points only.

\subsection{Case-checking of charts which are stabilized by lines}

We shall classify the charts which are stabilized partly because of lines according to the number 
of distinct marked points. Notice that the configuration of marked points is not stable if and only if
it doesn't contains a generic quadruple, in other words there is a line containing 
either all marked points or all points except one.

(6) If there are six distinct marked points, then all special lines are defined by them are defined uniquely,
so either this configuration is stable without adding special lines, or it is unstable even with special lines.

(5) If there are five distinct marked points, two of the marked points coincide, 
and the line passing through them adds additional information. Assume that $P_5$ and $P_6$ coincide, and 
$\ell_{5,6}$ gives additional information. Without $\ell_{5,6}$ the chart wouldn't be stable, 
so all points except one are collinear. We may present the map is such a way that the line containing
most of the marked points is the infinite line.

If the only marked point outside the infinite line is $P_5$ which coincides with $P_6$, then homotheties 
centered at that point with any coefficient preserve the configuration, so it is not stable.

In the remaining case $P_5$ is at the infinite line, and the configuration might turn out to be stable, 
but not interesting. Indeed, let us say that $P_1$ is not on the infinite line, and we may move it to be the origin. 
Projective transformations preserving the infinite line and the origin are the linear transformations.
Using rotation and homothety we may place the line $\ell_{5,6}$ to be $x=1$ (unless $\ell_{5,6}$ is either 
the infinite line or passes through zero, in which case the configuration is non-stable because of homotheties).
So, in this case the configuration is stable, but all the information is uniquely defined by 
the location of the points $P_2,...,P_6$ on the infinite line. A linear transformation defined 
by a triangular matrix can place $P_2$ at the horizontal infinite point, while $P_5$ and $P_6$ are 
at the vertical infinite point; vertical dilation can place $P_3$ at any other infinite point 
(say $(1:1:0)$) and the location of $P_4$ on the infinite line after that is the only remaining parameter;
however it is governed by some cross-ratio, so the chart is not interesting.

(4) We have four distinct marked points, either two of them are twice marked or one of them is thrice marked, 
and either 3 of them or all four of them are collinear.

(4.1) If all the marked points are collinear, we may place the line containing all of them at infinity. 
The transformation respecting all the marked points are translations and homotheties. 
We have either two lines or three parallel lines which are supposed to produce stability.

If those are 3 parallel lines, there is a family of translations parallel to these lines 
which contradicts stability.

If there are two non-parallel lines, there are homotheties from their intersection points, 
so there is no stability.

If one of these two lines turns out to be the infinite line, there are plenty 
of translations preserving the picture, so it is still not stable.

(4.2) The other case: we have four distinct marked points, 3 of them are on one line, the fourth is outside
that line; two of the points are marked twice and two extra lines through them are added, 
or one of the points is marked thrice and three extra lines through it are added. 
Which leads us to several subcases, and we shall discuss them all, but first we list them:

\begin{itemize}
	\item 4.2.1. Three simple marked points on one line, triple point outside that line.
	\item 4.2.2. Two simple and one triple points on one line, and one more point outside that line.
	\item 4.2.3. Two double and one simple point on one line.
	\item 4.2.4. Two simple marked points and a double point on one line, another double point outside that line.
\end{itemize}

(4.2.1) We assume $P_1,P_2,P_3$ are on the infinite line, $P_4,P_5,P_6$ are at the origin,
lines $\ell_{4,5},\ell_{4,6},\ell_{5,6}$ are added. 
However homotheties centered at the origin with any coefficient preserve the configurations, 
so it is not stable.

(4.2.2) Points $P_1,P_2,P_3$ coincide, and are on one line with $P_4$ and $P_5$, 
we may assume it is the infinite line. We may assume also that $P_1$ is the infinite point in 
the horizontal direction, and $P_4$ is in vertical direction, $~{P_5=(1:1:0)}$, and $P_6$ is at the origin.
The configuration would not be stable without the additional lines because of homotheties from the origin,
however there are 3 additional horizontal lines $\ell_{1,2}, \ell_{1,3}, \ell_{2,3}$.
The only parameters in this picture are the ratios of oriented distances from $P_6$ to those horizontal lines.
Those ratios are actually cross-ratios of the type $\left[P_2,P_3;P_6,P_4\right]_{P_1}$ and similar,
so they are defined by some other charts which are stabilized by the points, so this chart is not interesting.

(4.2.3) This is actually an interesting case, so we shall discuss it and its subcases after we finish the list.
 
(4.2.4) We assume points $P_1$ and $P_2$ coincide, points $P_1,P_2,P_3,P_4$ are on one line, which is the infinite line. 
We may assume that $P_1=P_2=(1:0:0)$, $P_3=(0:1:0)$, $P_4=(1:1:0)$.
Points $P_5$ and $P_6$ are at the origin. 
The lines $\ell_{1,2}$ and $\ell_{5,6}$ are supposed to stabilize the configuration.
If $\ell_{1,2}$ is the infinite line, or is passing through the origin, 
the configuration is preserved by the homotheties. 
So we may assume it is a generic horizontal line; by homothety we may bring it to $y=1$.
The direction of the line $\ell_{5,6}$ is the only remaining parameter in this configuration, 
and it is governed by the cross-ratio $\left[P_1,P_3;P_6,P_4\right]_{P_5}$, so this chart is not interesting.

(3) We have three distinct marked points. The multiplicities are either 4+1+1 or 3+2+1 or 2+2+2.
The points themselves might be either collinear or not, in any case the points don't stabilize 
the configuration without the extra lines. It leads us to 6 cases.

(3.1) Assume the multiplicities are 4+1+1, and the points are not collinear.
We may also assume that $P_1=P_2=P_3=P_4=(1:0:0)$, $P_5=(0:1:0)$, and $P_6=(0:0:1)$.
Then horizontal dilations preserve the configuration, even with the extra horizontal lines.

(3.2) Assume the multiplicities are 4+1+1, and all points are collinear; we may assume 
they are on the infinite line. We may also assume that $P_1=P_2=P_3=P_4=(1:0:0)$.
Then the configuration is preserved by all horizontal translations, so it is non-stable.

(3.3) Assume the multiplicities are 3+2+1, and the points are not collinear.
We may assume $P_1=P_2=P_3=(1:0:0)$, $P_4=P_5=(0:1:0)$, and $P_6=(0:0:1)$.
The extra lines which are added for stabilization are the horizontal lines 
$\ell_{1,2},\ell_{1,3},\ell_{2,3}$ and the vertical $\ell_{4,5}$.
If $\ell_{4,5}$ contains $P_6$ or is the infinite line, the configuration is preserved by horizontal dilations, 
so it is non-stable. Hence by horizontal dilation we may assume $\ell_{4,5}$ is $x=1$.
We also assume that one of the vertical lines, for instance $\ell_{1,2}$ doesn't contain $P_6$ and doesn't 
coincide with the infinite line, so we may assume it is $y=1$ because of horizontal dilations.
The parameters that remain in this chart are the locations of the other two vertical lines, 
which are defined by the cross-ratios $\left[P_6,P_4;P_1,P_3\right]_{P_2}$ and
$\left[P_6,P_4;P_2,P_3\right]_{P_1}$, so all the parameters in this chart are 
defined by the cross-ratios and the chart is not interesting.

(3.4) Assume the multiplicities are 3+2+1, and all points are collinear.
We may assume that $P_1=P_2=P_3=(1:0:0)$ (the horizontal infinite point), 
$P_4=P_5=(0:1:0)$ (the horizontal infinite point), $P_6=(0:0:1)$ (the origin).
The configuration is stabilized because of the 3 horizontal lines $\ell_{1,2},\ell_{1,3},\ell_{2,3}$,
and the vertical line $\ell_{4,5}$. The vertical line might be brought by horizontal dilation
to be at $x=1$, and one of the horizontal lines, say $\ell_{1,2}$ can be brought by a vertical dilation
to be $y_1$. The only parameters in this configuration are the locations of the other two horizontal lines,
they are specified by $\left[P_1,P_4;P_3,P_6\right]_{P_2}$ and $\left[P_2,P_4;P_3,P_6\right]_{P_1}$
so the chart is not interesting.

(3.5) Assume the multiplicities are 2+2+2, and the points are not collinear.
This is an extremely interesting type of configuration, 
so we shall discuss it after we finish the list.

(3.6) Assume the multiplicities are 2+2+2, and all points are collinear.
So the configuration is defined by 4 lines: the line containing all the points, 
and the 3 extra lines. If the four lines are placed in some degenerate way, this configuration 
is non-stable; if they are a generic quadruple, the configuration has no parameters, 
so it is not interesting.

(2) If we have only two distinct marked points, we may assume that one of the points is $(1:0:0)$ and the other is $(0:1:0)$.
The multiplicities are 5+1 or 4+2 or 3+3.

(2.1) If the multiplicities are 5+1, we have several horizontal extra lines, so the configuration
is preserved by horizontal translations. So the configuration is non-stable.

(2.2) If the multiplicities are 4+2, we have several extra lines: some are horizontal and one is vertical.
Then horizontal dilations around that line preserve the configuration, so it is still not stable.

(2.3) The last configuration with two distinct points is that both points are of multiplicity 3.
Assuming \\
\centerline{ $P_1=P_2=P_3=(1:0:0)$ and $P_4=P_5=P_6=(0:1:0)$,}\\
 we have 3 extra horizontal lines 
and 3 extra vertical lines. If at least two of 3 lines in each direction are distinct and not infinite,
The configuration is actually stable, but it is not interesting for a non-obvious reason.  
At least two of the horizontal lines have to be distinct and not infinite, otherwise the chart looses 
stability because of vertical dilations.
At least two of the vertical lines have to be distinct and not infinite, otherwise the chart looses 
stability because of horizontal dilations.
There are two parameters for this type of chart: the ratio at which the three vertical lines intersect 
any horizontal line, and the ratio at which the three horizontal lines intersect any vertical line.
Since they are similar, consider one of the parameters. 
So, on the line $\ell_{4,5}$ there are 4 points: intersections with the horizontal lines 
$\ell_{1,2},\ell_{1,3},\ell_{2,3}$ and the infinite point $P_4$, and we consider the cross-ratio 
of those four points on $\ell_{4,5}$, which is different from the cross-ratios of 4 lines which we usually 
consider. Anyway, similarly to those cross-ratios, it is a projective invariant, 
so it is the same in all the charts of the same Chow orbit in which it is computable.

Consider the chart with points $P_1,P_2,P_3,P_4$ are in generic position. Then the lines
of $\ell_{1,2},\ell_{1,3},\ell_{2,3}$ and $\ell_{4,5}$ are distinct, so the intersections of these lines
are well-defined. Considering those three intersection points with $P_4$, we see no more than 2 coincide, 
so the cross-ratio is uniquely defined in this chart. The cross-ratio is defined by the choice of the line 
$\ell_{4,5}$ through $P_4$, so one of the two parameters in the previous chart is defined by 
the chart which is stabilized by 4 points. 

Similarly the other parameter is defined by a chart in which $P_3,P_4,P_5,P_6$ are placed generically.
So the chart of this type might be omitted.

(1) All configurations is which all points coincide are non-stable even with extra lines, 
because the homotheties centered at that point preserve the configuration.

\subsection{The interesting cases stabilized by lines}
Among the cases in the previous subsection, we have shown that most might be disregarded, 
however the two interesting cases remain: (4.2.3) and (3.5). We shall try to understand something
about those cases, so that later when we organize the case-checking according to the charts which 
are stabilized by points, we shall be able to recognize when adding one of these two types of charts is 
relevant.

We start by considering (4.2.3). We assume $P_1,...,P_5$ are on one line: points $P_1$ and $P_2$ 
coincide, points $P_3$ and $P_4$ also coincide at another point, and point $P_5$ is on the same line;
the point $P_6$ is outside that line.
Two extra lines are added: $\ell_{1,2}$ and $\ell_{3,4}$.
Each of these lines can be generic or degenerate (contain more than 2 marked points). 
For instance if $\ell_{1,2}$ contains either $P_6$ or all the other points, we shall say it is degenerate,
otherwise we say it is generic. 
It is easy to see that if both extra lines are degenerate, then the chart is not stable and might be disregarded.
Even one of the extra lines is degenerate and another is generic, the chart is stable, but it contains
no parameters. If both extra lines are degenerate, the chart is not stable. So it is enough to consider the case in which both extra lines are generic.

We may assume that $P_1P_3P_6$ is a standard equilateral triangle, and the intersection of 
$\ell_{1,2}$ with  $\ell_{3,4}$ is the center of this triangle, then the only parameter in this chart
is the location of $P_5$ on the line $\ell_{1,3}$
In other words, the only parameter defining this chart is the triple ratio 
$\left\{P_1,P_3,P_6; P_2,P_4,P_5\right\}$.

This chart doesn't contain all the information necessary to reconstruct all the other charts, 
but we can reconstruct some.

If we place points $P_1,P_2,P_3,P_4$ generically, we see that $P_6$ is the intersection 
of $P_1P_2$ and $P_3P_4$. Indeed, $\left[P_2,P_3;P_6,P_4\right]_{P_1}=0$, so $P_6$ has to be on $\ell_{1,2}$;
similarly it has to be on $\ell_{1,2}$.

If we place points $P_1,P_2,P_3,P_6$ generically, the points $P_4$ and $P_5$ coincide with $P_3$.
Indeed, for $i \in \left\{ 4,5 \right\}$, we have 
\[\left[P_2,P_3;P_6,P_i\right]_{P_1}=0=\left[P_1,P_3;P_6,P_i\right]_{P_2},\]
so $P_4$ and $P_5$ on this chart have to be on $\ell_{1,3}$ and on $\ell_{2,3}$.

Similarly, on a chart with $P_1,P_3,P_4,P_6$ placed generically, points $P_2$ and $P_5$ coincide with $P_1$.

\begin{center}
\includegraphics{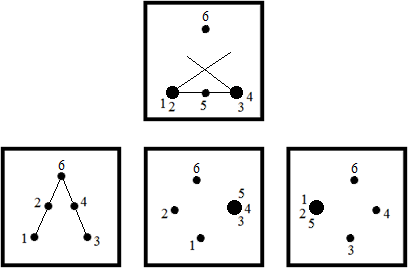}
\end{center}
So, when we shall list configurations of stabilized by points, in order to understand which configurations
are related to this type by a Chow orbit, we shall look for having two configurations with  
triple points, and another chart with two lines of 3 points.

Now we shall discuss the configuration of type (3.5).
We have a triangle of double points $P_1P_2P_3$; the pairs $P_1$ and $P_4$, 
$P_2$ and $P_5$, $P_3$ and $P_6$ collide, and 3 extra lines are added. 
The only parameter is the triple ratio $\left\{P_1,P_2,P_3;P_4,P_5,P_6\right\}$.

\begin{center}
\includegraphics{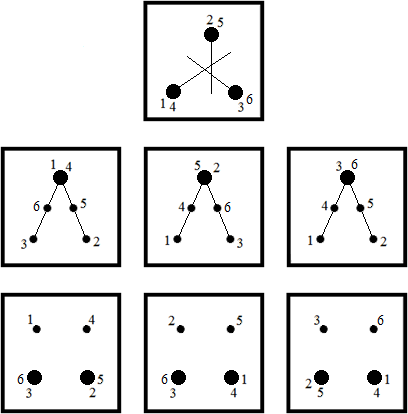}
\end{center}

Consider the chart in which $P_2,P_5,P_3,P_4$ are placed generically. 
Since $\left[P_6,P_2;P_1,P_5\right]_{P_3}=0$ we conclude $P_1$ is on $P_3P_6$;
similarly we see that $P_1$ is on $P_2P_5$. In the same way we see that $P_4$ 
is at the intersection of $P_2P_5$ and $P_3P_6$.

Similarly, we get the other two charts in the second row in the picture.

Now consider the chart in which $P_1,P_2,P_3,P_4$ are placed generically.
Since $\left[P_2,P_3;P_5,P_4\right]_{P_1}=0=\left[P_2,P_3;P_5,P_1\right]_{P_4}$,
the point $P_5$ must be on the lines $P_1P_2$ and $P_4P_2$, so $P_2=P_5$.
Similarly we see that $P_3=P_6$. 

In a similar way, we get two more charts with two double points.

In the type of configuration we considered, there is a degenerate case we should mention,
when the triple ratio $\left\{P_1,P_2,P_3;P_4,P_5,P_6\right\}$ is either zero or infinity.
In this case we get 3 charts instead of one: we may send 2 of the extra lines to the medians 
of the triangle, and the last extra line coincides with one of the sides of the triangle 
(either clockwise or counter-clockwise), but the other conclusions are the same.

To summarize, when we list the configurations according to the charts stabilized by points, 
we might have to add in some cases charts stabilized by lines; one of the cases may occur 
when there are two triple points in some chart, another might happen when we have 
a lot of degenerate charts like in the last picture.

\subsection{Listing the Chow orbits for the proof of smoothness of $\widetilde{Y_6}$.}

First family of cases starts with 5 points placed generically, and point 6 inserted anywhere.
In all those cases there are no 4 points on one line, so we may skip them.

\subsubsection{Adding point 6 to a simple degeneration.}

The second family of cases starts with 5 points in a simple degenerated configuration:

\begin{center}
	\includegraphics{5_degen_simple.png}
\end{center}

If in the left chart $P_6$ is located away from the line $P_3P_4$,
Then no 4 points are on one line, and we may skip these cases. 
If in the right chart $P_6$ is located away from the lines $P_1P_3$, $P_1P_4$, $P_1P_5$,
then in the right chart we have no 4 colinear points, so we may skip these cases as well.

\begin{center}
	\includegraphics{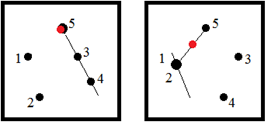}
\end{center}

Since until now $P_3$, $P_4$, and $P_5$ had symmetric roles, we may assume that the point $P_6$
is located in the line $P_1P_5$ in the right chart. 

So in the left chart, points $P_5$ and $P_6$ have to coincide. Indeed, the cross-ratios $\left[P_3,P_5;P_4,P_6\right]_{P_1}$
and $\left[P_3,P_5;P_4,P_6\right]_{P_2}$ are zeroes, so in the left chart also, $P_6$ has to be on the lines $P_1P_5$ and $P_2P_5$.

In the left chart we have two colliding points: $P_5$ and $P_6$. So the direction of $\ell_{5,6}$ has several different possibilities:
\newpage
\begin{itemize}
    \item ($\mathcal{A}$) It might be a generic direction.
    
    		\begin{center}
		\includegraphics{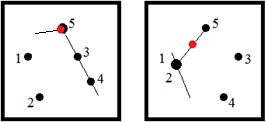}
		\end{center}
    \item ($\mathcal{B}$) It might coincide with the direction of $P_5P_1$ (or $P_5P_2$, but that is symmetric).
    
    		\begin{center}
		\includegraphics{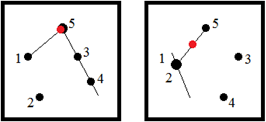}
		\end{center}
    \item ($\mathcal{C}$) it might coincide with the direction of $P_5P_3$.
    
    		\begin{center}
		\includegraphics{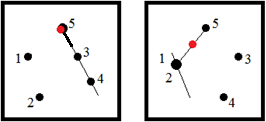}
		\end{center}
\end{itemize}

In the right chart, we also have two colliding points, $P_1$ and $P_2$. 
However, in this case the direction of the line $\ell_{1,2}$ is uniquely defined by $\left[P_2,P_3;P_4,P_5\right]_{P_1}$ which is given in the left chart.

The possible combinatorial ways to place $P_6$ on $P_1P_5$ in the right chart are: 
\begin{itemize}
    \item (1) A generic point on $P_1P_5$.
    \item (2) An intersection point of $P_1P_5$ and $P_3P_4$.
    \item (3) Colliding with $P_1$ (and $P_2$).
    \item (4) Colliding with $P_5$. This case might be separated into subcases according 
		to the direction of $\ell_{5,6}$.
\end{itemize}

We shall consider all those cases and discuss the extra charts appearing in them. 

\subsubsection{$\mathcal{A}.1$}
We start with a first chart having $P_1, P_2,P_3, P_4$ generically places, $P_5$ a generic point on $P_3P_4$, 
and $P_6$ colliding with $P_5$ in a generic direction. 

The second chart has points $P_1,P_3, P_4, P_5$,  placed generically, 
$P_6$ a generic point on $P_3P_5$, and $P_2$ colliding with $P_1$ 
in a generic direction. Notice that its direction can be 
reconstructed from the first chart by $\left[P_2,P_3;P_4,P_5\right]_{P_1}$ 
which should be equal in both charts.

We need to introduce yet another chart in this case, since there should be a chart in which the quadruple $P_1,P_2,P_5,P_6$ is placed generically.
Notice that $\left[P_2,P_5;P_3,P_6\right]_{P_1}$ and $\left[P_2,P_5;P_4,P_6\right]_{P_1}$ are still zeroes as in the first charts, so point $P_3$ and $P_4$
have to be on the line $P_1P_2$, as in the first chart. Also, $\left[P_2,P_5;P_3,P_6\right]_{P_1}$. We may also notice from the first chart that 
$\left[P_1,P_2;P_3,P_6\right]_{P_5}=\left[P_1,P_2;P_4,P_6\right]_{P_5}$, hence the lines $P_3P_5$ and $P_4P_5$ coincide on this new chart, hence $P_3$ and $P_4$
collide. The line $\ell_{3,4}$ on this new chart is uniquely defined on this new chart by $\left[P_1,P_5;P_4,P_6\right]_{P_3}$. 
\begin{center}
	\includegraphics{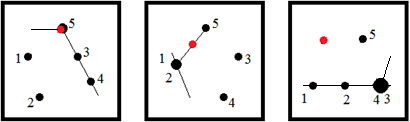}
\end{center}
It is easy to verify that each quadruple of points is generic on precisely one chart.
This case has a lot of nice combinatorial symmetries.

There are no charts with only 4 distinct points, so we don't have to add a chart stabilized by lines.

\subsubsection{$\mathcal{A}.2$}
In the first chart $P_1, P_2, P_3, P_4$ are generically placed, $P_5$ is a generic point on $P_3P_4$, 
and $P_6$ is colliding with $P_5$ in a generic direction. 

In the second chart has points $P_1,P_3, P_4, P_5$ are generically placed, while $P_2$ collides with $P_1$ 
in a generic direction, and $P_6$ is the intersection point on $P_1P_5$ and $P_3P_4$.

It is required to add a chart, on which $P_1,P_2,P_5, P_6$ are placed generically.
Since $\left[P_2,P_5;P_3,P_6\right]_{P_1}=\left[P_2,P_5;P_4,P_6\right]_{P_1}=0$ in the first chart, 
on this this new chart $P_3$ and $P_4$ must be placed on the line $P_1P_2$. 
Since $\left[P_1,P_2;P_3,P_6\right]_{P_5}=\left[P_1,P_2;P_4,P_6\right]_{P_5}$ is given by the first chart, 
the location of points $P_3$ and $P_4$ is uniquely defined as a generic point on $P_1P_2$. Since on this chart points 
$P_3$ and $P_4$ collide, the direction of the line $\ell_{3.4}$ should be discussed also. It follows from the second chart that 
$\left[P_1,P_3;P_5,P_6\right]_{P_4}=0$, so $\ell_{3,4}$ has to pass through $P_6$ on this new chart.
\begin{center}
	\includegraphics{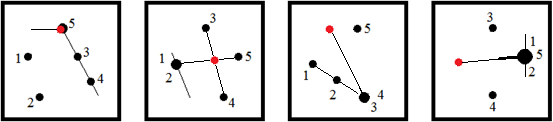}
\end{center}
Yet another chart should be added, in which $P_3,P_4,P_5,P_6$ are placed generically.  From the second chart we may conclude that 
\[\left[P_1,P_4;P_5,P_6\right]_{P_3}=0=\left[P_2,P_4;P_5,P_6\right]_{P_3},\]
 so $P_1$ and $P_2$ have to be on $P_3P_5$.
Similarly, by interchanging 3 and 4 in the last argument, $P_1$ and $P_2$ have to be on $P_4P_5$. So, $P_1$ and $P_2$ collide with $P_5$ on this chart.
Hence, the directions of $\ell_{1,5},\ell_{2,5}$ and $\ell_{1,2}$ have to be specified. It is easy to see from the second chart that 
$\left[P_3,P_5;P_4,P_6\right]_{P_1}=0=\left[P_3,P_5;P_4,P_6\right]_{P_2}$, so in the new chart $\ell_{1,5}$ and $\ell_{2,5}$ pass through $P_6$.
As for line $\ell_{1,2}$ it has a generic direction it may be specified by $\left[P_2,P_3;P_4,P_6\right]_{P_1}$ which is given by both the first and the second chart
in a consistent way.

There is only one chart with only 4 distinct points, so we don't have to add any charts stabilized by lines.

\subsubsection{$\mathcal{A}.3$}
In the first chart $P_1, P_2,P_3, P_4$ are generically placed, $P_5$ is a generic point on $P_3P_4$, 
and $P_6$ is colliding with $P_5$ in a generic direction. 

In the second chart has points $P_1,P_3, P_4, P_5$ are generically placed, while $P_2$ collides with $P_1$ 
in a generic direction, and $P_6$ also collides with $P_1$ and $P_2$. Lines $\ell_{1,6}$ and $\ell_{2,6}$ are passing through $P_5$,
since on he first chart we see that $\left[P_3,P_5;P_4,P_6\right]_{P_1}=0=\left[P_3,P_5;P_4,P_6\right]_{P_2}$.

We have to add more charts. In the third chart charts  $P_1, P_2, P_5, P_6$ are placed generically.
Since $\left[P_3,P_5;P_2,P_6\right]_{P_1}=0=\left[P_4,P_5;P_2,P_6\right]_{P_1}$, points $P_3$ and $P_4$ have to be on $P_1P_2$.
Since $\left[P_3,P_2;P_4,P_1\right]_{P_5}=0$, points $P_3$ and $P_4$ must even coincide. So this chart contains a double point, 
and hence we need to specify the direction of $\ell_{3,4}$, but since $\left[P_4,P_1;P_5,P_6\right]_{P_3}=0$ as can be seen from the second chart, 
so on this new chart $\ell_{3,4}$ must pass through $P_5$.

\begin{center}
	\includegraphics{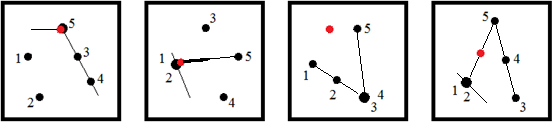}
\end{center}

It still remains to add one more chart in which $P_1, P_3, P_4, P_6$ are placed generically. 
It is easy to see from the first chart that $\left[P_3,P_5;P_4,P_6\right]_{P_1}=0$
and from the second chart that $\left[P_1,P_4;P_6,P_5\right]_{P_3}=0$, therefore $P_5$ must be the intersection of $P_3P_4$ and $P_1P_6$.
It is also easy to see that since $\left[P_1,P_3;P_2,P_4\right]_{P_6}=0=\left[P_1,P_5;P_2,P_6\right]_{P_3}$, 
so $P_1$ has to coincide with $P_2$,
and the direction of $\ell_{1,2}$ is given by $\left[P_2,P_3;P_4,P_5\right]_{P_1}$ which is given by the first chart.

Now, all quadruples appear in generic position in precisely one chart.
However, if we switch the names of $P_5$ and $P_6$, we get the previous case $\mathcal{A}.2$, so we may skip this case.

\subsubsection{$\mathcal{A}.4$}
In the first chart $P_1, P_2,P_3, P_4$ are generically placed, $P_5$ is a generic point on $P_3P_4$, 
and $P_6$ is colliding with $P_5$ in a generic direction. 

In the second chart has points $P_1,P_3, P_4, P_5$ are generically placed, while $P_2$ collides with $P_1$ 
in a generic direction, and $P_6$ collides with $P_5$. The direction of $\ell_{5.6}$ might be of different types:

\begin{itemize}
    \item $P_1P_5$,
    \item $P_3P_5$,
    \item $P_4P_5$, which is symmetric to the previous,
    \item generic, i. e. non of the above.
\end{itemize}
However, only the first of those types consider further discussion, since in the other 3 in the second chart we may compute any cross-ratio, 
(so this is similar to the case when no 4 points on some chart are collinear, even if there are 4 collinear points.

So, that is the case we will further refer to as $\mathcal{A}.4$. It turns out we have to add two more charts. 
The third chart is almost the same as in the previous case $\mathcal{A}.3$ (and for the same reasons), except the direction of $\ell_{3,4}$. 
Since $\left[P_1,P_5;P_4,P_6\right]_{P_3}=0$, the line $\ell_{3,4}$ has to contain $P_1$.

There is only one chart with only 4 distinct points, so we don't have to add any charts stabilized by lines.

\begin{center}
	\includegraphics{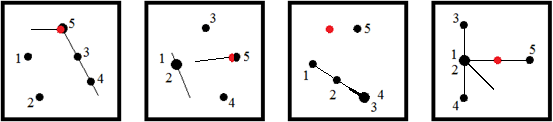}
\end{center}
The only quadruple which doesn't appear as generic in the first 4 charts is $P_3,P_4, P_5, P_6$.
Therefore, an extra chart is required. It is easy to see that $P_1$ and $P_2$ are at the intersection of $P_3P_4$ and $P_5P_6$. 
That might be seen from 
\[ \left[P_6,P_3;P_1,P_4\right]_{P_5}=\left[P_1,P_5;P_4,P_6\right]_{P_3}=\]
\[=\left[P_6,P_3;P_2,P_4\right]_{P_5}=\left[P_2,P_5;P_4,P_6\right]_{P_3}=0 \]
It gives rise to a question, what is $\ell_{1,2}$, and in fact it might be anything.
A good way to realize it is to think of lines from $P_1$ to the other marked points as a point in $X_5$, and different charts of 
$\widetilde{Y_6}$ induce different charts of $X_5$. When one considers the lines from $P_1$ as different charts of $X_5$, 
one might see that neither choice of $\ell_{1,2}$ in the last chart leads to a contradiction. 

\subsubsection{$\mathcal{B}.1$}
In the first chart, $P_1,P_2,P_3,P_4$ are placed generically, $P_5$ is on the line $P_3P_4$, $P_6$ coincides with $P_5$,
$\ell_{5,6}$ contains $P_1$. In the second chart, $P_1,P_3,P_4,P_5$ are placed generically, 
points $P_1$ and $P_2$ coincide, $P_6$ is a generic point on $P_1P_5$.

We have to add a third chart, on which $P_2,P_3,P_5,P_6$ are placed generically.
As $\left[P_1,P_5;P_2,P_6\right]_{P_3}=0=\left[P_1,P_3;P_6,P_2\right]_{P_5}$ 
from the first chart, so in the new chart $P_1$ 
has to be at the intersection of $P_2P_3$ and $P_5P_6$. 
Since $\left[P_1,P_3;P_2,P_4\right]_{P_6}=0=\left[P_5,P_3;P_6,P_4\right]_{P_2}$ 
as might be seen from the first chart, $P_4$ must coincide with $P_3$ on this chart. 
Finally, the direction of $\ell_{3,4}$ 
is defined by $\left[P_1,P_4;P_5,P_6\right]_{P_3}$ to be generic 
(distinct from $0,1,\infty$) from the second chart.

We have to add yet another chart in which $P_1,P_2,P_5,P_6$ are placed generically.
Notice that $\left[P_1,P_2;P_6,P_3\right]_{P_5}=0=\left[P_1,P_2;P_5,P_3\right]_{P_6}$. 
So on this new chart $P_3$ coincides with $P_2$.
For a similar reason (replace 3 with 4), also $P_4$ coincides with $P_2$. 
Since the last chart has a triple point, we should also discuss the location of 3 lines there:
$\ell_{2,3},\ell_{2,4}$, and $\ell_{3,4}$.
Since $\left[P_1,P_5;P_2,P_6\right]_{P_3}=0=\left[P_1,P_5;P_2,P_6\right]_{P_4}$, the lines $\ell_{2,3},\ell_{2,4}$ have to coincide with $P_1P_2$.
The direction $\ell_{3,4}$ is specified by $\left[P_1,P_4;P_5,P_6\right]_{P_3}$, which is given by the second chart.
\begin{center}
	\includegraphics{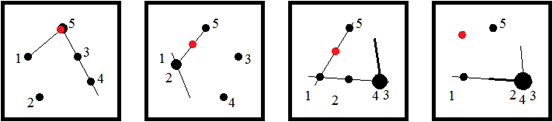}
\end{center}

Actually, this case is the same as $\mathcal{A}.2$, if we rename $P_1,P_2,...,P_6$ 
by $~{P_6,P_5,P_1,P_2,P_3,P_4}$ respectively, so we may skip it in further discussion.

\subsubsection{$\mathcal{B}.2$}
In the first chart, $P_1,P_2,P_3,P_4$ are placed generically, $P_5$ is on the line $P_3P_4$, $P_6$ coincides with $P_5$,
$\ell_{5,6}$ contains $P_1$. In the second chart, $P_1,P_3,P_4,P_5$ are placed generically, 
points $P_1$ and $P_2$ coincide, $P_6$ is the intersection of $P_1P_5$ and $P_3P_4$.

It is required to add a third chart, on which $P_2,P_3,P_5,P_6$ are placed generically.
Since $\left[P_1,P_5;P_3,P_6\right]_{P_2}=0=\left[P_4,P_5;P_3,P_6\right]_{P_2}$ in the first chart, 
the points $P_1$ and $P_4$ must be on a line $P_2P_3$ in the new chart.
Since $\left[P_1,P_2;P_6,P_3\right]_{P_5}=0$ in the first chart, $P_1$ must also be on $P_5P_6$ in the new chart,
so it is at the intersection of $P_5P_6$ and $P_2P_3$.
Since from the second chart $\left[P_3,P_2;P_4,P_6\right]_{P_5}=0$, the lines $\ell_{3,5}$ and $\ell_{4,5}$ coincide, so $P_4$ collides with $P_1$ in the third chart.
As for the direction of $\ell_{3,4}$, from the second chart we see that $\left[P_1,P_3;P_5,P_6\right]_{P_4}=0$, 
hence $\ell_{3,4}$ contains $P_6$ in the third chart.

A fourth chart is required should have $P_1, P_3, P_5, P_6$ placed generically. 
It is easy to see even by comparing with the third chart that points $P_2,P_3,P_4$ coincide;
because for $i\in \{2,4\}$ and  $\{k,m\}=\{5,6\}$ we have $\left[P_1,P_3;P_k,P_i\right]_{P_m}=0$.
It is also easy to see that $\ell_{2,3}$ and $\ell_{2,4}$ coincide with $P_3P_1$ since for $i\in \{3,4\}$ we have $\left[P_1,P_5;P_2,P_6\right]_{P_i}=0$.
Since $\left[P_1,P_4;P_5,P_6\right]_{P_3}=0$, the line $\ell_{3,4}$ has to contain $P_6$.

\begin{center}
	\includegraphics{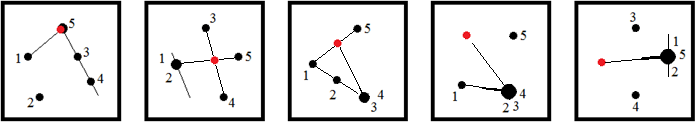}
\end{center}

A fifth should be added, in which $P_3,P_4,P_5,P_6$ are placed generically. It is precisely the same as the last chart in the case $\mathcal{A}.2$ 
and for precisely the same reasons, so we won't repeat them.

Among these charts, we have two triple points. So we might have to add a chart stabilized by lines
of type (4.2.3): with two double points $P_3=P_4$ and $P_1=P_5$, and $P_2$ on the same line and $P_6$
outside that line. Such a chart is defined by a triple ratio $\left\{P_1,P_3,P_6; P_4,P_5,P_2\right\}$.
This triple ratio is undefined in the first four charts, but it is defined by the fifth chart 
(as zero), so we may skip that chart.

\subsubsection{$\mathcal{B}.3$}
In the first chart, $P_1,P_2,P_3,P_4$ are placed generically, $P_5$ is on the line $P_3P_4$, $P_6$ coincides with $P_5$,
$\ell_{5,6}$ contains $P_1$. In the second chart, $P_1,P_3,P_4,P_5$ are placed generically, 
points $P_1$ and $P_2$ and $P_6$ coincide. The direction of the lines at the triple point in the second chart are defined by 
the first chart uniquely. 
The direction of $\ell_{1,2}$ is generic and is defined by $\left[P_2,P_3;P_4,P_5\right]_{P_1}$
which is not zero.
The direction of $\ell_{1,5}$ and $\ell_{2,5}$ are chosen to contain $P_5$ since 
$\left[P_3,P_5;P_4,P_6\right]_{P_1}=0=\left[P_3,P_5;P_4,P_6\right]_{P_2}$.

We have to add a third chart with $P_2,P_3,P_5, P_6$ placed generically.
Since $\left[P_3,P_5;P_1,P_6\right]_{P_2}=0=\left[P_2,P_1;P_3,P_6\right]_{P_5}$, 
the point $P_1$ has to be at the intersection of lines $P_5P_6$ and $P_2P_3$ in the third chart.

Since $\left[P_3,P_5;P_4,P_6\right]_{P_2}=0=\left[P_2,P_3;P_5,P_4\right]_{P_6}$, 
the point $P_4$ has to coincide with $P_3$. Thus we have a double point, so the direction of 
$\ell_{3,4}$ has to be discussed, but as $\left[P_4,P_1;P_5,P_6\right]_{P_3}=0$ in the second chart,
$\ell_{3,4}$ has to contain $P_5$.

We have to add a fourth chart, with $P_1,P_3,P_4, P_5$ placed generically. 
It is easy to see from the third chart that $P_5$ must be the intersection point of $P_3P_4$ and $P_1P_6$, and that $P_1$ and $P_2$ collide.
The direction of $\ell_{1,2}$ is obtained from $\left[P_2,P_3;P_4,P_5\right]_{P_1}=0$, which is given by chart 1.

We have to add a fifth chart with $P_1,P_2,P_5,P_6$ placed generically. 
For $i\in \{3,4\}$ and $\{k,m\}=\{5,6\}$ we have $\left[P_1,P_2;P_m,P_i\right]_{P_k}=0$, so $P_2$ coincides with $P_3$ and $P_4$ in the fifth chart.
The line $\ell_{3,4}$ has to contain $P_5$ as $\left[P_1,P_4;P_6,P_5\right]_{P_3}=0$.
The lines $\ell_{2,3}$ and $\ell_{2,4}$ have to contain $P_1$, as $\left[P_1,P_5;P_i,P_6\right]_{P_2}=0$.

\begin{center}
	\includegraphics{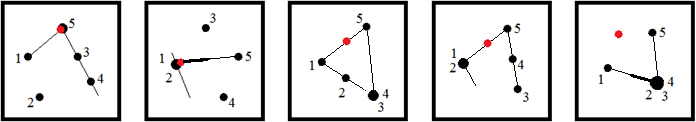}
\end{center}
However, if we swap the names of $P_5$ and $P_6$, we get the previous case $\mathcal{B}.2$ so this case may be omitted.

\subsubsection{$\mathcal{B}.4$}
In the first chart, $P_1,P_2,P_3,P_4$ are placed generically, 
$P_5$ is on the line $P_3P_4$, $P_6$ coincides with $P_5$,
$\ell_{5,6}$ contains $P_1$. In the second chart, $P_1,P_3,P_4,P_5$ are placed generically, 
points $P_1$ and $P_2$ coincide, also points $P_5$ and $P_6$ coincide.
In the second chart, there are two double points, so the directions of $\ell_{1,2}$ and of $\ell_{5,6}$ must be discussed.
The direction of $\ell_{1,2}$ is generic and is specified by $\left[P_2,P_3;P_4,P_5\right]_{P_1}$ which is given by the first chart.
The direction of $\ell_{5,6}$ can be anything. However if it is anything else then $P_1P_5$, all cross-ratios are well-defined in the second chart, 
and we may finish the discussion in the same way as we did for the case of charts where there are no more than 3 collinear points.
So, we shall consider only the case in which $\ell_{5,6}$ contains $P_1$ in the second chart.

The third chart has $P_1,P_2,P_5,P_6$ placed generically. For $i\in \{3,4\}$ we have $\left[P_i,P_5;P_2,P_6\right]_{P_1}=0$,
hence points $3,4$ are on the line $P_1P_2$ in the third chart. 
Notice that for $i=3,4$ we have $\left[P_i,P_5;P_2,P_6\right]_{P_1}=0$ and also $\left[P_i,P_1;P_2,P_6\right]_{P_5}=0$. 
Hence $P_3$ and $P_4$ coincide with $P_2$. The lines defined by pairs in this triple points have to contain $P_1$,
since $\left[P_1,P_5;P_4,P_6\right]_{P_3}=0$ and also $\left[P_1,P_5;P_i,P_6\right]_{P_2}=0$ for $i\in \{3,4\}$.

Consider the fourth chart in which $P_3,P_4,P_5,P_6$ are placed generically.
Since $\left[P_3,P_i;P_4,P_5\right]_{P_6}=0=\left[P_i,P_5;P_4,P_6\right]_{P_3}$ for $i\in\{1,2\}$, the points $P_1$ and $P_2$ must be the intersection 
point of $P_3P_4$ and $P_5P_6$. Hence $\ell_{1,2}$ on this chart has to be specified. However, $\left[P_3,P_2;P_4,P_5\right]_{P_1} \neq 0$, 
and the only way for it to be non-zero in the fourth chart if it is undefined on this chart, meaning $\ell_{1,2}$ contains $P_3$ and $P_4$.

We still need a fifth chart, on which $P_2,P_3,P_5,P_6$ are placed generically. 
Since $\left[P_1,P_5;P_3,P_6\right]_{P_2}=0=\left[P_2,P_1;P_3,P_5\right]_{P_6}$, the point $P_1$ 
is the intersection of $P_2P_3$ and $P_5P_6$. 
Since $\left[P_3,P_2;P_4,P_5\right]_{P_6}=0=\left[P_3,P_2;P_4,P_6\right]_{P_5}$, the points $P_3$ and $P_4$ coincide. 
The line $\ell_{3,4}$ has to contain $P_2$ since $\left[P_2,P_5;P_4,P_6\right]_{P_3}=0$.

\begin{center}
	\includegraphics{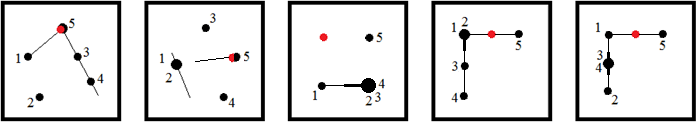}
\end{center}

Here we have 2 charts with only four distinct points, but they are of different types 
(one with the two double points, another with a triple point), so charts stabilized by 
lines don't have to be added.

\subsubsection{$\mathcal{C}.1$}
In the first chart, $P_1,P_2,P_3,P_4$ are placed generically, $P_5$ is on the line $P_3P_4$, $P_6$ coincides with $P_5$,
$\ell_{5,6}$ contains $P_3$. In the second chart, $P_1,P_3,P_4,P_5$ are placed generically, 
points $P_1$ and $P_2$ coincide, $P_6$ is a generic point on $P_1P_5$.

We might discuss two more charts which are required (and it turns out it is the same as $\mathcal{A}.4$), 
but anyway, the second chart is equivalent to the starting point of cases $\mathcal{A}$ up to renaming of the marked points, 
therefore we may skip this sub-case.


\begin{center}
	\includegraphics{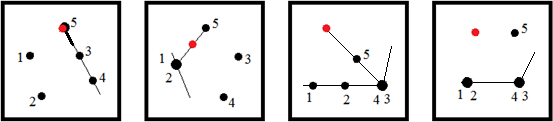}
\end{center}

\subsubsection{$\mathcal{C}.2$}
In the first chart, $P_1,P_2,P_3,P_4$ are placed generically, $P_5$ is on the line $P_3P_4$, $P_6$ coincides with $P_5$,
$\ell_{5,6}$ contains $P_3$. In the second chart, $P_1,P_3,P_4,P_5$ are placed generically, 
points $P_1$ and $P_2$ coincide, $P_6$ is the intersection of $P_1P_5$ and $P_3P_4$.

The third chart must be added in which $P_1,P_2, P_5,P_6$ are placed generically. Notice that for $i\in \{3,4\}$,
$\left[P_2,P_5;P_i,P_6\right]_{P_1}=0=\left[P_6,P_1;P_i,P_2\right]_{P_5}$. Hence $P_3$ and $P_4$ are at the intersection of $P_1P_2$ 
and $P_5P_6$ in the third chart. The line $\ell_{3,4}$ contains $P_5$, 
since $\left[P_4,P_1;P_6,P_5\right]_{P_3}=0$ in the second chart (otherwise it would by $\infty$ on this chart,
this way $\left[P_4,P_1;P_6,P_5\right]_{P_3}$ it is undefined on this chart which is the only possibility).

The fourth chart has points $P_1,P_3,P_5,P_6$ placed generically. 
Notice that $\left[P_1,P_3;P_2,P_6\right]_{P_5}=0=\left[P_1,P_3;P_2,P_5\right]_{P_6}$, so points $P_1$ and $P_2$ coincide. 
Also $\left[P_1,P_3;P_6,P_4\right]_{P_5}=0=\left[P_1,P_3;P_5,P_4\right]_{P_6}$, so points $P_3$ and $P_4$ coincide. 
The line $\ell_{1,2}$ has to contain $P_3$, since $\left[P_1,P_5;P_3,P_6\right]_{P_2}=0$.
The line $\ell_{3,4}$ has to contain $P_6$, since $\left[P_1,P_4;P_5,P_6\right]_{P_3}=0$.

The fifth chart has to be added to have $P_3, P_4, P_5, P_6$ placed generically. 
For $i\in\{1,2\}$ we have $\left[P_i,P_6;P_5,P_4\right]_{P_3}=0=\left[P_i,P_6;P_5,P_3\right]_{P_4}$, 
hence $P_2$ and $P_1$ coincide with $P_5$. The lines $\ell_{5,i}$ for $i\in\{1,2\}$ contain $P_5$, 
since $\left[P_i,P_3;P_6,P_4\right]_{P_5}=0$. 
The direction of $\ell_{1,2}$ is generic, consistent with the first chart so that $\left[P_2,P_3;P_4,P_5\right]_{P_1}$ would be the same.
\begin{center}
	\includegraphics{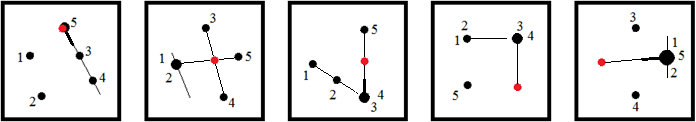}
\end{center}

There are only two charts with only 4 distinct points, but they are of different types 
(one is with two double points, another with a triple point), 
so we don't have to add any charts stabilized by lines.

\subsubsection{$\mathcal{C}.3$}
In the first chart, $P_1,P_2,P_3,P_4$ are placed generically, $P_5$ is on the line $P_3P_4$, $P_6$ coincides with $P_5$,
$\ell_{5,6}$ contains $P_3$. In the second chart, $P_1,P_3,P_4,P_5$ are placed generically, 
points $P_1$ and $P_2$ and $P_6$ coincide.
The direction of the lines at the triple point in the second chart are defined by 
the first chart uniquely. The direction of $\ell_{1,2}$ is generic and is defined by $\left[P_2,P_3;P_4,P_5\right]_{P_1}$.
The lines $\ell_{1,5}$ and $\ell_{2,5}$ must contain $P_5$ since 
\[\left[P_3,P_5;P_4,P_6\right]_{P_1}=0=\left[P_3,P_5;P_4,P_6\right]_{P_2}.\]

A third chart must be added in which $P_1,P_2, P_5,P_6$ are placed generically. Notice that for $i\in \{3,4\}$
$\left[P_2,P_5;P_i,P_6\right]_{P_1}=0=\left[P_6,P_1;P_i,P_2\right]_{P_5}$. Hence $P_3$ and $P_4$ are at the intersection of $P_1P_2$ 
and $P_5P_6$ in the third chart. The line $\ell_{3,4}$ contains $P_5$ 
since $\left[P_4,P_1;P_5,P_6\right]_{P_3}=0$.

In the fourth chart points $P_1,P_3,P_4,P_6$ are placed generically. 
Notice that $\left[P_1,P_3;P_2,P_6\right]_{P_4}=0=\left[P_1,P_3;P_2,P_5\right]_{P_6}$, so points $P_1$ and $P_2$ coincide. 
Also, $\left[P_1,P_4;P_6,P_5\right]_{P_3}=0=\left[P_5,P_3;P_6,P_4\right]_{P_1}$, so $P_5$ should be the intersection point of $P_3P_4$ and $P_1P_6$. 
The line $\ell_{1,2}$ is defined $\left[P_2,P_3;P_4,P_5\right]_{P_1}$ to be the same as in the first chart.

In the fifth chart points $P_1, P_3, P_5, P_6$ are placed generically. Points $P_1$ and $P_2$ must coincide 
because $\left[P_1,P_3;P_2,P_6\right]_{P_5}=0=\left[P_1,P_3;P_2,P_5\right]_{P_6}$. 
Points $P_3$ and $P_4$ must coincide 
because $\left[P_1,P_3;P_6,P_4\right]_{P_5}=0=\left[P_1,P_3;P_5,P_4\right]_{P_6}$. 
The lines through the double points: line $\ell_{1,2}$ must contain $P_3$ because 
$\left[P_2,P_5;P_3,P_6\right]_{P_1}$; line $\ell_{3,4}$ must contain $P_5$ because 
$\left[P_1,P_3;P_6,P_5\right]_{P_4}$.
\begin{center}
	\includegraphics{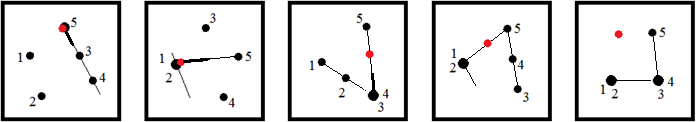}
\end{center}

However, if we swap $P_5$ and $P_6$, we get the previous case $\mathcal{C}.2$, so this case may be omitted.

\subsubsection{$\mathcal{C}.4$}
In the first chart, $P_1,P_2,P_3,P_4$ are placed generically, $P_5$ is on the line $P_3P_4$, $P_6$ coincides with $P_5$,
$\ell_{5,6}$ contains $P_3$. In the second chart, $P_1,P_3,P_4,P_5$ are placed generically, 
points $P_1$ and $P_2$ coincide, also points $P_5$ and $P_6$ coincide.

As in the cases $\mathcal{A}.4$ and $\mathcal{B}.4$, the direction of $\ell_{5,6}$ 
in the second chart might be anything, and as in those cases, 
the only case for which smoothness is not obvious is when $\ell_{5,6}$ contains $P_1$.

In the third chart, points $P_1,P_2,P_5,P_6$ are placed generically. For $i\in \{3,4\}$
we have 
$\left[P_2,P_5;P_i,P_6\right]_{P_1}=0=\left[P_i,P_1;P_6,P_2\right]_{P_5}$, 
hence $P_3$ and $P_4$ must coincide with the intersection point of $P_1P_2$ and $P_5P_6$. 
The direction of $\ell_{3,4}$ might be arbitrary.

In the fourth chart, $P_1,P_3,P_5,P_6$ are placed generically. Notice that for $\{s,t\}=\{5,6\}$, we have 
$\left[P_t,P_1;P_3,P_2\right]_{P_s}=0$, hence $P_1$ and $P_2$ coincide, also 
$\left[P_t,P_3;P_1,P_4\right]_{P_s}=0$, hence $P_3$ and $P_4$ coincide. 
The lines which are related to the double points $\ell_{1,2}$ and $\ell_{3,4}$ have to coincide with $P_1P_3$ because
$\left[P_1,P_5;P_4,P_6\right]_{P_3}=0=\left[P_2,P_5;P_3,P_6\right]_{P_1}$.

In the fifth chart, points $P_3,P_4,P_5,P_6$ are placed generically. For any $i\in \{1, 2\}$
we have $\left[P_i,P_5;P_4,P_6\right]_{P_3}=0=\left[P_i,P_3;P_6,P_4\right]_{P_5}$, 
hence $P_1$ and $P_2$ are at the intersection point of $P_3P_4$ and $P_5P_6$.
The line $\ell_{1,2}$ coincides with $P_3P_4$, 
otherwise we would have $\left[P_2,P_3;P_5,P_4\right]_{P_1}=0$, which is not the case.

\begin{center}
	\includegraphics{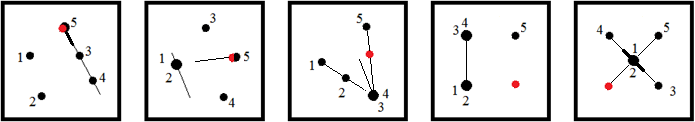}
\end{center}

Regarding possible extra charts stabilized by points:
there are two charts with only 4 distinct points. 
Neither of them has a triple point, so the type (4.2.3) is irrelevant. 
The type (3.5) requires 3 charts with two double points each, and we have only two such charts.
So we don't have to add any charts stabilized by lines.

\subsection{Adding point 6 to a deep degeneration.}
Now we take the most degenerate configuration of 5 points, and attempt to add a sixth point.
\begin{center}
	\includegraphics[scale=0.5]{5_degen_deep.png}
\end{center}
On any chart, if we have no 4 points on one line, we know that the neighborhood of this Chow orbit is smooth. 
Therefore, a new point must belong to the line $P_1P_2$ or to the line $P_3P_4$ in the middle chart. Let us assume it is the line $P_3P_4$,
since it is symmetric. There are several options for that:

\begin{itemize}
	\item ($\mathcal{D}$) $P_6$ might be a generic point on the line $P_3P_4$.
	\item ($\mathcal{E}$) $P_6$ might coincide with $P_3$ (or $P_4$, but it is symmetric). 
	\item ($\mathcal{F}$) $P_6$ might coincide with $P_5$ in the middle chart.
\end{itemize}
We shall discuss all these cases separately. 

\subsubsection{$\mathcal{D}$}
We start with the case when $P_6$ is a generic point on the line $P_3P_4$ in the seconds chart in the previous picture. 
In the third chart, we are reduced to the case that $P_6$ is on 
a generic line through the double point $P_1=P_2$, 
since $\left[P_2,P_3;P_6,P_4\right]_{P_1}$ should be the same as in the second chart. 
In the third chart, the point $P_6$ might be anywhere on that line; however, unless it coincides with $P_1$ and $P_2$ 
there are no more than 3 point on one line on this chart, the neighborhood of that configuration is smooth.
Therefore, we have to discuss only the case when in the third chart, point $P_6$ coincides with $P_1$ and $P_2$.

In the first chart points $P_3$ and $P_4$ have to coincide with $P_6$, because $\left[P_3,P_2;P_6,P_5\right]_{P_1}=0=\left[P_3,P_1;P_6,P_5\right]_{P_2}$.
All the lines through the triple point in the first chart: $\ell_{3,4},\ell_{3,6}$ and $\ell_{4,6}$ contain $P_5$, because for $\{i,j\} \subset \{3,4,6\}$
we have $\left[P_i,P_1;P_5,P_2\right]_{P_j}=0$.

We also have a triple point $P_1=P_2=P_6$ in the third chart, so three lines have to be discussed. As $\left[P_2,P_3;P_5,P_4\right]_{P_1}=0$, so 
$\ell_{1,2}$ contains $P_5$. The lines $\ell_{1,6}$ and $\ell_{2,6}$ coincide, and are specified 
by $\left[P_3,P_4;P_5,P_6\right]_{P_1}=\left[P_3,P_4;P_5,P_6\right]_{P_2}$.

A fourth chart should be added in which $P_1,P_3,P_4,P_6$ are placed generically.
The point $P_5$ must be on $P_3P_4$, because $\left[P_1,P_4;P_6,P_5\right]_{P_3}=0$. Within $P_3P_4$, the point $P_5$ is specified by $\left[P_3,P_4;P_5,P_6\right]_{P_1}$.
The line $\ell_{1,2}$ contains $P_5$, because $\left[P_2,P_3;P_5,P_4\right]_{P_1}=0$. 

\begin{center}
	\includegraphics{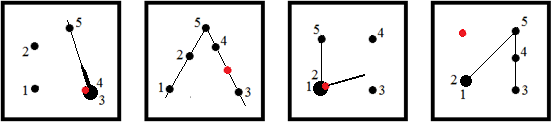}
\end{center}

As in all cases, one must consider if there should be extra charts stabilized by lines.
In this case, we have two charts with triple points, so it looks similar to (4.2.3).
The triples forming the triple points are $P_1,P_2,P_6$ and $P_3,P_4,P_6$, so we might expect 
a chart with a double point $P_1=P_2$, another double point $P_3=P_4$, point $P_6$ 
on the line connecting these two double points, and point $P_5$ outside the line.
This extra chart is defined by the triple ratio $\left\{ P_1,P_3,P_5; P_2,P_4,P_6 \right\}$.
This triple ratio is undefined in the first, the second, and the fourth charts, however 
it is defined to be zero by the third chart, so this chart (even though it exists) 
is redundant. 

\subsubsection{$\mathcal{E}$}
We start with 3 charts: in the second chart $P_1,P_2$ and $P_5$ are collinear, also $P_3,P_4$ and $P_5$ are collinear.
In the first chart $P_1,P_2, P_3$ and $P_5$ are placed generically, $P_4$ coincides with $P_3$, $\ell_{3,4}$ contains $P_5$. 
In the third chart $P_1,P_3, P_4$ and $P_5$ are placed generically, $P_1$ coincides with $P_2$, $\ell_{1,2}$ contains $P_5$. 

We add point 6, which coincides with point 3 in the second chart. In the first chart, point $P_3$ and $P_6$ coincide as well, 
because $\left[P_3,P_2;P_6,P_5\right]_{P_1}=0=\left[P_3,P_1;P_6,P_5\right]_{P_2}$. The line $\ell_{4,6}$ in the first chart contains $P_5$,
because $\left[P_1,P_6;P_2,P_5\right]_{P_4}=0$. The lines $\ell_{3,6}$ in the first and the second charts should be chosen consistently in the first 
and the second charts, and we shall do the case separation a bit later.

In the third chart, point $P_6$ must be added on the line $P_1P_3$, 
because $\left[P_3,P_4;P_6,P_5\right]_{P_1}=0$.

\begin{center}
	\includegraphics{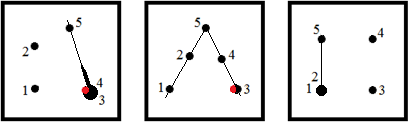}
\end{center}

Hence $P_6$ on the third chart might be one of the following:
\begin{itemize}
	\item ($a$) A generic point of $P_1P_3$.
	\item ($b$) Coincide with $P_1$ and $P_2$. 
	\item ($c$) Coincide with $P_3$.
	\item ($d$) The intersection of $\ell_{1,3}$ and $\ell_{4,5}$.
\end{itemize}

The case ($a$) may be omitted, since it is the same as adding $P_2$ to a configuration 
of the other 5 points with simple degeneration, which is what we have 
listed in the sub-cases of $\mathcal{A,B,C}$.

So we must look at the other 3 cases. Since the described charts don't contain all the information, 
those cases will be further divided into sub-cases.

\subsubsection{$\mathcal{E}.b$}

In the third chart the lines $\ell_{1,6}$ and $\ell_{2,6}$ must contain $P_3$ because
\[\left[P_3,P_4;P_6,P_5\right]_{P_1}=0=\left[P_3,P_4;P_6,P_5\right]_{P_2}\]

\begin{center}
	\includegraphics{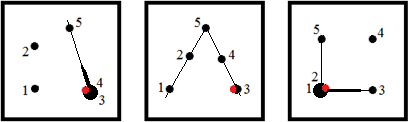}
\end{center}

Now regarding the direction of $\ell_{3,6}$ in the second chart (which is consistent with the first chart), it might:
\begin{itemize}
	\item ($\mathcal{E}.1$) Be generic.
	\item ($\mathcal{E}.2$) Contain $P_1$ (or $P_2$, but it is similar). 
	\item ($\mathcal{E}.3$) Contain $P_4$ and $P_5$.
\end{itemize}
We shall discuss these cases separately. 

\subsubsection{$\mathcal{E}.1$}
We start with the sub-case that $\ell_{3,6}$ has generic direction.

A fourth chart is added in which points $P_1,P_3,P_4, P_6$ are placed generically.
Since $\left[P_1,P_4;P_6,P_5\right]_{P_3}=0=\left[P_1,P_4;P_3,P_5\right]_{P_6}$, points $P_4$ and $P_5$ coincide. 
Since $\left[P_1,P_3;P_6,P_5\right]_{P_4}=0$, the line $\ell_{4,5}$ contains $P_3$.
The point $P_2$ has to be on the line $P_1P_4$, because $\left[P_2,P_3;P_4,P_5\right]_{P_1}=0$.
The specific location of $P_2$ on that line is defined by $\left[P_1,P_2;P_3,P_5\right]_{P_6}$; is distinct from $P_1$ and $P_4$.

A fifth chart has to be added, so that $P_1,P_4,P_5,P_6$ are placed generically.
Since $\left[P_1,P_5;P_2,P_6\right]_{P_4}=0=\left[P_1,P_4;P_2,P_6\right]_{P_5}$, the points $P_1$ and $P_2$ coincide.
Since $\left[P_3,P_4;P_6,P_5\right]_{P_1}=0=\left[P_3,P_1;P_4,P_6\right]_{P_5}$, the point $P_3$ must be on $P_1P_6$ and $P_5P_4$.
A line $\ell_{1,2}$ must contain $P_5$, because $\left[P_4,P_2;P_6,P_5\right]_{P_1}=0$.
\begin{center}
	\includegraphics{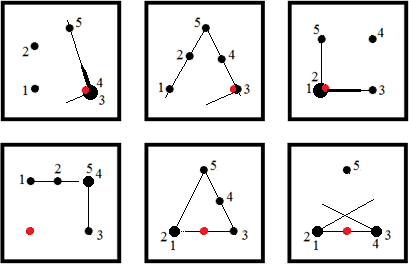}
\end{center}

Now regarding the charts stabilized by the lines: 
we have here two charts with triple points: the first chart
with $P_3=P_4=P_6$ and the third chart with $P_1=P_2=P_6$.
So it might be related to the 4.2.3 type, and  we might expect an extra chart with double 
points $P_1=P_2$ and $P_3=P_4$, point $P_6$
on the line connecting the two double points, 
and $P_5$ outside that line.

This extra chart is defined by the triple ratio $\left\{ P_1,P_3,P_5; P_2,P_4,P_6 \right\}$.
This is the only chart on which $\ell_{1,2},\ell_{3,4},\ell_{1,5},\ell_{3,5}$ are placed generically.

\subsubsection{$\mathcal{E}.2$}
Consider sub-case of the case $\mathcal{E.b}$, in which $\ell_{3,6}$ contains $P_1$ in the first two charts
(which happens simultaneously because of  $\left[P_1,P_2;P_6,P_5\right]_{P_3}$).

In the fourth chart, on which points $P_2,P_3,P_4,P_6$ are placed generically.
Since $\left[P_6,P_4;P_2,P_5\right]_{P_3}=0=\left[P_3,P_4;P_2,P_5\right]_{P_6}$, points $P_4$ and $P_5$ coincide. 
Since $\left[P_1,P_3;P_4,P_6\right]_{P_2}=0=\left[P_1,P_2;P_6,P_4\right]_{P_3}$, the point $P_1$
is placed at the intersection of $\ell_{2,4}$ and $\ell_{3,6}$.
Since $\left[P_2,P_3;P_6,P_5\right]_{P_4}=0$, the line $\ell_{4,5}$ has to contain $P_3$.

In the fifth chart,  $P_1,P_4,P_5,P_6$ are placed generically. Points $P_1$ and $P_2$ coincide, 
because $\left[P_1,P_5;P_2,P_6\right]_{P_4}=0=\left[P_1,P_4;P_2,P_6\right]_{P_5}$.
Notice also that $\left[P_3,P_4;P_6,P_5\right]_{P_1}=0=\left[P_3,P_1;P_4,P_6\right]_{P_5}$, 
so $P_3$ must be on $P_1P_6$ and $P_5P_4$.
The line $\ell_{1,2}$ must contain $P_5$, because $\left[P_4,P_2;P_6,P_5\right]_{P_1}=0$.

In the sixth chart, $P_1, P_2, P_3,P_6$ are placed generically.
For $\{i,j\}\in \{3,6\}$ and $k\in \{4,5\}$, we have $\left[P_k,P_1;P_2,P_j\right]_{P_i}=0$, so
$P_4$ and $P_5$ coincide with $P_2$.
The line $\ell_{4,5}$ contains $P_3$ because $\left[P_3,P_1;P_5,P_6\right]_{P_4}=0$.
The lines $\ell_{2,4}$ and $\ell_{2,5}$ contain $P_1$, because $\left[P_3,P_1;P_6,P_i\right]_{P_2}=0$ for $i\in\{4,5\}$.

\begin{center}
	\includegraphics{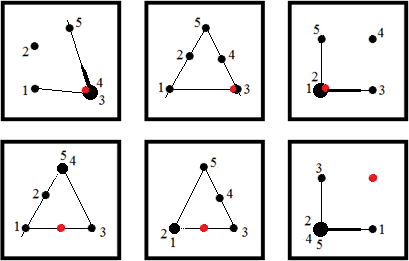}
\end{center}

This case has a lot of combinatorial symmetry.

Now we have to check regarding the charts which are stabilized by lines, whether they exist and 
whether they are interesting. The cue is the charts with only 4 distinct points, 
and here we have three such charts, each pair hints at a chart stabilized by lines of type (4.2.3),
so we start with the pair consisting of the first and the third chart with triple points 
$P_3=P_4=P_6$ and $P_1=P_2=P_6$. So, we should consider a new chart with double points $P_1=P_2$ and
$P_3=P_4$, with point $P_6$ on the line connecting the double points, and a point $P_5$ outside that line.
The chart is governed by the triple ratio $\left\{1,3,5;2,4,6\right\}$. This triple ratio is undefined on all 
the charts we have, so this time we should really add it. Another way to see it is that we need a chart with
$\ell_{1,2},\ell_{3,4},\ell_{1,5},\ell_{3,5}$ placed generically.

By symmetry, we should add two more charts of that kind, governed by triple ratios  
$\left\{1,3,5;4,6,2\right\}$ and $\left\{1,3,5;6,2,4\right\}$.

\begin{center}
	\includegraphics{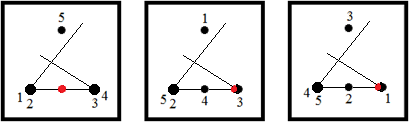}
\end{center}

\subsubsection{$\mathcal{E}.3$}
Consider sub-case of the case $\mathcal{E}.b$, in which $\ell_{3,6}$ 
contains $P_5$ in the first two charts.

In the fourth chart $P_1,P_2,P_3,P_6$ are placed generically.
For $i\in\{4,5\}$ we have $\left[P_2,P_3;P_i,P_6\right]_{P_1}=0=\left[P_1,P_i;P_2,P_6\right]_{P_3}$, 
therefore points $P_4$ and $P_5$ coincide with the intersection point of $P_1P_2$ and $P_3P_6$. 
Notice that $\left[P_1,P_3;P_6,P_4\right]_{P_5}=0$. In the fourth chart, if $\ell_{4,5}$ would be anything except $P_3P_4$, then  
$\left[P_1,P_3;P_6,P_4\right]_{P_5}$ would be infinity. Hence, $\ell_{4,5}$ contains $P_3$.

In the fifth chart $P_1,P_4,P_5,P_6$ are placed generically. 
Points $P_1$ and $P_2$ coincide, because
$\left[P_1,P_4;P_2,P_6\right]_{P_5}=0=\left[P_1,P_5;P_2,P_6\right]_{P_4}$. 
The point $P_3$ is the intersection of $P_1P_6$ and $P_4P_5$, since 
$\left[P_1,P_3;P_6,P_4\right]_{P_5}=0=\left[P_3,P_4;P_6,P_5\right]_{P_1}$. The line $\ell_{1,2}$ contains $P_5$ because $\left[P_2,P_3;P_5,P_4\right]_{P_1}=0$.

In the sixth chart $P_1,P_3,P_4,P_6$ are placed generically. 
Points $P_1$ and $P_2$ coincide,
since $\left[P_1,P_4;P_2,P_6\right]_{P_3}=0=\left[P_1,P_4;P_2,P_3\right]_{P_6}$.
Points $P_4$ and $P_5$ coincide,
because $\left[P_1,P_4;P_6,P_5\right]_{P_3}=0=\left[P_1,P_4;P_3,P_5\right]_{P_6}$.
The line $\ell_{1,2}$ contains $P_4$ in the sixth chart, because 
$\left[P_2,P_3;P_4,P_6\right]_{P_1}=0$.
The line $\ell_{4,5}$ contains $P_3$, because $\left[P_1,P_3;P_6,P_5\right]_{P_4}=0$.
\begin{center}
	\includegraphics{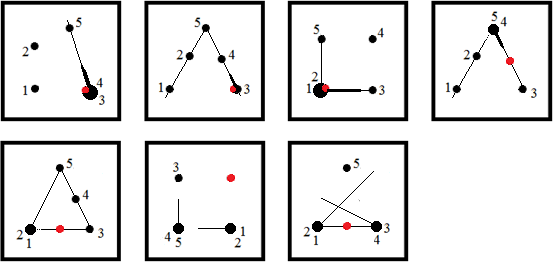}
\end{center}

We still have to check whether we need more charts stabilized by the lines.
We had 3 charts with at most 4 distinct points: two charts with a triple point 
and one with a double point. So, just one chart with two double points is 
not enough to get a chart of type (3.5). However we might have a chart of type (4.2.3),
with double points $P_1=P_2$ and $P_3=P_4$, points $P_6$ on the line connecting them,
and point $P_6$ outside that line. We need to have that chart because it is 
the only chart on which the four 
lines $\ell_{1,2},\ell_{1,5},\ell_{3,4},\ell_{3,5}$ are placed generically.
The chart is defined by the triple ratio $\left\{ P_1,P_3,P_5; P_2,P_4,P_6\right\}$, 
which happens to be undefined by all the other 6 charts,
so we should really add this chart.

\subsubsection{$\mathcal{E}.c$}
Like in all subtypes of E, on the first chart $P_1,P_2,P_3,P_5$ are placed generically,
while while $P_3$ and $P_4$ coincide; on the second chart $P_1,P_2,P_3,P_4$ are placed generically, while $P_5$ is at the intersection of $\ell_{1,2}$ and $\ell_{3,4}$;
on the third chat $P_1,P_3,P_4,P_5$ are placed generically, while $P_1$ and $P_2$ coincide.
In this particular case also $P_6$ coincides with $P_3$ on all these 3 charts.
Now regarding the lines related to the pairs of colliding points.
On the third chart, if $\ell_{3,6}$ does not contain $P_1$, then all cross-ratios are computable from that chart alone; hence the Chow orbit is 
a smooth point of the moduli space, 
like in those cases with a chart having at most 3 marked points on one line.
Hence we shall consider only the case where $\ell_{3,6}$ contains $P_1$ on the third chart.
Also on the third chart, there is $\ell_{1,2}$ which contains $P_5$
as might be seen from $\left[P_2,P_3;P_5,P_4\right]_{P_1}=0$ as in the second chart.

On the first chart, the lines $\ell_{3,4}$ and $\ell_{4,6}$ have to contain $P_5$
because $\left[P_1,P_5;P_2,P_4\right]_{P_3}=0=\left[P_1,P_5;P_2,P_4\right]_{P_6}$ 
as might be seen from the second chart.

\begin{center}
	\includegraphics{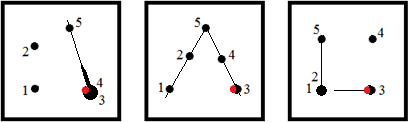}
\end{center}

The line which should yet be discussed is $\ell_{3,6}$ on the second chart, 
which should be chosen consistently with $\ell_{3,6}$ on the first chart 
since $\left[P_1,P_2;P_5,P_6\right]_{P_3}$ must be consistent.
There are different ways to choose the direction of $\ell_{3,6}$:
\begin{itemize}
	\item ($\mathcal{E}.4$) generic,
	\item ($\mathcal{E}.5$) containing $P_1$ (or $P_2$, but it is similar), 
	\item ($\mathcal{E}.6$) containing $P_4$ and $P_5$.
\end{itemize}

We shall list all these 3 cases, and in all cases some extra charts 
on which $P_3$ and $P_6$ are distinct, should be added.

\subsubsection{$\mathcal{E}.4$}
In addition to the three charts which we have already described, we shall add a fourth 
chart with $P_1,P_2,P_3, P_6$ placed generically. Since for $i\in\{4,5\}$ we have 
$\left[P_2,P_3;P_i,P_6\right]_{P_1}=0$, the line $\ell_{1,2}$ contains $P_4$ and $P_5$.
The location of $P_4$ and $P_5$ on $\ell_{1,2}$ is the same and is governed by 
the cross-ratio
$\left[P_1,P_2;P_4,P_6\right]_{P_3}=\left[P_1,P_2;P_5,P_6\right]_{P_3}\neq 0$. 

\begin{center}
	\includegraphics{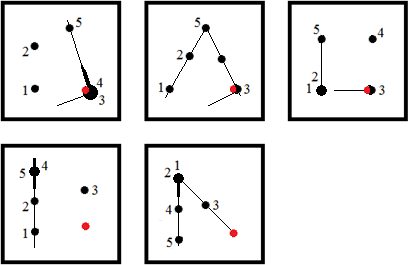}
\end{center}

We might further discuss the direction of $\ell_{4,5}$ in the fourth chart and 
also the fifth chart with $P_3,P_4,P_5,P_6$ placed generically; however we don't have to, 
since in the fourth chart, it is as if we add $P_5$ to 
a simple degeneration of the other points, which is the starting points of $\mathcal{A,B,C}$ cases 
which we have already listed (actually, in our case it is $\mathcal{C}2$ but it doesn't matter). 

\subsubsection{$\mathcal{E}.5$}
On the first chart, $P_1,P_2,P_3,P_5$ are placed generically,
$P_4$ and $P_6$ coincide with $P_3$; the lines $\ell_{3,4},\ell_{4,6}$ 
contain $P_5$, and $\ell_{3,6}$ contains $P_1$.

On the second chart, $P_1,P_2,P_3,P_4$ are placed generically,
$P_6$ coincides with $P_3$, and $P_5$ is at the intersection of
$\ell_{1,2}$ and $\ell_{3,4}$; the line $\ell_{3,6}$ contains $P_1$.

On the third chart, $P_1,P_3,P_4,P_5$ are placed generically, 
$P_2$ coincides with $P_1$ and $P_6$ coincides with $P_3$;
$\ell_{1,2}$ contains $P_5$ and $\ell_{3,6}$ contains $P_1$.

We have to add a fourth chart, with 
$P_1,P_2,P_3,P_6$ placed generically. For $i\in \{4,5\}$ we get 
$\left[P_2,P_3;P_i,P_6\right]_{P_1}$=0, so $P_i$ is on $\ell_{1,2}$;
also $\left[P_2,P_1;P_i,P_6\right]_{P_3}$, hence $P_i$ coincides with $P_2$.
All the special lines of the triple point contain $P_1$, since from the second 
and the third charts
\[\left[P_1,P_3;P_4,P_6\right]_{P_2}=\left[P_1,P_3;P_5,P_6\right]_{P_2}
=0=\left[P_1,P_3;P_5,P_6\right]_{P_4}.\]

There is a fifth chart, with 
$P_2 ,P_3, P_4, P_6$ placed generically. Then $P_1$ is at the intersection 
of $\ell_{2,4}$ and $\ell_{3,6}$, since
$\left[P_1,P_3;P_2,P_6\right]_{P_4}=0=\left[P_2,P_1;P_4,P_6\right]_{P_3}$.
Also, $P_5$ and $P_4$ coincide since 
$\left[P_4,P_2;P_5,P_6\right]_{P_3}=0=\left[P_4,P_2;P_5,P_3\right]_{P_6}$.
The line of the double point $\ell_{4,5}$ must contain $P_2$, since 
$\left[P_2,P_3;P_5,P_6\right]_{P_4}=0$.

We also have to add a sixth chart, with 
$P_3 ,P_4, P_5, P_6$ placed generically.
Here points $P_1$ and $P_2$ are located at the intersection of $\ell_{3,6}$
and $\ell_{4,5}$, because for $i\in \{1,2\}$, we have
$\left[P_i,P_4;P_6,P_5\right]_{P_3}=0=\left[P_i,P_3;P_5,P_6\right]_{P_4}$.
The line $\ell_{1,2}$ contains $P_4$, because otherwise
$\left[P_2,P_4;P_5,P_6\right]_{P_1}$ would be $\infty$ while it has 
to be zero (or undefined on this chart).

\begin{center}
	\includegraphics{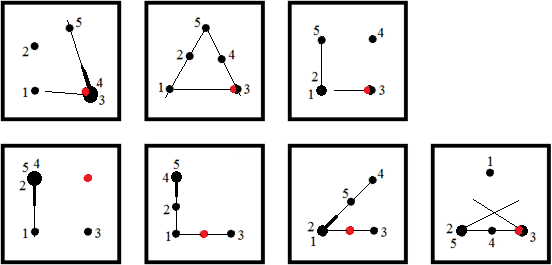}
\end{center}

No we must consider if we have to add any charts stabilized by lines.
There might not be configuration of type 3.5 since there is only one chart with 
two double points, and not 3 as required.
There might be a configuration of the type 4.2.3, since we do have a chart with 
two triple points. The respective roles of the points on this chart 
are easily determined from the general pattern of 4.2.3. This chart indeed exists,
since on all other charts the four lines 
$\ell_{3,6},\ell_{2,5},\ell_{1,2},\ell_{1,3}$
are not placed generically.

Anyway, this case is symmetric to $\mathcal{E}.3$, if we rename $P_1,P_2,...,P_6$
by $P_5,P_4, P_1, P_6, P_3, P_2$ respectively.

\subsubsection{$\mathcal{E}.6$}
We start with 3 given charts of this case.

In the first chart $P_1,P_2,P_3,P_5$ are placed generically,
points $P_4$ and $P_6$ coincide with $P_3$; 
the lines related to the triple point $\ell_{3,4}$, $\ell_{3,6}$,
and $\ell_{4,6}$ are all containing $P_5$.

In the second chart $P_1,P_2,P_3,P_4$ are placed generically,
$P_5$ is at the intersection of $\ell_{1,2}$ and $\ell_{3,4}$,
and $P_6$ coincides with $P_3$; the line of the double point $\ell_{3,6}$
contains $P_4$.

In the third chart $P_1,P_3,P_4,P_5$ are placed generically,
$P_2$ coincides with $P_1$, $P_6$ coincides with $P_3$; the line $\ell_{1,2}$
contains $P_5$, the line $\ell_{3,6}$ contains $P_1$.

\begin{center}
	\includegraphics{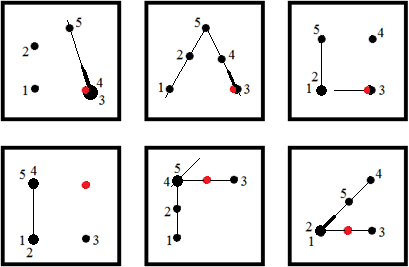}
\end{center}

Consider the fourth chart in which $P_1,P_3,P_4,P_6$ are placed generically. 
Then $P_2$ and $P_1$ coincide, because
$\left[P_1,P_6;P_2,P_5\right]_{P_3}=0=\left[P_1,P_3;P_2,P_5\right]_{P_6}$.
Also $P_4$ and $P_5$ coincide, because 
$\left[P_3,P_4;P_6,P_5\right]_{P_1}=0=\left[P_1,P_4;P_6,P_5\right]_{P_3}$.
Regarding the lines of the double points:
$\left[P_2,P_3;P_5,P_6\right]_{P_1}=0$, hence $\ell_{1,2}$ contains $P_5$; 
also $\left[P_1,P_3;P_5,P_6\right]_{P_4}=0$ hence $\ell_{4,5}$ contains $P_1$.

Consider the fifth chart in which 
$P_1,P_2,P_3,P_6$ are placed generically. 
Then $P_4$ and $P_5$ both coincide with 
the intersection of $\ell_{1,2}$ and 
$\ell_{3,6}$. 
Indeed, $\left[P_2,P_3;P_i,P_6\right]_{P_1}=0=
\left[P_i,P_1;P_6,P_2\right]_{P_3}$
for $i\in \{4,5\}$ as can be easily 
seen from the second chart.
The line $\ell_{4,5}$ might be anything;
it is governed by the triple ratio $\{P_1,P_3,P_4;P_6,P_2,P_5\}$
which is undefined on all previous charts.

The sixth chart has to be considered, with $P_3,P_4,P_5,P_6$
placed generically. Points $P_1$ and $P_2$ are placed at 
the intersection of $\ell_{3,6}$ and $\ell_{4,5}$,
because for $i\in\{1,2\}$, we have
$\left[P_i,P_4;P_6,P_5\right]_{P_3}=0=
\left[P_i,P_3;P_5,P_6\right]_{P_4}$.
The line $\ell_{1,2}$ has to contain $P_4$ and $P_5$, otherwise
$\left[P_2,P_4;P_3,P_5\right]_{P_1}$ would be 
zero, but on the third chart it is $\infty$.

It remains to verify whether there are charts
stabilized by lines. There is only one chart 
with a triple point hence it can not be 4.2.3.
There are only two charts with two double points,
hence it is not 3.5, which requires three charts of that kind.

\subsubsection{$\mathcal{E}.d$}

On the first chart, $P_1,P_2,P_3,P_5$ are placed generically,
$P_4$ and $P_6$ coincide with $P_3$. The lines $\ell_{3,4}$
and $\ell_{4,6}$ contain $P_5$.

On the second chart, $P_1,P_2,P_3,P_4$ are placed generically,
$P_5$ is at the intersection of $\ell_{1,2}$ with $\ell_{3,4}$,
$P_6$ coincided with $P_3$.

On third third chart, $P_1,P_3,P_4,P_5$ are placed generically,
$P_2$ coincides with $P_1$,
$P_6$ is at the intersection of $\ell_{1,3}$ with $\ell_{4,5}$;
the line $\ell_{1,2}$ contains $P_5$.

\begin{center}
	\includegraphics{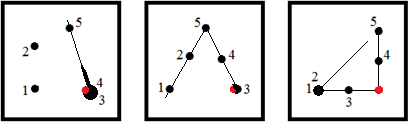}
\end{center}

The line which should yet be discussed is $\ell_{3,6}$ on the second chart, 
which should be chosen consistently with $\ell_{3,6}$ on the first chart 
since $\left[P_1,P_2;P_5,P_6\right]_{P_3}$ must be consistent.
There are different ways to choose the direction of $\ell_{3,6}$:
\begin{itemize}
	\item ($\mathcal{E}.7$) generic,
	\item ($\mathcal{E}.8$) containing $P_1$ (or $P_2$, but it is similar), 
	\item ($\mathcal{E}.9$) containing $P_4$ and $P_5$.
\end{itemize}

We shall list all these 3 cases.

\subsubsection{$\mathcal{E}.7$}
Consider the case that line $\ell_{3,6}$ is generic.

Consider a fourth chart on which $P_1,P_2, P_3, P_6$
are placed generically. It is easy to see from the second chart that 
$\left[P_2,P_3;P_i,P_6\right]_{P_1}=0$ for $i\in\{4,5\}$
so $P_4$ and $P_5$ are in $\ell_{1,2}$ in the fourth chart.
Notice also that 
\[\left[P_1,P_2;P_4,P_6\right]_{P_3}=\left[P_1,P_2;P_5,P_6\right]_{P_3}\notin\{0,1,\infty\}\]
hence $P_4$ and $P_5$ coincide at a generic point of $\ell_{1,2}$.

\begin{center}
	\includegraphics{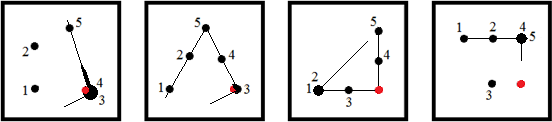}
\end{center}

We can also verify the direction of $\ell_{4,5}$ and we should add more charts,
but already at this moment we may say that on the fourth chart 
the 5 marked points except $P_5$ form a simple degeneration, 
so this case is symmetric to one of the cases in $\mathcal{A,B,C}$ families.

\subsubsection{$\mathcal{E}.8$}

Here we are given 3 charts, and several more charts have to be added.

On the first chart, $P_1,P_2,P_3,P_5$ are placed generically,
$P_4$ and $P_6$ coincide with $P_3$. The lines $\ell_{3,4}$
and $\ell_{4,6}$ contain $P_5$, and $\ell_{3,6}$ contains $P_1$.

On the second chart, $P_1,P_2,P_3,P_4$ are placed generically,
$P_5$ is at the intersection of $\ell_{1,2}$ with $\ell_{3,4}$,
$P_6$ coincided with $P_3$; the line $\ell_{3,6}$ contains $P_1$.

On third third chart, $P_1,P_3,P_4,P_5$ are placed generically,
$P_2$ coincides with $P_1$,
$P_6$ is at the intersection of $\ell_{1,3}$ with $\ell_{4,5}$;
the line $\ell_{1,2}$ contains $P_5$.

\begin{center}
	\includegraphics{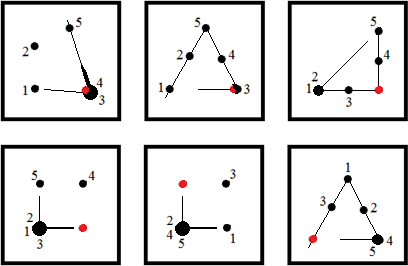}
\end{center}

On the fourth chart points $P_1,P_4,P_5,P_6$ are placed generically.
For $i\in\{1,2\}$ and $\{j,k\}\in\{4,5\}$ we can easily see on 
the third chart that $\left[P_k,P_3;P_6,P_i\right]_{P_j}=0$; 
therefore $P_3$ coincides with $P_1$ and $P_2$. 
Regarding the lines of this triple point: 
$\ell_{1,2}$ contains $P_5$ because $\left[P_2,P_4;P_5,P_6\right]_{P_1}=0$,
$\ell_{1,3}$ and $\ell_{2,3}$ contain $P_6$ because 
$\left[P_i,P_4;P_6,P_5\right]_{P_3}=0$ for $i\in\{1,2\}$.

On the fifth chart $P_1,P_2,P_3,P_6$ are placed generically.
Points $P_4$ and $P_5$ coincide with $P_2$ since for 
$i\in\{4,5\}$ and $\{j,k\} = \{3,6\}$ we can easily see from the second chart that 
$\left[P_k,P_2;P_1,P_i\right]_{P_j}=0$.
Regarding the lines related to this triple point: 
$\ell_{4,5}$ contains $P_6$ because $\left[P_1,P_4;P_3,P_6\right]_{P_5}=0$,
$\ell_{2,4}$ and $\ell_{2,5}$ contain $P_1$ because 
$\left[P_1,P_3;P_i,P_6\right]_{P_2}=0$ for $i\in\{4,5\}$.

On the sixth chart $P_2,P_3,P_4,P_6$ are placed generically;
$P_5$ coincides with $P_4$ because 
$\left[P_4,P_2;P_5,P_6\right]_{P_3}=0=\left[P_4,P_2;P_5,P_3\right]_{P_6}$.
The point $P_1$ coincides with the intersection of $\ell_{3,6}$ and $\ell_{2,4}$,
because 
\[\left[P_2,P_1;P_4,P_3\right]_{P_6}=0=\left[P_1,P_3;P_4,P_6\right]_{P_2}.\]

There are more charts stabilized by lines. Anyway, we can skip the discussion of this
case, since it is symmetric to $\mathcal{E}.2$ by swapping $P_3$ and $P_6$.

\subsubsection{$\mathcal{E}.9$}
On the first chart, $P_1,P_2,P_3,P_5$ are placed generically,
$P_4$ and $P_6$ coincide with $P_3$. The lines $\ell_{3,4}$,
$\ell_{4,6}$, $\ell_{3,6}$ contain $P_5$. 

On the second chart, $P_1,P_2,P_3,P_4$ are placed generically,
$P_5$ is at the intersection of $\ell_{1,2}$ with $\ell_{3,4}$,
$P_6$ coincided with $P_3$, the line $\ell_{3,6}$ contains $P_5$.

On third third chart, $P_1,P_3,P_4,P_5$ are placed generically,
$P_2$ coincides with $P_1$,
$P_6$ is at the intersection of $\ell_{1,3}$ with $\ell_{4,5}$;
the line $\ell_{1,2}$ contains $P_5$.

\begin{center}
	\includegraphics{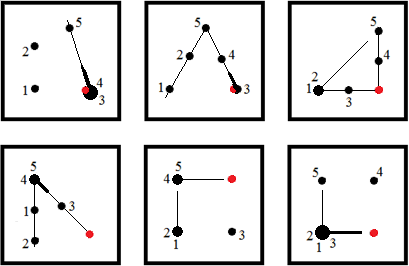}
\end{center}

We add a fourth chart with $P_1,P_2,P_3,P_6$ placed generically.
It is easy to see from the second chart that
$\left[P_i,P_3;P_2,P_6\right]_{P_1}=0=\left[P_1,P_i;P_2,P_6\right]_{P_3}$,
for $i\in\{4,5\}$,
therefore $P_4$ and $P_5$ should be placed on the intersection of $\ell_{1,2}$
and $\ell_{3,6}$. The line $\ell_{4,5}$ should contain $P_3$; otherwise
$\left[P_3,P_4;P_2,P_6\right]_{P_5}$ would be $\infty$ while it should be 0 
like in the third chart (or undefined like in this one).

We add a fifth chart with $P_1,P_3,P_4,P_6$ placed generically.
Then $P_2$ coincides with $P_1$ since 
$\left[P_1,P_4;P_2,P_6\right]_{P_3}=0=\left[P_1,P_4;P_2,P_3\right]_{P_6}$.
Furthermore, $P_4$ coincides with $P_5$, 
because $\left[P_4,P_1;P_5,P_6\right]_{P_3}=0=\left[P_4,P_3;P_5,P_6\right]_{P_1}$.
The line $\ell_{1,2}$ contains $P_5$, because $\left[P_2,P_3;P_5,P_6\right]_{P_1}=0$.
The line $\ell_{4,5}$ contains $P_6$, because $\left[P_1,P_5;P_3,P_6\right]_{P_4}=0$.

We also have to add a sixth chart with $P_1,P_4,P_5,P_6$ placed generically.
The points $P_1,P_2,P_3$ coincide, since 
for $i\in\{1,2\}$ and $\{j,k\}=\{4,5\}$ it can be easily seen 
from the third chart that $\left[P_j,P_i;P_6,P_3\right]_{P_k}=0$.
The lines $\ell_{1,3}$ and $\ell_{2,3}$ contain $P_6$, because 
$\left[P_3,P_4;P_6,P_5\right]_{P_i}=0$, for $i\in\{1,2\}$. 
The lines $\ell_{1,2}$ contains $P_5$, because$\left[P_2,P_4;P_5,P_6\right]_{P_1}=0$.

There is also a chart stabilized by lines. Anyway, we can skip 
the further discussion of this
case, since it is symmetric to $\mathcal{E}.3$ by swapping $P_3$ and $P_6$.

\subsubsection{$\mathcal{F}$}
In the second chart, points $P_1,P_2, P_3,P_4$ are placed generically, point $P_5$ is the intersection point of $P_1P2$ and $P_3P_4$,
and point $P_6$ coincides with point $P_5$. 
In the first chart, $P_1,P_2,P_3,P_5$ are placed generically, $P_4$ coincides with $P_3$, and $\ell_{3,4}$ contains $P_5$.
In the third chart, $P_1,P_3,P_4,P_5$ are placed generically, $P_2$ coincides with $P_1$, and $\ell_{1,2}$ contains $P_5$.

\begin{center}
	\includegraphics{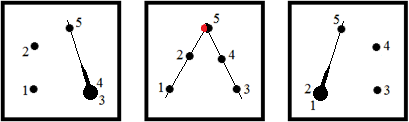}
\end{center}

The possibilities to place point $P_6$ will be further discussed, but before let us remark about the way we have organized our case-checking. 
We have chosen a chart on which some 5 marked points are distinct (like chart 2 in this case) and tried to add a sixth marked point. 
Now, we are in the last and the most degenerate case, which might further be divided into sub-cases. However, if we consider a sub-case
in which some 5 points are distinct, and together with the remaining point they are placed in a less degenerate way, it means that
we have already listed that sub-case, up to renaming of points. This idea will help us to skip some sub-cases.

Notice that $\left[P_5,P_1;P_6,P_2\right]_{P_3}=0$. Hence in the first chart, $\ell_{3,6}$ contains $P_5$.
Hence  $P_6$ must belong to the line $P_5P_3$. If $P_6$ is anything except $P_3$ and $P_5$, we get in the first chart 5 distinct marked points 
and a less degenerate configuration - even if $P_6$ is the intersection point of $P_1P_2$ and $P_3P_5$, it is the only point on that intersection.

So, the only cases we have to consider, is that in the first chart $P_6$ coincides with $P_5$ or $P_3$.
In the second case, we saw that $\ell_{3,6}$ contains $P_5$ and similarly, $\ell_{4,6}$ contains $P_5$. 

Within the case when in the first chart $P_5$ coincides with $P_6$, consider the sub-case that on $\ell_{5,6}$ does not contain $P_3$.
So $\ell_{5,6}$ in the first chart is either in a generic direction, or contains $P_1$, or $P_2$, but by symmetry we may assume it is not $P_2$.
There should be a fourth chart on which $P_2,P_3,P_5,P_6$ are placed generically. Then $P_4$ still coincides with $P_3$ because 
$\left[P_3,P_1;P_4,P_5\right]_{P_6}=0=\left[P_3,P_1;P_4,P_6\right]_{P_5}$. 
The point $P_1$ is on the line $P_2P_3$ in the fourth chart, 
since $\left[P_1,P_5;P_2,P_6\right]_{P_3}=0$. However, $P_2$ is distinct from the other marked points on that line, since 
$\left[P_1,P_2;P_3,P_6\right]_{P_5}\neq 0$. So, on this fourth chart we 5 distinct marked points in a less degenerate situation.

To summarize, if in the first chart $P_5$ coincides with $P_6$, then $\ell_{5,6}$ should contain $P_3$ unless it is similar to a previously discussed sub-case.

Whatever we discussed regarding the first chart, is relevant to the third chart as well. So, we may assume, that in the third chart,
either $P_6$ coincides with $P_5$ and $\ell_{5,6}$ contains $P_1$, or $P_6$ coincides with $P_1$ and $P_2$, and then $\ell_{1,6},\ell_{2,6}$ contain $P_5$.

The possibilities for first and third chart can combine into 4 cases, but by symmetry it is enough to consider 3 cases:

\begin{itemize}
	\item (1) Both in the first chart and the third chart, $P_6$ coincides with $P_5$.
	\item (2) In the first chart $P_6$ coincides with $P_5$, in the third chart $P_6$ coincides with $P_1$. 
	\item (3) In the first chart $P_6$ coincides with $P_3$, in the third chart $P_6$ coincides with $P_1$.
\end{itemize}
We shall discuss all these cases separately. 

\subsubsection{$\mathcal{F}.1$}
In the first three charts, $P_5$ coincides with $P_6$. In the first and the third chart, 
as we saw, all lines defined by pairs of colliding points have to contain the other double point.

In the fourth chart, points $P_1,P_2,P_5,P_6$ are placed generically. 
For $i\in\{3,4\}$ we have $\left[P_i,P_5;P_2,P_6\right]_{P_1}=0=\left[P_i,P_1;P_6,P_2\right]_{P_5}$,
so in $P_3$ and $P_4$ coincide with the intersection point of $P_1P_2$ and $P_5P_6$.

In the fifth chart $P_3,P_4,P_5,P_6$ are placed generically, while $P_1$ and $P_2$ coincide with 
the intersection point of lines $P_3P_4$ and $P_5P_6$ (and the reason is similar to the fourth chart so we don't repeat it).

In the sixth chart, points $P_1,P_3,P_5,P_6$ are placed generically. 
Since $\left[P_1,P_3;P_2,P_5\right]_{P_5}=0=\left[P_i,P_1;P_6,P_2\right]_{P_5}$, points $P_1$ and $P_2$ coincide.
Similarly, points $P_3$ and $P_4$ coincide. 
Notice that points $\left[P_2,P_5;P_3,P_6\right]_{P_1}=0$, hence $\ell_{1,2}$ coincides with $P_1P_3$ in the sixth chart. 
Similarly, $\ell_{3,4}$ is also the same line.

We have not specified the lines of coinciding points in the second, fourth and fifth charts; 
it turns out they may be anything without conflicting with any other charts.

\begin{center}
	\includegraphics{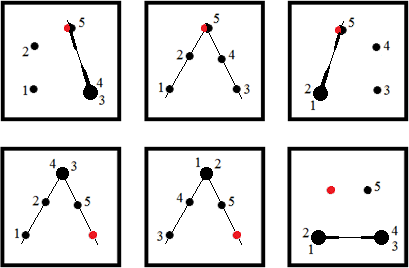}
\end{center}

We should add an extra chart stabilized by lines of type (3.5). 

Indeed, there should be a chart in which $\ell_{1,2}$, $\ell_{1,5}$, $\ell_{3,4}$, and $\ell_{5,3}$
are placed generically, and we don't have such a chart yet. So we can start with 
an equilateral triangle $P_1P_3P_5$,
and add two altitudes $\ell_{1,2}$ and $\ell_{3,4}$, even before choosing $P_2$ and $P_4$.
Then on those altitudes we have to choose $P_2$ and $P_4$ so that 
$\left[P_4,P_2;P_3,P_5\right]_{P_1}=0=\left[P_2,P_4;P_1,P_5\right]_{P_3}$,
hence $P_1=P_2$ and $P_3=P_4$. It remains to place $P_6$, however
$\left[P_6,P_2;P_5,P_3\right]_{P_1}=0=\left[P_6,P_1;P_5,P_4\right]_{P_3}$, so $P_6=P_5$. 
Hence it remains yet to choose the direction of $\ell_{5,6}$, which is governed by the
triple ratio $\left\{P_1,P_3,P_5;P_2,P_4,P_6 \right\}$, which is no defined by any of the previous charts.

Similarly we could place the quadruple $\left\{\ell_{1,2},\ell_{1,3},\ell_{5,6},\ell_{3,5}\right\}$ or the quadruple
$\left\{\ell_{1,3},\ell_{3,4},\ell_{1,5},\ell_{5,6}\right\}$ generically. If the triple ratio is 
not 0 or $\infty$, we get an equivalent chart, but if it is one of those two special cases, 
we get more charts.

\subsubsection{$\mathcal{F}.2$}
In the second chart, points $P_1,P_2, P_3,P_4$ are placed generically, $P_5$ is the intersection point of $P_1P2$ and $P_3P_4$,
and point $P_6$ coincides with point $P_5$. 
In the first chart, $P_1,P_2,P_3,P_5$ are placed generically, $P_4$ coincides with $P_3$, $P_6$ coincides with $P_5$, 
the lines $\ell_{3,4}$ and $\ell_{5,6}$ coincide with $P_3P_5$.
In the third chart, $P_1,P_3,P_4,P_5$ are placed generically, $P_2$ and $P_6$ coincide with $P_1$, the lines $\ell_{1,2}$, $\ell_{1,6}$ and $\ell_{2,6}$ contains $P_5$.
the third chart $P_6$ coincides with $P_1$. 


Consider a fourth chart, on which points $P_1,P_3,P_4,P_6$ are placed generically.
Since $\left[P_1,P_4;P_2,P_6\right]_{P_3}=0=\left[P_1,P_3;P_2,P_6\right]_{P_4}$, points $P_1$ and $P_2$ must coincide.
Since $\left[P_3,P_5;P_4,P_6\right]_{P_1}=0=\left[P_1,P_3;P_6,P_5\right]_{P_4}$, we conclude that $P_5$ is the point of intersection
of $P_1P_6$ and $P_3P_4$. We may discuss the direction of $\ell_{1,2}$ and some other charts, but anyway, in the fourth chart 
we have a configuration with six distinct marked points, which is less degenerate then case F, so we can skip this sub-case without further discussion.

\begin{center}
	\includegraphics{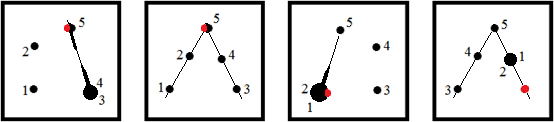}
\end{center}

\subsubsection{$\mathcal{F}.3$}
In the second chart, points $P_1,P_2, P_3,P_4$ are placed generically, $P_5$ is the intersection point of $P_1P_2$ and $P_3P_4$,
and $P_6$ coincides with $P_5$. 
In the first chart, $P_1,P_2,P_3,P_5$ are placed generically, $P_4$ coincides with $P_3$, $P_6$ coincides with $P_5$, 
the lines $\ell_{3,4}$ and $\ell_{5,6}$ coincide with $P_3P_5$.
In the third chart, $P_1,P_3,P_4,P_5$ are placed generically, $P_2$ and $P_6$ coincide with $P_1$, the lines $\ell_{1,2}$, $\ell_{1,6}$ and $\ell_{2,6}$ contains $P_5$.
the third chart $P_6$ coincides with $P_1$. 

Consider a fourth chart, with $P_1,P_3,P_4,P_6$ placed generically.
Since $\left[P_1,P_4;P_2,P_6\right]_{P_3}=0=\left[P_1,P_3;P_2,P_6\right]_{P_4}$, points $P_1$ and $P_2$ must coincide.
Since $\left[P_3,P_5;P_4,P_6\right]_{P_1}=0=\left[P_1,P_3;P_6,P_5\right]_{P_4}$, we conclude that $P_5$ is the point of intersection
of $P_1P_6$ and $P_3P_4$. We may discuss the direction of $\ell_{1,2}$ and some other charts, but anyway, in the fourth chart 
we have a configuration with five distinct marked points, 
which is less degenerate then case F, so we can skip this sub-case without further discussion.

\begin{center}
	\includegraphics{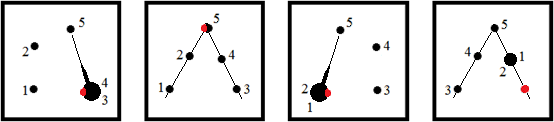}
\end{center}

\subsection{Summary of the case separation.}
So, after listing the cases for which smoothness is not entirely obvious, 
we get 13 distinct cases (actually 14 as $\mathcal{F}.1$ splits into sub-cases)
which correspond to some degenerate configurations: 
$\mathcal{A}.1$, $\mathcal{A}.2$, $\mathcal{A}.4$, $\mathcal{B}.2$, 
$\mathcal{B}.4$, $\mathcal{C}.2$, $\mathcal{C}.4$, $\mathcal{D}$, $\mathcal{E}.1$, 
$\mathcal{E}.2$, $\mathcal{E}.3$, $\mathcal{E}.6$, $\mathcal{F}.1$.

\section{Proof of smoothness of $\widetilde{Y}_6$ - case checking.}
For each chart in each of the above configurations, we may write down some 4 natural local parameters.
There will be some relations between these parameters, 
due to the fact that cross-ratios are the same in all maps.
When we consider a cotangent space, we have to verify it is 4-dimensional. 
So what we shall verify that the linear approximations of the relations on the coordinates 
are sufficient to decrease the dimension from $4k$ to $4$, where $k$ is the number of charts.

Before getting into the case-checking, we shall remark regarding the notations we shall use.

As each point appears on several charts, instead of denoting the points $P_1,P_2,...,P_6$,
we shall denote them $A,B,...,F$.
We shall have different manifestations of the same point on all charts, 
so the first marked point in chart $i$ will be denoted $A_i$, 
the second point in the chart $i$ will be denoted $B_i$ and so on.

The cross ratios such as $\left[A,B;C,D\right]_{E}$ should be independent of chart, 
and writing them down in coordinates in different charts will produce some relations 
between different coordinates we have. 

It is useful to have the formula 
\[u = \left[ {0,\infty ;u,1} \right] = {{\frac{u}{{u - \infty }}} \mathord{\left/
 {\vphantom {{\frac{u}{{u - \infty }}} {\frac{1}{{1 - \infty }}}}} \right.
 \kern-\nulldelimiterspace} {\frac{1}{{1 - \infty }}}}\]
readily available, to be able to translate coordinates into cross-ratios.

Another useful type of formula is 
\[{\left[ {A,B;C,D} \right]_E} = \frac{{\left| {A,C,E} \right| \cdot \left| {B,D,E} \right|}}
{{\left| {B,C,E} \right| \cdot \left| {A,D,E} \right|}}\]
where $\left| {A,C,E} \right|$ means determinant of the matrix with columns which 
are the projective coordinates of points $A, C, E$, etc.

We try to use points with many zeroes and ones as coordinates for the marked points.
If for instance $E=(0,0,1)$ or another vector of standard basis, 
the formula is reduced to $2\times 2$ determinants rather than $3\times 3$.

A few last notations before we delve into treating all special cases.

We shall write $u \sim v$ to denote that $u-v$ is near zero,
and we write $u \approx v$ to denote that $u-v$ is near zero up to second order.

It is useful when we investigate the dimension of cotangent space in coordinates,
for instance if $u=v(1+v)$ and $v \sim 0$ we may simplify and write $u \approx v$.

In some cases, the order of charts and of the names for the points will 
be different than in the previous list, to have a more symmetric picture.

\subsection{$\mathcal{A}.1$} 
In this case, there is a nice combinatorial symmetry, so we choose coordinates in a symmetric way.
\begin{center}
	\includegraphics{Case_A1.png}
\end{center}

\[\begin{array}{*{20}{c}}
\begin{array}{l}
{A_1} = \left( {1:0:0} \right)\\
{B_1} = \left( {1:1:1} \right)\,\\
{C_1} = \left( {0:1:0} \right)\\
{D_1} = \left( {{y_1}:1:{x_1}} \right)\\
{E_1} = \left( {0:0:1} \right)\\
{F_1} = \left( {{z_1}:{z_1}{t_1}:1} \right)
\end{array}&\vline& \begin{array}{l}
{A_2} = \left( {1:0:0} \right)\\
{B_2} = \left( {1:{z_2}:{z_2}{t_2}} \right)\,\\
{C_2} = \left( {0:1:0} \right)\\
{D_2} = \left( {1:1:1} \right)\\
{E_2} = \left( {0:0:1} \right)\\
{F_2} = \left( {{x_2}:{y_2}:1} \right)
\end{array}&\vline& \begin{array}{l}
{A_3} = \left( {1:0:0} \right)\\
{B_3} = \left( {1:{x_3}:{y_3}} \right)\\
{C_3} = \left( {0:1:0} \right)\\
{D_3} = \left( {{z_3}{t_3}:1:{z_3}} \right)\\
{E_3} = \left( {0:0:1} \right)\\
{F_3} = \left( {1:1:1} \right)
\end{array}
\end{array}\]
where 
\[ y_1,y_2,y_3,z_1,z_2,z_3 \sim 0,  \]
\[ x_1,x_2,x_3 \ne 0, \]
\[ t_1,t_2,t_3 \ne 0,1. \]

\[{x_1} = {\left[ {C,E;D,B} \right]_A} = {\left[ {\left( \begin{array}{l}
0\\
1\\
0
\end{array} \right),\left( \begin{array}{l}
0\\
0\\
1
\end{array} \right);\left( \begin{array}{l}
1\\
1\\
1
\end{array} \right),\left( \begin{array}{c}
1\\
{z_2}\\
{z_2}{t_2}
\end{array} \right)} \right]_
{\left( {\scriptstyle1\atop
        {\scriptstyle0\atop
         \scriptstyle0}} \right)}} =\]
\[ = \frac{{\left| {\left( {\begin{array}{*{20}{c}}
1&1\\
0&1
\end{array}} \right)} \right| \cdot \left| {\left( {\begin{array}{*{20}{c}}
0&{{z_2}}\\
1&{{z_2}{t_2}}
\end{array}} \right)} \right|}}{{\left| {\left( {\begin{array}{*{20}{c}}
1&{{z_2}}\\
0&{{z_2}{t_2}}
\end{array}} \right)} \right| \cdot \left| {\left( {\begin{array}{*{20}{c}}
0&1\\
1&1
\end{array}} \right)} \right|}} = \frac{1}{{{t_2}}}\]

\[{y_1} = {\left[ {C,A;D,B} \right]_E} = {\left[ {\left( \begin{array}{l}
0\\
1\\
0
\end{array} \right),\left( \begin{array}{l}
1\\
0\\
0
\end{array} \right);\left( \begin{array}{l}
1\\
1\\
1
\end{array} \right),\left( \begin{array}{c}
1\\
{z_2}\\
{z_2}{t_2}
\end{array} \right)} \right]_{\left( {\scriptstyle0\atop
{\scriptstyle0\atop
\scriptstyle1}} \right)}} =\]\[= \frac{{\left| {\left( {\begin{array}{*{20}{c}}
0&1\\
1&1
\end{array}} \right)} \right| \cdot \left| {\left( {\begin{array}{*{20}{c}}
1&1\\
0&{{z_2}}
\end{array}} \right)} \right|}}{{\left| {\left( {\begin{array}{*{20}{c}}
0&1\\
1&{{z_2}}
\end{array}} \right)} \right| \cdot \left| {\left( {\begin{array}{*{20}{c}}
1&1\\
0&1
\end{array}} \right)} \right|}} = {z_2}\]

\[{z_1} = {\left[ {E,A;F,B} \right]_C} = {\left[ {\left( \begin{array}{l}
0\\
0\\
1
\end{array} \right),\left( \begin{array}{l}
1\\
0\\
0
\end{array} \right);\left( \begin{array}{l}
{x_2}\\
{y_2}\\
1
\end{array} \right),\left( \begin{array}{c}
1\\
{z_2}\\
{z_2}{t_2}
\end{array} \right)} \right]_{\left( {\scriptstyle0\atop
{\scriptstyle1\atop
\scriptstyle0}} \right)}} =\]\[= \frac{{\left| {\left( {\begin{array}{*{20}{c}}
0&{{x_2}}\\
1&1
\end{array}} \right)} \right| \cdot \left| {\left( {\begin{array}{*{20}{c}}
1&1\\
0&{{z_2}{t_2}}
\end{array}} \right)} \right|}}{{\left| {\left( {\begin{array}{*{20}{c}}
0&1\\
1&{{z_2}{t_2}}
\end{array}} \right)} \right| \cdot \left| {\left( {\begin{array}{*{20}{c}}
1&{{x_2}}\\
0&1
\end{array}} \right)} \right|}} = {x_2}{z_2}{t_2}\]

\[{t_1} = {\left[ {A,C;F,B} \right]_E} = {\left[ {\left( \begin{array}{l}
1\\
0\\
0
\end{array} \right),\left( \begin{array}{l}
0\\
1\\
0
\end{array} \right);\left( \begin{array}{l}
{x_2}\\
{y_2}\\
1
\end{array} \right),\left( \begin{array}{c}
1\\
{z_2}\\
{z_2}{t_2}
\end{array} \right)} \right]_{\left( {\scriptstyle0\atop
{\scriptstyle0\atop
\scriptstyle1}} \right)}} =\]\[= \frac{{\left| {\left( {\begin{array}{*{20}{c}}
1&{{x_2}}\\
0&{{y_2}}
\end{array}} \right)} \right| \cdot \left| {\left( {\begin{array}{*{20}{c}}
0&1\\
1&{{z_2}}
\end{array}} \right)} \right|}}{{\left| {\left( {\begin{array}{*{20}{c}}
1&1\\
0&{{z_2}}
\end{array}} \right)} \right| \cdot \left| {\left( {\begin{array}{*{20}{c}}
0&{{x_2}}\\
1&{{y_2}}
\end{array}} \right)} \right|}} = \frac{{{y_2}}}{{{x_2}{z_2}}}\]

Since we have chosen all notations symmetrically, we may write two more similar quadruples of formulas:
\[\begin{array}{*{20}{c}}
\begin{array}{l}
{x_1}{t_2} = 1\\
{y_1} = {z_2}\\
{z_1} = {x_2}{z_2}{t_2}\\
{t_1}{x_2}{z_2} = {y_2}
\end{array}&\vline& \begin{array}{l}
{x_2}{t_3} = 1\\
{y_2} = {z_3}\\
{z_2} = {x_3}{z_3}{t_3}\\
{t_2}{x_3}{z_3} = {y_3}
\end{array}&\vline& \begin{array}{l}
{x_3}{t_1} = 1\\
{y_3} = {z_1}\\
{z_3} = {x_1}{z_1}{t_1}\\
{t_3}{x_1}{z_1} = {y_1}
\end{array}
\end{array}\]

From the first two lines we see that $t_i$ can be excluded and replaced by $\frac{1}{x_{i-1}}$ (and as 
$t_i\neq 0$, we can express $dt_i$ in terms of $dx_{i-1}$), and $z_i$ may be excluded as well.
So we can reduce to 6 coordinates $x_i,y_i$; the last two lines of formulas turn into 
\[{x_1}{y_3} = {x_2}{y_1} = {x_3}{y_2}\]
Since $y_i\sim 0 \neq x_i$, we conclude 
\[x_1 \cdot dy_3 =  x_2\cdot dy_1 = x_3\cdot dy_2\]
Hence the cotangent space is spanned by $dx_1,dx_2,dx_3$ and $dy_1$ and hence it is four-dimensional.

\subsection{$\mathcal{A}.2$} 
\begin{center}
	\includegraphics{Case_A2.png}
\end{center}
We choose local coordinates on the first 3 charts, as the last chart doesn't contain new parameters.

\[\small \begin{array}{*{10}{c}}
\begin{array}{l}
{A_1} = \left( {1:0:0} \right)\\
{B_1} = \left( {1:1:1} \right)\\
{C_1} = \left( {0:1:0} \right)\\
{D_1} = \left( {{x_1}:1:{y_1}} \right)\\
{E_1} = \left( {0:0:1} \right)\\
{F_1} = \left( {{z_1}:{z_1}{t_1}:1} \right)\\
{x_1},{z_1} \sim 0 \ne {y_1},{t_1}
\end{array}&\vline& \begin{array}{l}
{A_2} = \left( {1:0:0} \right)\\
{B_2} = \left( {1:{x_2}:{x_2}{y_2}} \right)\\
{C_2} = \left( {0:1:0} \right)\\
{D_2} = \left( {0:0:1} \right)\\
{E_2} = \left( {1:1:1} \right)\,\\
{F_2} = \left( {{z_2}:{1+t_2}:1} \right)\\
{x_2},{z_2},{t_2} \sim 0 \ne {y_2}
\end{array}&\vline& \begin{array}{l}
{A_3} = \left( {1:0:0} \right)\\
{B_3} = \left( {{x_3}:1:{y_3}} \right)\,\\
{C_3} = \left( {0:1:0} \right)\\
{D_3} = \left( {{z_3}{t_3}:1:{z_3}} \right)\,\\
{E_3} = \left( {1:1:1} \right)\\
{F_3} = \left( {0:0:1} \right)\\
{y_3},{z_3},{t_3} \sim 0 \ne {x_3}
\end{array}&\vline&  * 
\end{array}\]

The relations between the first and the second charts are:
\[\begin{array}{l}
{x_1} = {\left[ {C,A;D,B} \right]_E} = {\left[ {\left( \begin{array}{c}
0\\
1\\
0
\end{array} \right),\left( \begin{array}{c}
1\\
0\\
0
\end{array} \right);\left( \begin{array}{c}
0\\
0\\
1
\end{array} \right),\left( \begin{array}{c}
1\\
{x_2}\\
{x_2}{y_2}
\end{array} \right)} \right]_{\left( {\scriptstyle1\atop
{\scriptstyle1\atop
\scriptstyle1}} \right)}} = \end{array}\]
\[ = \frac{{\left| {\left( {\begin{array}{*{20}{c}}
0&0&1\\
1&0&1\\
0&1&1
\end{array}} \right)} \right| \cdot \left| {\left( {\begin{array}{*{20}{c}}
1&1&1\\
0&{{x_2}}&1\\
0&{{x_2}{y_2}}&1
\end{array}} \right)} \right|}}{{\left| {\left( {\begin{array}{*{20}{c}}
0&1&1\\
1&{{x_2}}&1\\
0&{{x_2}{y_2}}&1
\end{array}} \right)} \right| \cdot \left| {\left( {\begin{array}{*{20}{c}}
1&0&1\\
0&0&1\\
0&1&1
\end{array}} \right)} \right|}} = 
\frac{{{x_2}\left( {1 - {y_2}} \right)}}{1 - {x_2}{y_2} } \]

\[\begin{array}{l}
{y_1} = {\left[ {C,E;D,B} \right]_A} = {\left[ {\left( \begin{array}{c}
0\\
1\\
0
\end{array} \right),\left( \begin{array}{l}
1\\
1\\
1
\end{array} \right);\left( \begin{array}{c}
0\\
0\\
1
\end{array} \right),\left( \begin{array}{c}
1\\
{x_2}\\
{x_2}{y_2}
\end{array} \right)} \right]_{\left( {\scriptstyle1\atop
{\scriptstyle0\atop
\scriptstyle0}} \right)}}\,\, = \end{array} \]\[
 = \frac{{\left| {\left( {\begin{array}{*{20}{c}}
1&0\\
0&1
\end{array}} \right)} \right| \cdot \left| {\left( {\begin{array}{*{20}{c}}
1&{{x_2}}\\
1&{{x_2}{y_2}}
\end{array}} \right)} \right|}}{{\left| {\left( {\begin{array}{*{20}{c}}
1&{{x_2}}\\
0&{{x_2}{y_2}}
\end{array}} \right)} \right| \cdot \left| {\left( {\begin{array}{*{20}{c}}
1&0\\
1&1
\end{array}} \right)} \right|}} = \frac{{{x_2}\left( {{y_2} - 1} \right)}}{{{x_2}{y_2}}} 
= \frac{{{y_2} - 1}}{{{y_2}}}\]

\[\begin{array}{l}
{z_1} = {\left[ {E,A;F,B} \right]_C} = {\left[ {\left( \begin{array}{l}
1\\
1\\
1
\end{array} \right),\left( \begin{array}{c}
1\\
0\\
0
\end{array} \right);\left( \begin{array}{c}
{z_2}\\
{1+t_2}\\
1
\end{array} \right),\left( \begin{array}{c}
1\\
{x_2}\\
{x_2}{y_2}
\end{array} \right)} \right]_{\left( {\scriptstyle0\atop
{\scriptstyle1\atop
\scriptstyle0}} \right)}} = \end{array}\]\[
\, = \frac{{\left| {\left( {\begin{array}{*{20}{c}}
1&{{z_2}}\\
1&1
\end{array}} \right)} \right| \cdot \left| {\left( {\begin{array}{*{20}{c}}
1&1\\
0&{{x_2}{y_2}}
\end{array}} \right)} \right|}}{{\left| {\left( {\begin{array}{*{20}{c}}
1&1\\
1&{{x_2}{y_2}}
\end{array}} \right)} \right| \cdot \left| {\left( {\begin{array}{*{20}{c}}
1&{{z_2}}\\
0&1
\end{array}} \right)} \right|}} = \frac{{\left( {1 - {z_2}} \right){x_2}{y_2}}}{{{x_2}{y_2} - 1}}\]

\[\begin{array}{l}
{t_1} = {\left[ {A,C;F,B} \right]_E} = \,{\left[ {\left( \begin{array}{c}
1\\
0\\
0
\end{array} \right),\left( \begin{array}{c}
0\\
1\\
0
\end{array} \right);\left( \begin{array}{c}
{z_2}\\
1 + {t_2}\\
1
\end{array} \right),\left( \begin{array}{c}
1\\
{x_2}\\
{x_2}{y_2}
\end{array} \right)} \right]_{\left( {\scriptstyle1\atop
{\scriptstyle1\atop
\scriptstyle1}} \right)}} = \end{array}\]
\[ = \frac{{\left| {\left( {\begin{array}{*{20}{c}}
1&{{z_2}}&1\\
0&{1 + {t_2}}&1\\
0&1&1
\end{array}} \right)} \right| \cdot \left| {\left( {\begin{array}{*{20}{c}}
0&1&1\\
1&{{x_2}}&1\\
0&{{x_2}{y_2}}&1
\end{array}} \right)} \right|}}{{\left| {\left( {\begin{array}{*{20}{c}}
1&1&1\\
0&{{x_2}}&1\\
0&{{x_2}{y_2}}&1
\end{array}} \right)} \right| \cdot \left| {\left( {\begin{array}{*{20}{c}}
0&{{z_2}}&1\\
1&{1 + {t_2}}&1\\
0&1&1
\end{array}} \right)} \right|}} = \frac{{ - {t_2}\left( {1 - {x_2}{y_2}} \right)}}{{{x_2}\left( {1 - {y_2}} \right)\left( {1 - {z_2}} \right)}}\]

Three of the relations between the first chart and the third chart are:
\[\begin{array}{l}
{x_1} = {\left[ {C,A;D,B} \right]_E} = {\left[ {\left( \begin{array}{c}
0\\
1\\
0
\end{array} \right),\left( \begin{array}{c}
1\\
0\\
0
\end{array} \right);\left( \begin{array}{c}
{z_3}{t_3}\\
1\\
{z_3}
\end{array} \right),\left( \begin{array}{c}
{x_3}\\
1\\
{y_3}
\end{array} \right)} \right]_{\left( {\scriptstyle1\atop
{\scriptstyle1\atop
\scriptstyle1}} \right)}} = \end{array}\]\[
 = \frac{{\left| {\left( {\begin{array}{*{20}{c}}
0&{{z_3}{t_3}}&1\\
1&1&1\\
0&{{z_3}}&1
\end{array}} \right)} \right| \cdot \left| {\left( {\begin{array}{*{20}{c}}
1&{{x_3}}&1\\
0&1&1\\
0&{{y_3}}&1
\end{array}} \right)} \right|}}{{\left| {\left( {\begin{array}{*{20}{c}}
0&{{x_3}}&1\\
1&1&1\\
0&{{y_3}}&1
\end{array}} \right)} \right| \cdot \left| {\left( {\begin{array}{*{20}{c}}
1&{{z_3}{t_3}}&1\\
0&1&1\\
0&{{z_3}}&1
\end{array}} \right)} \right|}} = \frac{{{z_3}\left( {1 - {t_3}} \right)\left( {1 - {y_3}} \right)}}{{\left( {{y_3} - {x_3}} \right)\left( {1 - {z_3}} \right)}}\]

\[\begin{array}{l}
{y_1} = {\left[ {C,E;D,B} \right]_A} = {\left[ {\left( \begin{array}{c}
0\\
1\\
0
\end{array} \right),\left( \begin{array}{l}
1\\
1\\
1
\end{array} \right);\left( \begin{array}{c}
{z_3}{t_3}\\
1\\
{z_3}
\end{array} \right),\left( \begin{array}{c}
{x_3}\\
1\\
{y_3}
\end{array} \right)} \right]_{\left( {\scriptstyle1\atop
{\scriptstyle0\atop
\scriptstyle0}} \right)}}\,\, = \end{array}\]\[
 = \frac{{\left| {\left( {\begin{array}{*{20}{c}}
1&1\\
0&{{z_3}}
\end{array}} \right)} \right| \cdot \left| {\left( {\begin{array}{*{20}{c}}
1&{{x_3}}\\
1&{{y_3}}
\end{array}} \right)} \right|}}{{\left| {\left( {\begin{array}{*{20}{c}}
1&{{x_3}}\\
0&{{y_3}}
\end{array}} \right)} \right| \cdot \left| {\left( {\begin{array}{*{20}{c}}
1&1\\
1&{{z_3}}
\end{array}} \right)} \right|}} = \frac{{{z_3}\left( {{y_3} - {x_3}} \right)}}{{{y_3}\left( {{z_3} - 1} \right)}}\]

\[\begin{array}{l}
{t_1} = {\left[ {A,C;F,B} \right]_E} = \,{\left[ {\left( \begin{array}{c}
1\\
0\\
0
\end{array} \right),\left( \begin{array}{c}
0\\
1\\
0
\end{array} \right);\left( \begin{array}{c}
0\\
0\\
1
\end{array} \right),\left( \begin{array}{c}
{x_3}\\
1\\
{y_3}
\end{array} \right)} \right]_{\left( {\scriptstyle1\atop
{\scriptstyle1\atop
\scriptstyle1}} \right)}} = \, \end{array} \]
\[ = \frac{{\left| {\left( {\begin{array}{*{20}{c}}
1&0&1\\
0&0&1\\
0&1&1
\end{array}} \right)} \right| \cdot \left| {\left( {\begin{array}{*{20}{c}}
0&{{x_3}}&1\\
1&1&1\\
0&{{y_3}}&1
\end{array}} \right)} \right|}}{{\left| {\left( {\begin{array}{*{20}{c}}
1&{{x_3}}&1\\
0&1&1\\
0&{{y_3}}&1
\end{array}} \right)} \right| \cdot \left| {\left( {\begin{array}{*{20}{c}}
0&0&1\\
1&0&1\\
0&1&1
\end{array}} \right)} \right|}} = \frac{{{x_3} - {y_3}}}{{1 - {y_3}}}\]
We add one more relation between the second and the third chart:
\[{z_2} = {\left[ {D,A;F,E} \right]_C} = {\left[ {\left( \begin{array}{c}
{z_3}{t_3}\\
1\\
{z_3}
\end{array} \right),\left( \begin{array}{c}
1\\
0\\
0
\end{array} \right);\left( \begin{array}{c}
0\\
0\\
1
\end{array} \right),\left( \begin{array}{c}
1\\
1\\
1
\end{array} \right)} \right]_{\left( {\scriptstyle0\atop
{\scriptstyle1\atop
\scriptstyle0}} \right)}} = \]
\[ = \frac{{\left| {\left( {\begin{array}{*{20}{c}}
{{z_3}{t_3}}&0\\
{{z_3}}&1
\end{array}} \right)} \right| \cdot \left| {\left( {\begin{array}{*{20}{c}}
1&1\\
0&1
\end{array}} \right)} \right|}}{{\left| {\left( {\begin{array}{*{20}{c}}
{{z_3}{t_3}}&1\\
{{z_3}}&1
\end{array}} \right)} \right| \cdot \left| {\left( {\begin{array}{*{20}{c}}
1&0\\
0&1
\end{array}} \right)} \right|}} = \frac{{{z_3} \cdot \left| {\left( {\begin{array}{*{20}{c}}
{{t_3}}&0\\
1&1
\end{array}} \right)} \right|}}{{{z_3} \cdot \left| {\left( {\begin{array}{*{20}{c}}
{{t_3}}&1\\
1&1
\end{array}} \right)} \right|}} = \frac{{{t_3}}}{{{t_3} - 1}}\]

We multiply all formulas by their denominators, and omit the negligible terms, we get:
\[x_1 \approx  {x_2} \left( {1 - {y_2}} \right)\]
\[{y_1}{y_2} = {y_2} - 1\]
\[-{z_1} \approx {x_2}{y_2}\]
\[t_1(1-y_2) \approx  - t_2\]
\[ - {x_1}{x_3} \approx {z_3}\]
\[-y_1 y_3 \approx  - {x_3}{z_3}\]
\[(1-y_3)t_1 = {x_3} - {y_3}\]
\[-z_2 \approx t_3\]

It is easy to see that $y_1,y_2,x_3,t_1\neq 1$.
So on the level of differentials, 

\[dx_1 = (1-y_2) dx_2\]
\[0=(1-y_1)dy_2-y_2dy_1\]
\[-dz_1 = y_2dx_2\]
\[(1-y_2)dt_1-t_1dy_2 =  - dt_2\]
\[x_3 dx_1 + dz_3=0\]
\[x_3 dz_3 - y_1 dy_3=0\]
\[dt_1+(1-t_1)dy_3=dx_3\]
\[-dz_2 = dt_3\]

Using these 8 equation, we may exclude 8 out of 12 given differentials.
First we use the last equation to exclude $dt_3$. 
Then we use the seventh equation to exclude $dx_3$.
After that, we use equations 6 and 5 to exclude $dy_3$ and $dz_3$,
so we get rid of all coordinates of the third chart.
At last we may use the first 4 equations to get rid of all the coordinates 
of the first chart. 
Hence the cotangent space is spanned by $dx_2, dy_2, dz_2$ and $dt_2$.

\subsection{$\mathcal{A}.4$} 

\begin{center}
	\includegraphics{Case_A4.png}
\end{center}

We choose the coordinates for the charts:
\[\small
\begin{array}{*{20}{c}}
\begin{array}{l}
{A_1} = \,\left( {1:0:0} \right)\\
{B_1} = \,\left( {0:1:0} \right)\,\\
{C_1} = \left( {1:1:1} \right)\\
{D_1} = \left( {{x_1}:{y_1}:1} \right)\,\\
{E_1} = \left( {0:0:1} \right)\\
{F_1} = \left( {{z_1}:{z_1}{t_1}:1} \right)
\end{array}&\vline&  * &\vline& \begin{array}{l}
{A_3} = \,\,\left( {1:0:0} \right)\\
{B_3} = \,\,\left( {1:{x_3}:{y_3}} \right)\\
{C_3} = \,\,\left( {0:1:0} \right)\\
{D_3} = \,\,\left( {{z_3}:1:{z_3}{t_3}} \right)\\
{E_3} = \,\,\left( {0:0:1} \right)\\
{F_3} = \,\,\left( {1:1:1} \right)
\end{array}&\vline& \begin{array}{l}
{A_4} = \,\,\,\left( {1:1 + {x_4}:{y_4}} \right)\\
{B_4} = \,\,\left( {1:1 + {x_4} + {z_4}:} \right.\\
\left. {\,\,\,\,\,\,\,\,\,\,\,\,\,\,\,\,\,\,\,\,\,\,\,\,\,\,\,{y_4} + {t_4}{z_4}} \right)\\
{C_4} = \,\,\left( {1:0:0} \right)\\
{D_4} = \,\,\left( {0:1:0} \right)\\
{E_4} = \,\,\left( {0:0:1} \right)\\
{F_4} = \,\,\left( {1:1:1} \right)
\end{array}
\end{array}\]
We don't need the coordinates on the second chart, since all things related to the second chart 
follow easily from the other charts. On the fourth charts we have chosen coordinates in such a way, 
that small changes of these coordinates would cause small changes in the location of points and in 
the direction of $\ell_{1,2}$.
On the first chart $x_1\sim y_1 \neq 0$ since $D$ is near $CE$, 
also $z_1 \sim 0$, $t_1\neq 0$ since $F$ is near $E$ in generic direction.
On the third chart $x_3\neq 0$, $y_3,z_3,t_3\sim 0$.
On the fourth chart $x_4,y_4,z_4\sim 0$ and $t_4$ may be anything.

Here are some formulas relating the coordinates of the first chart and the third chart:
\[\begin{array}{l}
{x_1} = {\left[ {E,A;D,C} \right]_B} = {\left[ {\left( \begin{array}{l}
0\\
0\\
1
\end{array} \right),\left( \begin{array}{l}
1\\
0\\
0
\end{array} \right);\left( \begin{array}{c}
{z_3}\\
1\\
{z_3}{t_3}
\end{array} \right),\left( \begin{array}{l}
0\\
1\\
0
\end{array} \right)} \right]_{\left( {\scriptstyle1\atop
{\scriptstyle{x_3}\atop
\scriptstyle{y_3}}} \right)}} = \end{array}\] 
\[ = \frac{{\left| {\left( {\begin{array}{*{20}{c}}
0&{{z_3}}&1\\
0&1&{{x_3}}\\
1&{{z_3}{t_3}}&{{y_3}}
\end{array}} \right)} \right| \cdot \left| {\left( {\begin{array}{*{20}{c}}
1&0&1\\
0&1&{{x_3}}\\
0&0&{{y_3}}
\end{array}} \right)} \right|}}{{\left| {\left( {\begin{array}{*{20}{c}}
0&0&1\\
0&1&{{x_3}}\\
1&0&{{y_3}}
\end{array}} \right)} \right| \cdot \left| {\left( {\begin{array}{*{20}{c}}
1&{{z_3}}&1\\
0&1&{{x_3}}\\
0&{{z_3}{t_3}}&{{y_3}}
\end{array}} \right)} \right|}} = \frac{{\left( {1 - {x_3}{z_3}} \right){y_3}}}{{{y_3} - {x_3}{z_3}{t_3}}}\]

\[\begin{array}{l}
{y_1} = {\left[ {E,B;D,C} \right]_A} = {\left[ {\left( \begin{array}{l}
0\\
0\\
1
\end{array} \right),\left( \begin{array}{c}
1\\
{x_3}\\
{y_3}
\end{array} \right);\left( \begin{array}{c}
{z_3}\\
1\\
{z_3}{t_3}
\end{array} \right),\left( \begin{array}{l}
0\\
1\\
0
\end{array} \right)} \right]_{\left( {\scriptstyle1\atop
{\scriptstyle0\atop
\scriptstyle0}} \right)}} = \end{array}\]
\[ = \frac{{\left| {\left( {\begin{array}{*{20}{c}}
0&1\\
1&{{z_3}{t_3}}
\end{array}} \right)} \right| \cdot \left| {\left( {\begin{array}{*{20}{c}}
{{x_3}}&1\\
{{y_3}}&0
\end{array}} \right)} \right|}}{{\left| {\left( {\begin{array}{*{20}{c}}
0&1\\
1&0
\end{array}} \right)} \right| \cdot \left| {\left( {\begin{array}{*{20}{c}}
{{x_3}}&1\\
{{y_3}}&{{z_3}{t_3}}
\end{array}} \right)} \right|}} = \frac{{{y_3}}}{{{y_3} - {x_3}{z_3}{t_3}}}\]

\[\begin{array}{l}
{z_1} = {\left[ {E,A;F,C} \right]_B} = {\left[ {\left( \begin{array}{l}
0\\
0\\
1
\end{array} \right),\left( \begin{array}{l}
1\\
0\\
0
\end{array} \right);\left( \begin{array}{l}
1\\
1\\
1
\end{array} \right),\left( \begin{array}{l}
0\\
1\\
0
\end{array} \right)} \right]_
{\left( {\scriptstyle1\atop
{\scriptstyle x_3\atop
\scriptstyle y_3}} \right)}} = \end{array}\]
\[ = \frac{{\left| {\left( {\begin{array}{*{20}{c}}
0&1&1\\
0&1&{{x_3}}\\
1&1&{{y_3}}
\end{array}} \right)} \right| \cdot \left| {\left( {\begin{array}{*{20}{c}}
1&0&1\\
0&1&{{x_3}}\\
0&0&{{y_3}}
\end{array}} \right)} \right|}}{{\left| {\left( {\begin{array}{*{20}{c}}
0&0&1\\
0&1&{{x_3}}\\
1&0&{{y_3}}
\end{array}} \right)} \right| \cdot \left| {\left( {\begin{array}{*{20}{c}}
1&1&1\\
0&1&{{x_3}}\\
0&1&{{y_3}}
\end{array}} \right)} \right|}} = \frac{{\left( {1 - {x_3}} \right){y_3}}}{{{y_3} - {x_3}}}\]

\[\begin{array}{l}
{t_1} = {\left[ {A,B;F,C} \right]_E} = {\left[ {\left( \begin{array}{l}
1\\
0\\
0
\end{array} \right),\left( \begin{array}{c}
1\\
{x_3}\\
{y_3}
\end{array} \right);\left( \begin{array}{l}
1\\
1\\
1
\end{array} \right),\left( \begin{array}{l}
0\\
1\\
0
\end{array} \right)} \right]_{\left( {\scriptstyle0\atop
{\scriptstyle0\atop
\scriptstyle1}} \right)}} = \end{array}\]\[
 = \frac{{\left| {\left( {\begin{array}{*{20}{c}}
1&1\\
0&1
\end{array}} \right)} \right| \cdot \left| {\left( {\begin{array}{*{20}{c}}
1&0\\
{{x_3}}&1
\end{array}} \right)} \right|}}{{\left| {\left( {\begin{array}{*{20}{c}}
1&0\\
0&1
\end{array}} \right)} \right| \cdot \left| {\left( {\begin{array}{*{20}{c}}
1&1\\
{{x_3}}&1
\end{array}} \right)} \right|}} = \frac{1}{{1 - {x_3}}}\]

We may write similar formulas relating the third chart and the fourth chart:
\[\begin{array}{l}
{x_3} = {\left[ {A,C;B,F} \right]_E} = {\left[\footnotesize {\left( \begin{array}{c}
1\\
1 + {x_4}\\
{y_4}
\end{array} \right),\left( \begin{array}{l}
1\\
0\\
0
\end{array} \right);\left( \begin{array}{c}
1\\
1 + {x_4} + {z_4}\\
{y_4} + {t_4}{z_4}
\end{array} \right),\left( \begin{array}{l}
1\\
1\\
1
\end{array} \right)} \right]_{\left( {\scriptstyle0\atop
{\scriptstyle0\atop
\scriptstyle1}} \right)}} = \end{array}\]
\[ = \frac{{\left| {\left( {\begin{array}{*{20}{c}}
1&1\\
{1 + {x_4}}&{1 + {x_4} + {z_4}}
\end{array}} \right)} \right| \cdot \left| {\left( {\begin{array}{*{20}{c}}
1&1\\
0&1
\end{array}} \right)} \right|}}{{\left| {\left( {\begin{array}{*{20}{c}}
1&1\\
{1 + {x_4}}&1
\end{array}} \right)} \right| \cdot \left| {\left( {\begin{array}{*{20}{c}}
1&1\\
0&{1 + {x_4} + {z_4}}
\end{array}} \right)} \right|}} = \frac{{{z_4}}}{{ - {x_4}\left( {1 + {x_4} + {z_4}} \right)}}
\]

\[\begin{array}{l}
{y_3} = {\left[ {A,E;B,F} \right]_C} = {\left[\footnotesize {\left( \begin{array}{c}
1\\
1 + {x_4}\\
{y_4}
\end{array} \right),\left( \begin{array}{l}
0\\
0\\
1
\end{array} \right);\left( \begin{array}{c}
1\\
1 + {x_4} + {z_4}\\
{y_4} + {t_4}{z_4}
\end{array} \right),\left( \begin{array}{l}
1\\
1\\
1
\end{array} \right)} \right]_{\left( {\scriptstyle1\atop
{\scriptstyle0\atop
\scriptstyle0}} \right)}} = \\
 = \frac{{\left| {\left( {\begin{array}{*{20}{c}}
{1 + {x_4}}&{1 + {x_4} + {z_4}}\\
{{y_4}}&{{y_4} + {t_4}{z_4}}
\end{array}} \right)} \right| \cdot \left| {\left( {\begin{array}{*{20}{c}}
0&1\\
1&1
\end{array}} \right)} \right|}}{{\left| {\left( {\begin{array}{*{20}{c}}
{1 + {x_4}}&1\\
{{y_4}}&1
\end{array}} \right)} \right| \cdot \left| {\left( {\begin{array}{*{20}{c}}
0&{1 + {x_4} + {z_4}}\\
1&{{y_4} + {t_4}{z_4}}
\end{array}} \right)} \right|}} = \end{array}\]
\[ =\frac{{\left| {\left( {\begin{array}{*{20}{c}}
{1 + {x_4}}&{{z_4}}\\
{{y_4}}&{{t_4}{z_4}}
\end{array}} \right)} \right| }}{\left| {\left( {\begin{array}{*{20}{c}}
{1 + {x_4}-y_4}&0\\
{{y_4}}&1
\end{array}} \right)} \right| \cdot  {\left(1 + {x_4} + {z_4} \right) }} 
 = \frac{{{z_4}\left( {{t_4} + {x_4}{t_4} - {y_4}} \right)}}
 {{\left( {1 + {x_4} - {y_4}} \right)\left( {1 + {x_4} + {z_4}} \right)}} \]

\[\begin{array}{l}
{z_3} = {\left[ {C,A;D,F} \right]_E} = {\left[ {\left( \begin{array}{l}
1\\
0\\
0
\end{array} \right),\left( \begin{array}{c}
1\\
1 + {x_4}\\
{y_4}
\end{array} \right);\,\left( \begin{array}{l}
0\\
1\\
0
\end{array} \right),\left( \begin{array}{l}
1\\
1\\
1
\end{array} \right)} \right]_{\left( {\scriptstyle0\atop
{\scriptstyle0\atop
\scriptstyle1}} \right)}} = \end{array}\]
\[ = \frac{{\left| {\left( {\begin{array}{*{20}{c}}
1&0\\
0&1
\end{array}} \right)} \right| \cdot \left| {\left( {\begin{array}{*{20}{c}}
1&1\\
{1+x_4}&1
\end{array}} \right)} \right|}}{{\left| {\left( {\begin{array}{*{20}{c}}
1&1\\
0&1
\end{array}} \right)} \right| \cdot \left| {\left( {\begin{array}{*{20}{c}}
1&0\\
{1 + {x_4}}&1
\end{array}} \right)} \right|}} = -x_4\]

\[\begin{array}{l}
{t_3} = {\left[ {A,E;D,F} \right]_C} = {\left[ {\left( \begin{array}{c}
1\\
1 + {x_4}\\
{y_4}
\end{array} \right),\left( \begin{array}{l}
0\\
0\\
1
\end{array} \right);\,\left( \begin{array}{l}
0\\
1\\
0
\end{array} \right),\left( \begin{array}{l}
1\\
1\\
1
\end{array} \right)} \right]_{\left( {\scriptstyle1\atop
{\scriptstyle0\atop
\scriptstyle0}} \right)}} = \end{array}\]
\[ = \frac{{\left| {\left( {\begin{array}{*{20}{c}}
{1 + {x_4}}&1\\
{{y_4}}&0
\end{array}} \right)} \right| \cdot \left| {\left( {\begin{array}{*{20}{c}}
0&1\\
1&1
\end{array}} \right)} \right|}}{{\left| {\left( {\begin{array}{*{20}{c}}
{1 + {x_4}}&1\\
{{y_4}}&1
\end{array}} \right)} \right| \cdot \left| {\left( {\begin{array}{*{20}{c}}
0&1\\
1&0
\end{array}} \right)} \right|}} = \frac{{{y_4}}}{{{y_4} - {x_4} - 1}}\]

We multiply by the denominators:
\[{x_1}\left( {{y_3} - {x_3}{z_3}{t_3}} \right) = \left( {1 - {x_3}{z_3}} \right){y_3}\]
\[{y_1}\left( {{y_3} - {x_3}{z_3}{t_3}} \right) = {y_3}\]
\[{z_1}\left( {{y_3} - {x_3}} \right) = \left( {1 - {x_3}} \right){y_3}\]
\[{t_1}\left( {1 - {x_3}} \right) = 1\]
\[{x_3}{x_4}\left( {1 + {x_4} + {z_4}} \right) =  - {z_4}\]
\[{y_3}\left( {1 + {x_4} - {y_4}} \right)\left( {1 + {x_4} + {z_4}} \right) = {z_4}\left( {{t_4} + {x_4}{t_4} - {y_4}} \right)\]
\[z_3 = -x_4\]
\[{t_3}\left( {{y_4} - {x_4} - 1} \right) = {y_4}\]

We may disregard terms of second order:
since $z_1, y_3,z_3,t_3, x_4,y_4,z_4 \sim 0$, the product of any 2 of those is negligible.

\[ 0 \approx  y_3 (1-x_1)\]
\[-y_1 y_3 \approx y_3\]
\[- z_1 x_3 \approx \left( {1 - x_3} \right)y_3\]
\[{t_1}\left( {1 - {x_3}} \right) = 1\]
\[ x_3x_4 \approx - z_4 \]
\[ y_3  \approx z_4 t_4 \]
\[ z_3 = -x_4\]
\[-t_3 \approx y_4\]

Notice that $D$ and $C$ on the first chart don't coincide, so $x_1 \neq 1 $, hence by the first of the last 8 equations,
$y_3 \approx  0$. This allows to simplify several of the other equations. 
The second equation becomes unimportant. The third equation becomes $z_1 x_3 \approx 0$, 
but since $z_1 \sim 0$ and $x_3 \neq 0$, it means $z_1\approx 0$. 
Therefore:
\[dy_3 = 0 = dz_1\]
\[ (1-x_3)dt_1 - t_1 dx_3= 0\]
\[ x_3\cdot dx_4 = - dz_4 \]
\[ 0  =  t_4\cdot dz_4 \]
\[ dz_3 = -dx_4\]
\[ dt_3 = -dy_4\]

With these equations, we may express $dx_3,dy_3,dz_3, dt_3$ using the coordinates from charts one and four,
and remain with 7 generators for the cotangent space (since $dz_1 =0$) and a relation
\[ 0  =  t_4\cdot dz_4 \]
Clearly, we need more relations. We may relate first and fourth charts directly:

\[\begin{array}{l}
{x_1} = {\left[ {E,A;D,C} \right]_B} = {\left[ {\left( \begin{array}{l}
0\\
0\\
1
\end{array} \right),\left( \begin{array}{c}
1\\
1 + {x_4}\\
{y_4}
\end{array} \right);\left( \begin{array}{l}
0\\
1\\
0
\end{array} \right),\left( \begin{array}{l}
1\\
0\\
0
\end{array} \right)} \right]_{\left( {\scriptstyle1\atop
{\scriptstyle1 + {x_4} + {z_4}\atop
\scriptstyle{y_4} + {t_4}{z_4}}} \right)}} = \\
 = {\scriptsize\frac{{\left| {\left( {\begin{array}{*{20}{c}}
0&0&1\\
0&1&{1 + {x_4} + {z_4}}\\
1&0&{{y_4} + {t_4}{z_4}}
\end{array}} \right)} \right| \cdot \left| {\left( {\begin{array}{*{20}{c}}
1&1&1\\
{1 + {x_4}}&0&{1 + {x_4} + {z_4}}\\
{{y_4}}&0&{{y_4} + {t_4}{z_4}}
\end{array}} \right)} \right|}}{{\left| {\left( {\begin{array}{*{20}{c}}
0&1&1\\
0&0&{1 + {x_4} + {z_4}}\\
1&0&{{y_4} + {t_4}{z_4}}
\end{array}} \right)} \right| \cdot \left| {\left( {\begin{array}{*{20}{c}}
1&0&1\\
{1 + {x_4}}&1&{1 + {x_4} + {z_4}}\\
{{y_4}}&0&{{y_4} + {t_4}{z_4}}
\end{array}} \right)} \right|}}}
 = \frac{{\left| {\left( {\begin{array}{*{20}{c}}
{1 + {x_4}}&{{z_4}}\\
{{y_4}}&{{t_4}{z_4}}
\end{array}} \right)} \right|}}{{\left( {1 + {x_4} + {z_4}} \right)\left| {\left( {\begin{array}{*{20}{c}}
1&0\\
{{y_4}}&{{t_4}{z_4}}
\end{array}} \right)} \right|}} = \end{array}\]
\[ = \frac{{\left( {1 + {x_4}} \right){t_4}{z_4} - {y_4}{z_4}}}{{\left( {1 + {x_4} + {z_4}} \right){t_4}{z_4}}} = \frac{{\left( {1 + {x_4}} \right){t_4} - {y_4}}}{{\left( {1 + {x_4} + {z_4}} \right){t_4}}}\]

\[
{y_1} = {\left[ {E,B;D,C} \right]_A} = {\left[\small {\left( \begin{array}{l}
0\\
0\\
1
\end{array} \right),\left( \begin{array}{c}
1\\
1 + {x_4} + {z_4}\\
{y_4} + {t_4}{z_4}
\end{array} \right);\left( \begin{array}{l}
0\\
1\\
0
\end{array} \right),\left( \begin{array}{l}
1\\
0\\
0
\end{array} \right)} \right]_{\left( {\scriptstyle1\atop
{\scriptstyle1 + {x_4}\atop
\scriptstyle{y_4}}} \right)}} = \]

\[ = \frac{{\left| {\left( {\begin{array}{*{20}{c}}
0&0&1\\
0&1&{1 + {x_4}}\\
1&0&{{y_4}}
\end{array}} \right)} \right| \cdot \left| {\left( {\begin{array}{*{20}{c}}
1&1&1\\
{1 + {x_4} + {z_4}}&0&{1 + {x_4}}\\
{{y_4} + {t_4}{z_4}}&0&{{y_4}}
\end{array}} \right)} \right|}}{{\left| {\left( {\begin{array}{*{20}{c}}
0&1&1\\
0&0&{1 + {x_4}}\\
1&0&{{y_4}}
\end{array}} \right)} \right| \cdot \left| {\left( {\begin{array}{*{20}{c}}
1&0&1\\
{1 + {x_4} + {z_4}}&1&{1 + {x_4}}\\
{{y_4} + {t_4}{z_4}}&0&{{y_4}}
\end{array}} \right)} \right|}} = \]

\[ = \frac{{\left| {\left( {\begin{array}{*{20}{c}}
{{z_4}}&{1 + {x_4}}\\
{{t_4}{z_4}}&{{y_4}}
\end{array}} \right)} \right|}}{{\left( {1 + {x_4}} \right) \cdot \left| {\left( {\begin{array}{*{20}{c}}
0&1\\
{{t_4}{z_4}}&{{y_4}}
\end{array}} \right)} \right|}} = \frac{{\left| {\left( {\begin{array}{*{20}{c}}
1&{1 + {x_4}}\\
{{t_4}}&{{y_4}}
\end{array}} \right)} \right|}}{{\left( {1 + {x_4}} \right) \cdot \left| {\left( {\begin{array}{*{20}{c}}
0&1\\
{{t_4}}&{{y_4}}
\end{array}} \right)} \right|}} = \frac{{\left( {1 + {x_4}} \right){t_4} - {y_4}}}{{\left( {1 + {x_4}} \right)\left( t_4 \right)}}
\]

\[\begin{array}{l}
{t_1} = {\left[ \small {A,B;F,C} \right]_E} = {\left[\footnotesize {\left( \begin{array}{c}
1\\
1 + {x_4}\\
{y_4}
\end{array} \right),\left( \begin{array}{c}
1\\
1 + {x_4} + {z_4}\\
{y_4} + {t_4}{z_4}
\end{array} \right);\left( \begin{array}{l}
1\\
1\\
1
\end{array} \right),\left( \begin{array}{l}
1\\
0\\
0
\end{array} \right)} \right]_{\left( {\scriptstyle0\atop
{\scriptstyle0\atop
\scriptstyle1}} \right)}} = \end{array}\]
\[ = \frac{{\left| {\left( {\begin{array}{*{20}{c}}
1&1\\
{1 + {x_4}}&1
\end{array}} \right)} \right| \cdot \left| {\left( {\begin{array}{*{20}{c}}
1&1\\
{1 + {x_4} + {z_4}}&0
\end{array}} \right)} \right|}}{{\left| {\left( {\begin{array}{*{20}{c}}
1&1\\
{1 + {x_4}}&0
\end{array}} \right)} \right| \cdot \left| {\left( {\begin{array}{*{20}{c}}
1&1\\
{1 + {x_4} + {z_4}}&1
\end{array}} \right)} \right|}} = \frac{{{x_4}\left( {1 + {x_4} + {z_4}} \right)}}
{{\left( {1 + {x_4}} \right)\left( {{x_4} + {z_4}} \right)}}\]

Hence 
\[\begin{array}{l}
{x_1}\left( {1 + {x_4} + {z_4}} \right){t_4} = \left( {1 + {x_4}} \right){t_4} - {y_4}\\
{y_1}\left( {1 + {x_4}} \right)  t_4  = \left( {1 + {x_4}} \right){t_4} - {y_4}\\
{t_1}\left( {1 + {x_4}} \right)\left( {{x_4} + {z_4}} \right) = {x_4}\left( {1 + {x_4} + {z_4}} \right)
\end{array}\]

We may use the second equation (among the last triple) to express $dy_1$ from the other coordinates 
(since the coefficient of $dy_1$ will be approximately 1).
Then we may use the first equation to express $dx_1$ using the other coordinates. 
Then we use the last equation to express $dz_4$, since $z_1 \neq 0$.
So, the cotangent space is generated by $dx_1,dy_1, dt_1, dt_4$, so its dimension is not greater than 4. 

\newpage
\subsection{$\mathcal{B}.2$} 
\begin{center}
	\includegraphics{Case_B2.png}
\end{center}
The local coordinates for the first 3 charts are 
\[\begin{array}{*{20}{c}}
\begin{array}{l}
{A_1} = \left( {1:0:0} \right)\\
{B_1} = \left( {1:1:1} \right)\\
{C_1} = \left( {0:1:0} \right)\\
{D_1} = \left( {{x_1}:1:{y_1}} \right)\\
{E_1} = \left( {0:0:1} \right)\\
{F_1} = \left( {{z_1}:{z_1}{t_1}:1} \right)\\
{x_1},{z_1},{t_1} \sim 0 \ne {y_1}
\end{array}&\vline& \begin{array}{l}
{A_2} = \left( {1:0:0} \right)\\
{B_2} = \left( {1:{x_2}:{x_2}{y_2}} \right)\\
{C_2} = \left( {0:1:0} \right)\\
{D_2} = \left( {0:0:1} \right)\\
{E_2} = \left( {1:1:1} \right)\,\\
{F_2} = \left( {{z_2}:1 + {t_2}:1} \right)\\
{x_2},{z_2},{t_2} \sim 0 \ne {y_2}
\end{array}&\vline& \begin{array}{l}
{A_3} = \left( {1:1 + {x_3}:{y_3}} \right)\\
{B_3} = \left( {1:0:0} \right)\,\\
{C_3} = \left( {0:1:0} \right)\\
{D_3} = \left( {{z_3}{t_3}:1:{z_3}} \right)\,\\
{E_3} = \left( {1:1:1} \right)\\
{F_3} = \left( {0:0:1} \right)\\
{x_3},{y_3},{z_3},{t_3} \sim 0
\end{array}
\end{array}\]
We won't write coordinates for the last two charts, since they follow easily from the first 3 charts.
We may relate the coordinates of the second chart to the first chart:
\[\begin{array}{l}
{x_2} = {\left[ {A,C;B,E} \right]_D} = {\left[ {\left( \begin{array}{c}
1\\
0\\
0
\end{array} \right),\left( \begin{array}{c}
0\\
1\\
0
\end{array} \right);\left( \begin{array}{l}
1\\
1\\
1
\end{array} \right),\left( \begin{array}{c}
0\\
0\\
1
\end{array} \right)} \right]_{\left( {\scriptstyle{x_1}\atop
{\scriptstyle1\atop
\scriptstyle{y_1}}} \right)}} = \end{array}\]
\[ = \frac{{\left| {\left( {\begin{array}{*{20}{c}}
1&1&{{x_1}}\\
0&1&1\\
0&1&{{y_1}}
\end{array}} \right)} \right| \cdot \left| {\left( {\begin{array}{*{20}{c}}
0&0&{{x_1}}\\
1&0&1\\
0&1&{{y_1}}
\end{array}} \right)} \right|}}{{\left| {\left( {\begin{array}{*{20}{c}}
1&0&{{x_1}}\\
0&0&1\\
0&1&{{y_1}}
\end{array}} \right)} \right| \cdot \left| {\left( {\begin{array}{*{20}{c}}
0&1&{{x_1}}\\
1&1&1\\
0&1&{{y_1}}
\end{array}} \right)} \right|}} = \frac{{\left( {1 - {y_1}} \right){x_1}}}{{{x_1} - {y_1}}}\]

\[\begin{array}{l}
{y_2} = {\left[ {C,D;B,E} \right]_A} = {\left[ {\left( \begin{array}{c}
0\\
1\\
0
\end{array} \right),\left( \begin{array}{c}
{x_1}\\
1\\
{y_1}
\end{array} \right);\left( \begin{array}{l}
1\\
1\\
1
\end{array} \right),\left( \begin{array}{c}
0\\
0\\
1
\end{array} \right)} \right]_{\left( {\scriptstyle1\atop
{\scriptstyle0\atop
\scriptstyle0}} \right)}} = \end{array}\]
\[ = \frac{{\left| {\left( {\begin{array}{*{20}{c}}
1&1\\
0&1
\end{array}} \right)} \right| \cdot \left| {\left( {\begin{array}{*{20}{c}}
1&0\\
{{y_1}}&1
\end{array}} \right)} \right|}}{{\left| {\left( {\begin{array}{*{20}{c}}
1&0\\
0&1
\end{array}} \right)} \right| \cdot \left| {\left( {\begin{array}{*{20}{c}}
1&1\\
{{y_1}}&1
\end{array}} \right)} \right|}} = \frac{1}{{1 - {y_1}}}\]

\[\begin{array}{l}
{z_2} = {\left[ {D,A;F,E} \right]_C} = {\left[ {\left( \begin{array}{c}
{x_1}\\
1\\
{y_1}
\end{array} \right),\left( \begin{array}{c}
1\\
0\\
0
\end{array} \right);\left( \begin{array}{c}
{z_1}\\
{z_1}{t_1}\\
1
\end{array} \right),\left( \begin{array}{c}
0\\
0\\
1
\end{array} \right)} \right]_{\left( {\scriptstyle0\atop
{\scriptstyle1\atop
\scriptstyle0}} \right)}} = \end{array}\]
\[ = \frac{{\left| {\left( {\begin{array}{*{20}{c}}
{{x_1}}&{{z_1}}\\
{{y_1}}&1
\end{array}} \right)} \right| \cdot \left| {\left( {\begin{array}{*{20}{c}}
1&0\\
0&1
\end{array}} \right)} \right|}}{{\left| {\left( {\begin{array}{*{20}{c}}
{{x_1}}&0\\
{{y_1}}&1
\end{array}} \right)} \right| \cdot \left| {\left( {\begin{array}{*{20}{c}}
1&{{z_1}}\\
0&1
\end{array}} \right)} \right|}} = \frac{{{x_1} - {y_1}{z_1}}}{{{x_1}}}\]

\[\begin{array}{l}
1 + {t_2} = {\left[ {D,C;F,E} \right]_A} = {\left[ {\left( \begin{array}{c}
{x_1}\\
1\\
{y_1}
\end{array} \right),\left( \begin{array}{c}
0\\
1\\
0
\end{array} \right);\left( \begin{array}{c}
{z_1}\\
{z_1}{t_1}\\
1
\end{array} \right),\left( \begin{array}{c}
0\\
0\\
1
\end{array} \right)} \right]_{\left( {\scriptstyle1\atop
{\scriptstyle0\atop
\scriptstyle0}} \right)}} = \end{array} \]\[
 = \frac{{\left| {\left( {\begin{array}{*{20}{c}}
1&{{z_1}{t_1}}\\
{{y_1}}&1
\end{array}} \right)} \right| \cdot \left| {\left( {\begin{array}{*{20}{c}}
1&0\\
0&1
\end{array}} \right)} \right|}}{{\left| {\left( {\begin{array}{*{20}{c}}
1&0\\
{{y_1}}&1
\end{array}} \right)} \right| \cdot \left| {\left( {\begin{array}{*{20}{c}}
1&{{z_1}{t_1}}\\
0&1
\end{array}} \right)} \right|}} = 1 - {y_1}{z_1}{t_1}\]

Similarly, we may relate second chart to the third chart:
\[\begin{array}{l}
{x_2} = {\left[ {A,C;B,E} \right]_D} = {\left[ {\left( \begin{array}{c}
1\\
1 + {x_3}\\
{y_3}
\end{array} \right),\left( \begin{array}{c}
0\\
1\\
0
\end{array} \right);\left( \begin{array}{c}
1\\
0\\
0
\end{array} \right),\left( \begin{array}{l}
1\\
1\\
1
\end{array} \right)} \right]_{\left( {\scriptstyle{z_3}{t_3}\atop
{\scriptstyle1\atop
\scriptstyle{z_3}}} \right)}} = \end{array}\]
\[ = \frac{{\left| {\left( {\begin{array}{*{20}{c}}
1&1&{{z_3}{t_3}}\\
{1 + {x_3}}&0&1\\
{{y_3}}&0&{{z_3}}
\end{array}} \right)} \right| \cdot \left| {\left( {\begin{array}{*{20}{c}}
0&1&{{z_3}{t_3}}\\
1&1&1\\
0&1&{{z_3}}
\end{array}} \right)} \right|}}{{\left| {\left( {\begin{array}{*{20}{c}}
1&1&{{z_3}{t_3}}\\
{1 + {x_3}}&1&1\\
{{y_3}}&1&{{z_3}}
\end{array}} \right)} \right| \cdot \left| {\left( {\begin{array}{*{20}{c}}
0&1&{{z_3}{t_3}}\\
1&0&1\\
0&0&{{z_3}}
\end{array}} \right)} \right|}} = \]\[
 = \frac{{\left( {{z_3} + {z_3}{x_3} - {y_3}} \right)\left( {{t_3} - 1} \right)}}
 {{ {{z_3} - 1 + {y_3} - {y_3}{z_3}{t_3} + \left( {1 + {x_3}} \right){z_3}\left( {{t_3} - 1} \right)} }}\]

\[\begin{array}{l}
{y_2} = {\left[ {C,D;B,E} \right]_A} = {\left[ {\left( \begin{array}{c}
0\\
1\\
0
\end{array} \right),\left( \begin{array}{c}
{z_3}{t_3}\\
1\\
{z_3}
\end{array} \right);\left( \begin{array}{c}
1\\
0\\
0
\end{array} \right),\left( \begin{array}{l}
1\\
1\\
1
\end{array} \right)} \right]_{\left( {\scriptstyle1\atop
{\scriptstyle1 + {x_3}\atop
\scriptstyle{y_3}}} \right)}} = \end{array}\]
\[= \frac{{\left| {\left( {\begin{array}{*{20}{c}}
0&1&1\\
1&0&{1 + {x_3}}\\
0&0&{{y_3}}
\end{array}} \right)} \right| \cdot \left| {\left( {\begin{array}{*{20}{c}}
{{z_3}{t_3}}&1&1\\
1&1&{1 + {x_3}}\\
{{z_3}}&1&{{y_3}}
\end{array}} \right)} \right|}}{{\left| {\left( {\begin{array}{*{20}{c}}
0&1&1\\
1&1&{1 + {x_3}}\\
0&1&{{y_3}}
\end{array}} \right)} \right| \cdot \left| {\left( {\begin{array}{*{20}{c}}
{{z_3}{t_3}}&1&1\\
1&0&{1 + {x_3}}\\
{{z_3}}&0&{{y_3}}
\end{array}} \right)} \right|}} \]\[
 = \frac{{{y_3}\left( {{z_3} - 1 + {y_3} - {y_3}{z_3}{t_3} + \left( {1 + {x_3}} \right){z_3}\left( {{t_3} - 1} \right)} \right)}}{{\left( {1 - {y_3}} \right)\left( {{z_3} + {x_3}{z_3} - {y_3}} \right)}}\]

\[\begin{array}{l}
{z_2} = {\left[ {D,A;F,E} \right]_C} = {\left[ {\left( \begin{array}{c}
{z_3}{t_3}\\
1\\
{z_3}
\end{array} \right),\left( \begin{array}{c}
1\\
1 + {x_3}\\
{y_3}
\end{array} \right);\left( \begin{array}{c}
0\\
0\\
1
\end{array} \right),\left( \begin{array}{l}
1\\
1\\
1
\end{array} \right)} \right]_{\left( {\scriptstyle0\atop
{\scriptstyle1\atop
\scriptstyle0}} \right)}} = \end{array} \]
\[ = \frac{{\left| {\left( {\begin{array}{*{20}{c}}
{{z_3}{t_3}}&0\\
{{z_3}}&1
\end{array}} \right)} \right| \cdot \left| {\left( {\begin{array}{*{20}{c}}
1&1\\
{{y_3}}&1
\end{array}} \right)} \right|}}{{\left| {\left( {\begin{array}{*{20}{c}}
{{z_3}{t_3}}&1\\
{{z_3}}&1
\end{array}} \right)} \right| \cdot \left| {\left( {\begin{array}{*{20}{c}}
1&0\\
{{y_3}}&1
\end{array}} \right)} \right|}} = \frac{{{t_3}\left( {1 - {y_3}} \right)}}{{ {{t_3} - 1} }}\]

It turns out the equation for $t_2$ in terms of the third charts coordinates is long and after a lot of cancellation doesn't produce 
an additional dependency in the cotangent space, so we won't write it down. So we shall add one more relation between the first and the third charts.
\[\begin{array}{l}
{t_1} = {\left[ {A,C;F,B} \right]_E} = {\left[ {\left( \begin{array}{c}
1\\
1 + {x_3}\\
{y_3}
\end{array} \right),\left( \begin{array}{c}
0\\
1\\
0
\end{array} \right);\left( \begin{array}{c}
0\\
0\\
1
\end{array} \right),\left( \begin{array}{c}
1\\
0\\
0
\end{array} \right)} \right]_{\left( {\scriptstyle1\atop
{\scriptstyle1\atop
\scriptstyle1}} \right)}} = \\
 = \frac{{\left| {\left( {\begin{array}{*{20}{c}}
1&0&1\\
{1 + {x_3}}&0&1\\
{{y_3}}&1&1
\end{array}} \right)} \right| \cdot \left| {\left( {\begin{array}{*{20}{c}}
0&1&1\\
1&0&1\\
0&0&1
\end{array}} \right)} \right|}}{{\left| {\left( {\begin{array}{*{20}{c}}
1&1&1\\
{1 + {x_3}}&0&1\\
{{y_3}}&0&1
\end{array}} \right)} \right| \cdot \left| {\left( {\begin{array}{*{20}{c}}
0&0&1\\
1&0&1\\
0&1&1
\end{array}} \right)} \right|}} = \frac{{{x_3}}}{{1 + {x_3} - {y_3}}}
\end{array}\]

After multiplying by the denominators and omitting the negligible terms we get:
\[ - {y_1}{x_2} \approx \left( {1 - {y_1}} \right){x_1}\]
\[\left( {1 - {y_1}} \right){y_2} = 1\]
\[0 \approx {x_1} - {y_1}{z_1}\]
\[1 + {t_2} \approx 1\]
\[{x_2} \approx {z_3} - {y_3}\]
\[{y_2}\left( {{z_3} - {y_3}} \right) \approx  - 2{y_3}\]
\[ -z_2 = {t_3}\]
\[t_1 \approx x_3 \]

We conclude 
\[{y_1}d{x_2} + \left( {1 - {y_1}} \right)d{x_1} = 0\]
\[ - {y_2}d{y_1} + \left( {1 - {y_1}} \right)d{y_2} = 0\]
\[d{x_1} - {y_1}d{z_1} = 0\]
\[d{t_2} = 0\]
\[d{x_2} + d{y_3} - d{z_3} = 0\]
\[\left( {2 - {y_2}} \right)d{y_3} + {y_2}d{z_3} = 0\]
\[d{t_3} + d{z_2} = 0\]
\[dt_1 - dx_3 = 0\]

Considering the last 8 equation, we may use the third equation to eliminate $dz_1$,
then equations 1, 2, and 8 to eliminate $dx_1,dy_1$ and $dt_1$ respectively.
Then we may use equations 4, 5, and 7 to eliminate $dt_2,dx_2$ and $dz_2$.
Lastly, we may use the remaining equation 6 to eliminate $dz_3$.
So, the cotangent space is generated by $dy_2, dx_3,dy_3,dt_3$. 

\subsection{$\mathcal{B}.4$} 
\begin{center}
	\includegraphics{Case_B4.png}
\end{center}

We don't need coordinates for the second an third charts, 
as all information there can be expressed explicitly from the other charts.
We may choose local coordinates related the other charts.
\[\begin{array}{*{20}{c}}
\begin{array}{l}
{A_1} = \left( {1:0:0} \right)\\
{B_1} = \left( {1:1:1} \right)\\
{C_1} = \left( {0:1:0} \right)\\
{D_1} = \left( {{x_1}:1:{y_1}} \right)\\
{E_1} = \left( {0:0:1} \right)\\
{F_1} = \left( {{z_1}:{z_1}{t_1}:1} \right)\\
{x_1},{z_1},{t_1} \sim 0 \ne {y_1}
\end{array}&\vline& \begin{array}{l}
{A_4} = \left( {1 + {x_4}:1:{y_4}} \right)\\
{B_4} = \left( 1 + {x_4} + {z_4}:\right. \\
\hfill\left. 1:{y_4} + {z_4}{t_4} \right)\\
{C_4} = \left( {1:0:0} \right)\,\\
{D_4} = \left( {0:1:0} \right)\\
{E_4} = \left( {0:0:1} \right)\\
{F_4} = \left( {1:1:1} \right)\\
{x_4},{y_4},{z_4},{t_4} \sim 0
\end{array}&\vline& \begin{array}{l}
{A_5} = \left( {1:1 + {x_5}:{y_5}} \right)\\
{B_5} = \left( {1:0:0} \right)\\
{C_5} = \left( {0:1:0} \right)\\
{D_5} = \left( {{z_5}:1:{z_5}{t_5}} \right)\\
{E_5} = \left( {0:0:1} \right)\\
{F_5} = \left( {1:1:1} \right)\\
{x_5},{y_5},{z_5},{t_5} \sim 0
\end{array}
\end{array}\]

We may connect the third chart to the first chart.

\[\begin{array}{l}
{x_1} = {\left[ {C,A;D,B} \right]_E} = {\left[\footnotesize {\left( \begin{array}{c}
1\\
0\\
0
\end{array} \right),\left( \begin{array}{c}
1 + {x_4}\\
1\\
{y_4}
\end{array} \right);\left( \begin{array}{c}
0\\
1\\
0
\end{array} \right),\left( \begin{array}{c}
1 + {x_4} + {z_4}\\
1\\
{y_4} + {z_4}{t_4}
\end{array} \right)} \right]_{\left( {\scriptstyle0\atop
{\scriptstyle0\atop
\scriptstyle1}} \right)}} = \\
 = \frac{{\left| {\left( {\begin{array}{*{20}{c}}
1&0\\
0&1
\end{array}} \right)} \right| \cdot \left| {\left( {\begin{array}{*{20}{c}}
{1 + {x_4}}&{1 + {x_4} + {z_4}}\\
1&1
\end{array}} \right)} \right|}}{{\left| {\left( {\begin{array}{*{20}{c}}
1&{1 + {x_4} + {z_4}}\\
0&1
\end{array}} \right)} \right| \cdot \left| {\left( {\begin{array}{*{20}{c}}
{1 + {x_4}}&0\\
1&1
\end{array}} \right)} \right|}} = \frac{{ - {z_4}}}{{1 + {x_4}}}
\end{array}\]

\[\begin{array}{l}
{y_1} = {\left[ {C,E;D,B} \right]_A} = {\left[\footnotesize {\left( \begin{array}{c}
1\\
0\\
0
\end{array} \right),\left( \begin{array}{c}
0\\
0\\
1
\end{array} \right);\left( \begin{array}{c}
0\\
1\\
0
\end{array} \right),\left( \begin{array}{c}
1 + {x_4} + {z_4}\\
1\\
{y_4} + {z_4}{t_4}
\end{array} \right)} \right]_{\left( {\scriptstyle1 + {x_4}\atop
{\scriptstyle1\atop
\scriptstyle{y_4}}} \right)}} = \\
 = \frac{{\left| {\left( {\begin{array}{*{20}{c}}
1&0&{1 + {x_4}}\\
0&1&1\\
0&0&{{y_4}}
\end{array}} \right)} \right| \cdot \left| {\left( {\begin{array}{*{20}{c}}
0&{1 + {x_4} + {z_4}}&{1 + {x_4}}\\
0&1&1\\
1&{{y_4} + {z_4}{t_4}}&{{y_4}}
\end{array}} \right)} \right|}}{{\left| {\left( {\begin{array}{*{20}{c}}
1&{1 + {x_4} + {z_4}}&{1 + {x_4}}\\
0&1&1\\
0&{{y_4} + {z_4}{t_4}}&{{y_4}}
\end{array}} \right)} \right| \cdot \left| {\left( {\begin{array}{*{20}{c}}
0&0&{1 + {x_4}}\\
0&1&1\\
1&0&{{y_4}}
\end{array}} \right)} \right|}} = 
\end{array}\]
\[ = \frac{{{y_4} \cdot \left| {\left( {\begin{array}{*{20}{c}}
{1 + {x_4} + {z_4}}&{1 + {x_4}}\\
1&1
\end{array}} \right)} \right|}}{{ - \left( {1 + {x_4}} \right) \cdot \left| {\left( {\begin{array}{*{20}{c}}
1&1\\
{{y_4} + {z_4}{t_4}}&{{y_4}}
\end{array}} \right)} \right|}} = \frac{{{y_4}{z_4}}}{{\left( {1 + {x_4}} \right){z_4}{t_4}}} = \frac{{{y_4}}}{{\left( {1 + {x_4}} \right){t_4}}}\]

\[\begin{array}{l}
{z_1} = {\left[ {E,A;F,B} \right]_C} = {\left[ {\footnotesize \left( \begin{array}{c}
0\\
0\\
1
\end{array} \right),\left( \begin{array}{c}
1 + {x_4}\\
1\\
{y_4}
\end{array} \right);\left( \begin{array}{l}
1\\
1\\
1
\end{array} \right),\left( \begin{array}{c}
1 + {x_4} + {z_4}\\
1\\
{y_4} + {z_4}{t_4}
\end{array} \right)} \right]_{\left( {\scriptstyle1\atop
{\scriptstyle0\atop
\scriptstyle0}} \right)}}=\\
 = \frac{{\left| {\left( {\begin{array}{*{20}{c}}
0&1\\
1&1
\end{array}} \right)} \right| \cdot \left| {\left( {\begin{array}{*{20}{c}}
1&1\\
{{y_4}}&{{y_4} + {z_4}{t_4}}
\end{array}} \right)} \right|}}{{\left| {\left( {\begin{array}{*{20}{c}}
0&1\\
1&{{y_4} + {z_4}{t_4}}
\end{array}} \right)} \right| \cdot \left| {\left( {\begin{array}{*{20}{c}}
1&1\\
{{y_4}}&1
\end{array}} \right)} \right|}} = \frac{{{z_4}{t_4}}}{{1 - {y_4}}}
\end{array}\]

\[\begin{array}{l}
{t_1} = {\left[ {A,C;F,B} \right]_E} = {\footnotesize \left[ {\left( \begin{array}{c}
1 + {x_4}\\
1\\
{y_4}
\end{array} \right),\left( \begin{array}{c}
1\\
0\\
0
\end{array} \right);\left( \begin{array}{l}
1\\
1\\
1
\end{array} \right),\left( \begin{array}{c}
1 + {x_4} + {z_4}\\
1\\
{y_4} + {z_4}{t_4}
\end{array} \right)} \right]_{\left( {\scriptstyle0\atop
{\scriptstyle0\atop
\scriptstyle1}} \right)}} = \\
 = \frac{{\left| {\left( {\begin{array}{*{20}{c}}
{1 + {x_4}}&1\\
1&1
\end{array}} \right)} \right| \cdot \left| {\left( {\begin{array}{*{20}{c}}
1&{1 + {x_4} + {z_4}}\\
0&1
\end{array}} \right)} \right|}}{{\left| {\left( {\begin{array}{*{20}{c}}
{1 + {x_4}}&{1 + {x_4} + {z_4}}\\
1&1
\end{array}} \right)} \right| \cdot \left| {\left( {\begin{array}{*{20}{c}}
1&1\\
0&1
\end{array}} \right)} \right|}} =  - \frac{{{x_4}}}{{{z_4}}}
\end{array}\]

\[\begin{array}{l}
{x_1} = {\left[ {C,A;D,B} \right]_E} = {\left[ {\left( \begin{array}{c}
0\\
1\\
0
\end{array} \right),\left( \begin{array}{c}
1\\
1 + {x_5}\\
{y_5}
\end{array} \right);\left( \begin{array}{c}
{z_5}\\
1\\
{z_5}{t_5}
\end{array} \right),\left( \begin{array}{c}
1\\
0\\
0
\end{array} \right)} \right]_{\left( {\scriptstyle0\atop
{\scriptstyle0\atop
\scriptstyle1}} \right)}} = \\
 = \frac{{\left| {\left( {\begin{array}{*{20}{c}}
0&{{z_5}}\\
1&1
\end{array}} \right)} \right| \cdot \left| {\left( {\begin{array}{*{20}{c}}
1&1\\
{1 + {x_5}}&0
\end{array}} \right)} \right|}}{{\left| {\left( {\begin{array}{*{20}{c}}
0&1\\
1&0
\end{array}} \right)} \right| \cdot \left| {\left( {\begin{array}{*{20}{c}}
1&{{z_5}}\\
{1 + {x_5}}&1
\end{array}} \right)} \right|}} =  - \frac{{{z_5}\left( {1 + {x_5}} \right)}}{{1 - {z_5} - {x_5}{z_5}}}
\end{array}\]

\[\begin{array}{l}
{y_1} = {\left[ {C,E;D,B} \right]_A} = {\left[ {\left( \begin{array}{c}
0\\
1\\
0
\end{array} \right),\left( \begin{array}{c}
0\\
0\\
1
\end{array} \right);\left( \begin{array}{c}
{z_5}\\
1\\
{z_5}{t_5}
\end{array} \right),\left( \begin{array}{c}
1\\
0\\
0
\end{array} \right)} \right]_{\left( {\scriptstyle1\atop
{\scriptstyle1 + {x_5}\atop
\scriptstyle{y_5}}} \right)}} = \\
 = \frac{{\left| {\left( {\begin{array}{*{20}{c}}
0&{{z_5}}&1\\
1&1&{1 + {x_5}}\\
0&{{z_5}{t_5}}&{{y_5}}
\end{array}} \right)} \right| \cdot \left| {\left( {\begin{array}{*{20}{c}}
0&1&1\\
0&0&{1 + {x_5}}\\
1&0&{{y_5}}
\end{array}} \right)} \right|}}{{\left| {\left( {\begin{array}{*{20}{c}}
0&1&1\\
1&0&{1 + {x_5}}\\
0&0&{{y_5}}
\end{array}} \right)} \right| \cdot \left| {\left( {\begin{array}{*{20}{c}}
0&{{z_5}}&1\\
0&1&{1 + {x_5}}\\
1&{{z_5}{t_5}}&{{y_5}}
\end{array}} \right)} \right|}} = \frac{{{z_5}\left( {{t_5} - {y_5}} \right)\left( {1 + {x_5}} \right)}}{{{y_5}\left( {1 - {z_5} - {z_5}{x_5}} \right)}}
\end{array}\]


\[\begin{array}{l}
{t_1} = {\left[ {A,C;F,B} \right]_E} = {\left[ {\left( \begin{array}{c}
1\\
1 + {x_5}\\
{y_5}
\end{array} \right),\left( \begin{array}{c}
0\\
1\\
0
\end{array} \right);\left( \begin{array}{c}
1\\
1\\
1
\end{array} \right),\left( \begin{array}{c}
1\\
0\\
0
\end{array} \right)} \right]_{\left( {\scriptstyle0\atop
{\scriptstyle0\atop
\scriptstyle1}} \right)}} = \\
 = \frac{{\left| {\left( {\begin{array}{*{20}{c}}
1&1\\
{1 + {x_5}}&1
\end{array}} \right)} \right| \cdot \left| {\left( {\begin{array}{*{20}{c}}
0&1\\
1&0
\end{array}} \right)} \right|}}{{\left| {\left( {\begin{array}{*{20}{c}}
1&1\\
{1 + {x_5}}&0
\end{array}} \right)} \right| \cdot \left| {\left( {\begin{array}{*{20}{c}}
0&1\\
1&1
\end{array}} \right)} \right|}} = \frac{{{x_5}}}{{1 + {x_5}}}
\end{array}\]

\[\begin{array}{l}
{t_5} = {\left[ {B,E;D,F} \right]_C} = {\left[ {\left( \begin{array}{c}
1 + {x_4} + {z_4}\\
1\\
{y_4} + {z_4}{t_4}
\end{array} \right),\left( \begin{array}{c}
0\\
0\\
1
\end{array} \right);\left( \begin{array}{c}
0\\
1\\
0
\end{array} \right),\left( \begin{array}{c}
1\\
1\\
1
\end{array} \right)} \right]_{\left( {\scriptstyle1\atop
{\scriptstyle0\atop
\scriptstyle0}} \right)}} = \\
 = \frac{{\left| {\left( {\begin{array}{*{20}{c}}
1&1\\
{{y_4} + {z_4}{t_4}}&0
\end{array}} \right)} \right| \cdot \left| {\left( {\begin{array}{*{20}{c}}
0&1\\
1&1
\end{array}} \right)} \right|}}{{\left| {\left( {\begin{array}{*{20}{c}}
1&1\\
{{y_4} + {z_4}{t_4}}&1
\end{array}} \right)} \right| \cdot \left| {\left( {\begin{array}{*{20}{c}}
0&1\\
1&0
\end{array}} \right)} \right|}} =  - \frac{{{y_4} + {z_4}{t_4}}}{{1 - {y_4} - {z_4}{t_4}}}
\end{array}\]

After multiplying by the denominators and omitting the negligible terms we get:
\[{x_1} \approx  - {z_4}\]
\[{y_1}{t_4} \approx {y_4}\]
\[{z_1} \approx 0\]
\[0 \approx  - {x_4}\]
\[{x_1} \approx  - {z_5}\]
\[{y_1}{y_5} \approx 0\]
\[{t_1} \approx {x_5}\]
\[{t_5} \approx  - {y_4}\]
We may conclude 
\[d{x_1} =  - d{z_4}\]
\[{y_1} \cdot d{t_4} = d{y_4}\]
\[d{z_1} = 0\]
\[d{x_4} = 0\]
\[d{x_1} =  - d{z_5}\]
\[{y_1} \cdot d{y_5} = 0\]
\[d{t_1} = d{x_5}\]
\[dt_5=-dy_4\]

The last four equations allow us to eliminate $dx_5, dy_5, dz_5, dt_5$. 
The former four equation allow to eliminate 
$dz_1, dx_4,dy_4,dz_4$. SO the cotangent space is spanned by $dx_1,dy_1,dt_1,dt_4$.

\newpage
\subsection{$\mathcal{C}.2$} 
\begin{center}
	\includegraphics{Case_C2.png}
\end{center}

We shall not discuss the last two charts, as it is obvious all that the first 
three charts contain all the cross ratios and define the last charts completely.
On the first  three charts, we choose the following coordinates:
\[\begin{array}{*{20}{c}}
\begin{array}{l}
{A_1} = \left( {1:0:0} \right)\\
{B_1} = \left( {1:1:1} \right)\\
{C_1} = \left( {0:1:0} \right)\\
{D_1} = \left( {{x_1}:1:{y_1}} \right)\\
{E_1} = \left( {0:0:1} \right)\\
{F_1} = \left( {{z_1}{t_1}:{z_1}:1} \right)\\
{x_1},{z_1},{t_1} \sim 0 \ne {y_1}
\end{array}&\vline& \begin{array}{l}
{A_2} = \left( {1:0:0} \right)\\
{B_2} = \left( {1:{x_2}:{x_2}{y_2}} \right)\\
{C_2} = \left( {0:1:0} \right)\\
{D_2} = \left( {0:0:1} \right)\\
{E_2} = \left( {1:1:1} \right)\\
{F_2} = \left( {{z_2}:1 + {t_2}:1} \right)\\
{x_2},{z_2},{t_2} \sim 0 \ne {y_2}
\end{array}&\vline& \begin{array}{l}
{A_3} = \left( {1:0:0} \right)\\
{B_3} = \left( {0:1:0} \right)\\
{C_3} = \left( {1:1 + {x_3}:{y_3}} \right)\\
{D_3} = \left( 1:1 + {x_3} + {z_3}{t_3}\right. \\
\hfill \left.:{y_3} + {z_3} \right)\\
{E_3} = \left( {0:0:1} \right)\\
{F_3} = \left( {1:1:1} \right)\\
{x_3},{y_3},{z_3},{t_3} \sim 0
\end{array}
\end{array}\]

We write the formulas relating the coordinates of different charts.

\[\begin{array}{l}
{y_1} = {\left[ {C,E;D,B} \right]_A} = {\left[ {\left( \begin{array}{c}
0\\
1\\
0
\end{array} \right),\left( \begin{array}{l}
1\\
1\\
1
\end{array} \right);\left( \begin{array}{c}
0\\
0\\
1
\end{array} \right),\left( \begin{array}{c}
1\\
{x_2}\\
{x_2}{y_2}
\end{array} \right)} \right]_{\left( {\scriptstyle1\atop
{\scriptstyle0\atop
\scriptstyle0}} \right)}} = \end{array}\]
\[ = \frac{{\left| {\left( {\begin{array}{*{20}{c}}
1&0\\
0&1
\end{array}} \right)} \right| \cdot \left| {\left( {\begin{array}{*{20}{c}}
1&{{x_2}}\\
1&{{x_2}{y_2}}
\end{array}} \right)} \right|}}{{\left| {\left( {\begin{array}{*{20}{c}}
1&{{x_2}}\\
0&{{x_2}{y_2}}
\end{array}} \right)} \right| \cdot \left| {\left( {\begin{array}{*{20}{c}}
1&0\\
1&1
\end{array}} \right)} \right|}} = 
\frac{{{x_2}\left( {{y_2} - 1} \right)}}{{{x_2}{y_2}}} = \frac{{{y_2} - 1}}{{{y_2}}}\]

\[\begin{array}{l}
{z_1} = {\left[ {E,C;F,B} \right]_A} = {\left[ {\left( \begin{array}{l}
1\\
1\\
1
\end{array} \right),\left( \begin{array}{c}
0\\
1\\
0
\end{array} \right);\left( \begin{array}{c}
{z_2}\\
1 + {t_2}\\
1
\end{array} \right),\left( \begin{array}{c}
1\\
{x_2}\\
{x_2}{y_2}
\end{array} \right)} \right]_{\left( {\scriptstyle1\atop
{\scriptstyle0\atop
\scriptstyle0}} \right)}} = \end{array}\]
\[ = \frac{{\left| {\left( {\begin{array}{*{20}{c}}
1&{1 + {t_2}}\\
1&1
\end{array}} \right)} \right| \cdot \left| {\left( {\begin{array}{*{20}{c}}
1&{{x_2}}\\
0&{{x_2}{y_2}}
\end{array}} \right)} \right|}}{{\left| {\left( {\begin{array}{*{20}{c}}
1&{{x_2}}\\
1&{{x_2}{y_2}}
\end{array}} \right)} \right| \cdot \left| {\left( {\begin{array}{*{20}{c}}
1&{1 + {t_2}}\\
0&1
\end{array}} \right)} \right|}} = \frac{{{t_2} \cdot {x_2}{y_2}}}{{{x_2}\left( {1 - {y_2}} \right)}} = \frac{{{y_2}{t_2}}}{{1 - {y_2}}}\]

\[\begin{array}{l}
{t_1} = {\left[ {C,A;F,B} \right]_E} = {\left[ {\left( \begin{array}{c}
0\\
1\\
0
\end{array} \right),\left( \begin{array}{c}
1\\
0\\
0
\end{array} \right);\left( \begin{array}{c}
{z_2}\\
1 + {t_2}\\
1
\end{array} \right),\left( \begin{array}{c}
1\\
{x_2}\\
{x_2}{y_2}
\end{array} \right)} \right]_{\left( {\scriptstyle1\atop
{\scriptstyle1\atop
\scriptstyle1}} \right)}} = \end{array}\]
\[ = \frac{{\left| {\left( {\begin{array}{*{20}{c}}
0&{{z_2}}&1\\
1&{1 + {t_2}}&1\\
0&1&1
\end{array}} \right)} \right| \cdot \left| {\left( {\begin{array}{*{20}{c}}
1&1&1\\
0&{{x_2}}&1\\
0&{{x_2}{y_2}}&1
\end{array}} \right)} \right|}}{{\left| {\left( {\begin{array}{*{20}{c}}
0&1&1\\
1&{{x_2}}&1\\
0&{{x_2}{y_2}}&1
\end{array}} \right)} \right| \cdot \left| {\left( {\begin{array}{*{20}{c}}
1&{{z_2}}&1\\
0&{1 + {t_2}}&1\\
0&1&1
\end{array}} \right)} \right|}} = 
\frac{{\left( {1 - {z_2}} \right) \cdot {x_2}\left( {1 - {y_2}} \right)}}{{\left( {{x_2}{y_2} - 1} \right){t_2}}}\]

\[ \begin{array}{l}
{x_1} = {\left[ {C,A;D,B} \right]_E} = {\left[\footnotesize {\left( \begin{array}{c}
1\\
1 + {x_3}\\
{y_3}
\end{array} \right),\left( \begin{array}{c}
1\\
0\\
0
\end{array} \right);\left( \begin{array}{c}
1\\
1 + {x_3} + {z_3}{t_3}\\
{y_3} + {z_3}
\end{array} \right),\left( \begin{array}{c}
0\\
1\\
0
\end{array} \right)} \right]_{\left( {\scriptstyle0\atop
{\scriptstyle0\atop
\scriptstyle1}} \right)}} = \end{array}\]
\[ = \frac{{\left| {\left( {\begin{array}{*{20}{c}}
1&1\\
{1 + {x_3}}&{1 + {x_3} + {z_3}{t_3}}
\end{array}} \right)} \right| \cdot \left| {\left( {\begin{array}{*{20}{c}}
1&0\\
0&1
\end{array}} \right)} \right|}}{{\left| {\left( {\begin{array}{*{20}{c}}
1&0\\
{1 + {x_3}}&1
\end{array}} \right)} \right| \cdot \left| {\left( {\begin{array}{*{20}{c}}
1&1\\
0&{1 + {x_3} + {z_3}{t_3}}
\end{array}} \right)} \right|}} = \frac{{{z_3}{t_3}}}{{1 + {x_3} + {z_3}{t_3}}}\]

\[\begin{array}{l}
{y_1} = {\left[ {C,E;D,B} \right]_A} = {\left[\footnotesize {\left( \begin{array}{c}
1\\
1 + {x_3}\\
{y_3}
\end{array} \right),\left( \begin{array}{c}
0\\
0\\
1
\end{array} \right);\left( \begin{array}{c}
1\\
1 + {x_3} + {z_3}{t_3}\\
{y_3} + {z_3}
\end{array} \right),\left( \begin{array}{c}
0\\
1\\
0
\end{array} \right)} \right]_{\left( {\scriptstyle1\atop
{\scriptstyle0\atop
\scriptstyle0}} \right)}} = \\
 = \frac{{\left| {\left( {\begin{array}{*{20}{c}}
{1 + {x_3}}&{1 + {x_3} + {z_3}{t_3}}\\
{{y_3}}&{{y_3} + {z_3}}
\end{array}} \right)} \right| \cdot \left| {\left( {\begin{array}{*{20}{c}}
0&1\\
1&0
\end{array}} \right)} \right|}}{{\left| {\left( {\begin{array}{*{20}{c}}
{1 + {x_3}}&1\\
{{y_3}}&0
\end{array}} \right)} \right| \cdot \left| {\left( {\begin{array}{*{20}{c}}
0&{1 + {x_3} + {z_3}{t_3}}\\
1&{{y_3} + {z_3}}
\end{array}} \right)} \right|}} = \frac{{ - \left| {\left( {\begin{array}{*{20}{c}}
{1 + {x_3}}&{{z_3}{t_3}}\\
{{y_3}}&{{z_3}}
\end{array}} \right)} \right|}}{{{y_3}\left( {1 + {x_3} + {z_3}{t_3}} \right)}} = 
\end{array}\]
\[
=\frac{{ - {z_3}\left( {1 + {x_3} - {y_3}{t_3}} \right)}}{{{y_3}\left( {1 + {x_3} + {z_3}{t_3}} \right)}}
\]

\[\begin{array}{l}
{z_1} = {\left[ {E,C;F,B} \right]_A} = {\left[ {\left( \begin{array}{c}
0\\
0\\
1
\end{array} \right),\left( \begin{array}{c}
1\\
1 + {x_3}\\
{y_3}
\end{array} \right);\left( \begin{array}{l}
1\\
1\\
1
\end{array} \right),\left( \begin{array}{c}
0\\
1\\
0
\end{array} \right)} \right]_{\left( {\scriptstyle1\atop
{\scriptstyle0\atop
\scriptstyle0}} \right)}} = \\
 = \frac{{\left| {\left( {\begin{array}{*{20}{c}}
0&1\\
1&1
\end{array}} \right)} \right| \cdot \left| {\left( {\begin{array}{*{20}{c}}
{1 + {x_3}}&1\\
{{y_3}}&0
\end{array}} \right)} \right|}}{{\left| {\left( {\begin{array}{*{20}{c}}
0&1\\
1&0
\end{array}} \right)} \right| \cdot \left| {\left( {\begin{array}{*{20}{c}}
{1 + {x_3}}&1\\
{{y_3}}&1
\end{array}} \right)} \right|}} = -\frac{{{y_3}}}{{1 + {x_3} - {y_3}}}
\end{array}\]

\[\begin{array}{l}
{t_1} = {\left[ {C,A;F,B} \right]_E} = {\left[\small {\left( \begin{array}{c}
1\\
1 + {x_3}\\
{y_3}
\end{array} \right),\left( \begin{array}{c}
1\\
0\\
0
\end{array} \right);\left( \begin{array}{l}
1\\
1\\
1
\end{array} \right),\left( \begin{array}{c}
0\\
1\\
0
\end{array} \right)} \right]_{\left( {\scriptstyle0\atop
{\scriptstyle0\atop
\scriptstyle1}} \right)}} = \\
 = \frac{{\left| {\left( {\begin{array}{*{20}{c}}
1&1\\
{1 + {x_3}}&1
\end{array}} \right)} \right| \cdot \left| {\left( {\begin{array}{*{20}{c}}
1&0\\
0&1
\end{array}} \right)} \right|}}{{\left| {\left( {\begin{array}{*{20}{c}}
1&0\\
{1 + {x_3}}&1
\end{array}} \right)} \right| \cdot \left| {\left( {\begin{array}{*{20}{c}}
1&1\\
0&1
\end{array}} \right)} \right|}} =  - {x_3}
\end{array}\]

\[{z_2} = {\left[ {D,A;F,E} \right]_C} = {\left[\tiny {\left( \begin{array}{c}
1\\
1 + {x_3} + {z_3}{t_3}\\
{y_3} + {z_3}
\end{array} \right),\left( \begin{array}{c}
1\\
0\\
0
\end{array} \right);\left( \begin{array}{l}
1\\
1\\
1
\end{array} \right),\left( \begin{array}{c}
0\\
0\\
1
\end{array} \right)} \right]_{\left( {\scriptstyle1\atop
{\scriptstyle1 + {x_3}\atop
\scriptstyle{y_3}}} \right)}} = \]
\[ = \frac{{\left| {\left( {\begin{array}{*{20}{c}}
1&1&1\\
{1 + {x_3} + {z_3}{t_3}}&1&{1 + {x_3}}\\
{{y_3} + {z_3}}&1&{{y_3}}
\end{array}} \right)} \right| \cdot \left| {\left( {\begin{array}{*{20}{c}}
1&0&1\\
0&0&{1 + {x_3}}\\
0&1&{{y_3}}
\end{array}} \right)} \right|}}{{\left| {\left( {\begin{array}{*{20}{c}}
1&0&1\\
{1 + {x_3} + {z_3}{t_3}}&0&{1 + {x_3}}\\
{{y_3} + {z_3}}&1&{{y_3}}
\end{array}} \right)} \right| \cdot \left| {\left( {\begin{array}{*{20}{c}}
1&1&1\\
0&1&{1 + {x_3}}\\
0&1&{{y_3}}
\end{array}} \right)} \right|}} = \]
\[ = \frac{{\left| {\left( {\begin{array}{*{20}{c}}
0&1&1\\
{{z_3}{t_3}}&1&{1 + {x_3}}\\
{{z_3}}&1&{{y_3}}
\end{array}} \right)} \right| \cdot \left( {1 + {x_3}} \right)}}{{\left| {\left( {\begin{array}{*{20}{c}}
0&0&1\\
{{z_3}{t_3}}&0&{1 + {x_3}}\\
{{z_3}}&1&{{y_3}}
\end{array}} \right)} \right| \cdot \left( {1 + {x_3} - {y_3}} \right)}} =\]
\[= \frac{{\left| {\left( {\begin{array}{*{20}{c}}
0&1&1\\
{{z_3}{t_3}}&0&{{x_3}}\\
{{z_3}}&0&{{y_3} - 1}
\end{array}} \right)} \right| \cdot \left( {1 + {x_3}} \right)}}{{{z_3}{t_3} \cdot \left( {1 + {x_3} - {y_3}} \right)}}  = \frac{{\left( {{x_3} - {t_3}{y_3} + {t_3}} \right) \cdot \left( {1 + {x_3}} \right)}}{{{t_3} \cdot \left( {1 + {x_3} - {y_3}} \right)}}\]

After multiplying by the denominators and omitting the negligible terms, we get 
\[{y_1}{y_2} = {y_2} - 1\]
\[{z_1}\left( {1 - {y_2}} \right) = {y_2}{t_2}\]
\[0 \approx {x_2}\left( {1 - {y_2}} \right)\]
\[{x_1} \approx 0\]
\[{y_1}{y_3} \approx  - {z_3}\]
\[{z_1} \approx -{y_3}\]
\[{t_1} =  - {x_3}\]
\[0 \approx {x_3} + {t_3}\]

We conclude
\[y_2dy_1=(1-y_1)dy_2\]
\[\left( {1 - {y_2}} \right) dz_1 = {y_2}dt_2\]
\[0 = \left( {1 - {y_2}} \right)dx_2\]
\[dx_1 = 0\]
\[y_3dy_1 =  - dz_3\]
\[dz_1 = -dy_3\]
\[dt_1 =  - dx_3\]
\[0 = d{x_3} + d{t_3}\]

The last equation may be used to eliminate $dt_3$, the 3 equations before it may be used to exclude $dx_3,dy_3,dz_3$.
The first four among the last eight equations allow to exclude $dx_1,dx_2,dy_2,dt_2$. 
So the cotangent bundle is spanned $dy_1, dz_1, dt_1, dz_2$.

\newpage
\subsection{$\mathcal{C}.4$} 
\begin{center}
	\includegraphics{Case_C4.png}
\end{center}
The information in the charts 2 and 4 follows from charts 1, 3, and 5.
\[\begin{array}{*{20}{c}}
\begin{array}{l}
{A_1} = \left( {1:0:0} \right)\\
{B_1} = \left( {1:1:1} \right)\\
{C_1} = \left( {0:1:0} \right)\\
{D_1} = \left( {{x_1}:1:{y_1}} \right)\\
{E_1} = \left( {0:0:1} \right)\\
{F_1} = \left( {{z_1}{t_1}:{z_1}:1} \right)\\
{x_1},{z_1},{t_1} \sim 0 \ne {y_1}
\end{array}&\vline& \begin{array}{l}
{A_3} = \left( {1:0:0} \right)\\
{B_3} = \left( {0:1:0} \right)\\
{C_3} = \left( {1 + {x_3}:1:{y_3}} \right)\\
{D_3} = \left( 1 + {x_3} + {z_3}\right. \\
\hfill \left. :1:{y_3} + {z_3}{t_3} \right)\\
{E_3} = \left( {0:0:1} \right)\\
{F_3} = \left( {1:1:1} \right)\\
{x_3},{y_3},{z_3} \sim 0
\end{array}&\vline& \begin{array}{l}
{A_5} = \left( {1:1 + {x_5}:{y_5}} \right)\\
{B_5} = \left( {1:1 + {x_5} + {z_5}} \right.\\
\left. {\,\,\,\,\,\,\,\,\,\,\,\,\,\,\,\,\,\,:{y_5} + {z_5}{t_5}} \right)\\
{C_5} = \left( {1:0:0} \right)\\
{D_5} = \left( {0:1:0} \right)\\
{E_5} = \left( {0:0:1} \right)\\
{F_5} = \left( {1:1:1} \right)\\
{x_5},{y_5},{z_5},{t_5} \sim 0
\end{array}
\end{array}\]

The formulas connecting different charts:
\[\begin{array}{l}
{x_1} = {\left[ {C,A;D,B} \right]_E} = {\left[\footnotesize {\left( \begin{array}{c}
1 + {x_3}\\
1\\
{y_3}
\end{array} \right),\left( \begin{array}{c}
1\\
0\\
0
\end{array} \right);\left( \begin{array}{c}
1 + {x_3} + {z_3}\\
1\\
{y_3} + {z_3}{t_3}
\end{array} \right),\left( \begin{array}{c}
0\\
1\\
0
\end{array} \right)} \right]_{\left( {\scriptstyle0\atop
{\scriptstyle0\atop
\scriptstyle1}} \right)}} = 
\end{array}\]
\[ = \frac{{\left| {\left( {\begin{array}{*{20}{c}}
{1 + {x_3}}&{1 + {x_3} + {z_3}}\\
1&1
\end{array}} \right)} \right| \cdot \left| {\left( {\begin{array}{*{20}{c}}
1&0\\
0&1
\end{array}} \right)} \right|}}{{\left| {\left( {\begin{array}{*{20}{c}}
{1 + {x_3}}&0\\
1&1
\end{array}} \right)} \right| \cdot \left| {\left( {\begin{array}{*{20}{c}}
1&{1 + {x_3} + {z_3}}\\
0&1
\end{array}} \right)} \right|}} = \frac{{ - {z_3}}}{{1 + {x_3}}}
\]

\[\begin{array}{l}
{y_1} = {\left[ {C,E;D,B} \right]_A} = 
{\left[\footnotesize {\left( \begin{array}{c}
1 + {x_3}\\
1\\
{y_3}
\end{array} \right),\left( \begin{array}{c}
0\\
0\\
1
\end{array} \right);\left( \begin{array}{c}
1 + {x_3} + {z_3}\\
1\\
{y_3} + {z_3}{t_3}
\end{array} \right),\left( \begin{array}{c}
0\\
1\\
0
\end{array} \right)} \right]_{\left( {\scriptstyle1\atop
{\scriptstyle0\atop
\scriptstyle0}} \right)}} = 
\end{array}\]
\[ = \frac{{\left| {\left( {\begin{array}{*{20}{c}}
1&1\\
{{y_3}}&{{y_3} + {z_3}{t_3}}
\end{array}} \right)} \right| \cdot \left| {\left( {\begin{array}{*{20}{c}}
0&1\\
1&0
\end{array}} \right)} \right|}}{{\left| {\left( {\begin{array}{*{20}{c}}
1&1\\
{{y_3}}&0
\end{array}} \right)} \right| \cdot \left| {\left( {\begin{array}{*{20}{c}}
0&1\\
1&{{y_3} + {z_3}{t_3}}
\end{array}} \right)} \right|}} = \frac{{ - {z_3}{t_3}}}{{{y_3}}}
\]

\[\begin{array}{l}
{z_1} = {\left[ {E,C;F,B} \right]_A} = {\left[ {\left( \begin{array}{c}
0\\
0\\
1
\end{array} \right),\left( \begin{array}{c}
1 + {x_3}\\
1\\
{y_3}
\end{array} \right);\left( \begin{array}{l}
1\\
1\\
1
\end{array} \right),\left( \begin{array}{c}
0\\
1\\
0
\end{array} \right)} \right]_{\left( {\scriptstyle1\atop
{\scriptstyle0\atop
\scriptstyle0}} \right)}} = 
\end{array}\]
\[ = \frac{{\left| {\left( {\begin{array}{*{20}{c}}
0&1\\
1&1
\end{array}} \right)} \right| \cdot \left| {\left( {\begin{array}{*{20}{c}}
1&1\\
{{y_3}}&0
\end{array}} \right)} \right|}}{{\left| {\left( {\begin{array}{*{20}{c}}
0&1\\
1&0
\end{array}} \right)} \right| \cdot \left| {\left( {\begin{array}{*{20}{c}}
1&1\\
{{y_3}}&1
\end{array}} \right)} \right|}} = \frac{{{y_3}}}{{{y_3} - 1}}
\]

\[\begin{array}{l}
{t_1} = {\left[ {C,A;F,B} \right]_E} = {\left[ {\left( \begin{array}{c}
1 + {x_3}\\
1\\
{y_3}
\end{array} \right),\left( \begin{array}{c}
1\\
0\\
0
\end{array} \right);\left( \begin{array}{l}
1\\
1\\
1
\end{array} \right),\left( \begin{array}{c}
0\\
1\\
0
\end{array} \right)} \right]_{\left( {\scriptstyle0\atop
{\scriptstyle0\atop
\scriptstyle1}} \right)}} = \end{array}\]
\[ = \frac{{\left| {\left( {\begin{array}{*{20}{c}}
{1 + {x_3}}&1\\
1&1
\end{array}} \right)} \right| \cdot \left| {\left( {\begin{array}{*{20}{c}}
1&0\\
0&1
\end{array}} \right)} \right|}}{{\left| {\left( {\begin{array}{*{20}{c}}
{1 + {x_3}}&0\\
1&1
\end{array}} \right)} \right| \cdot \left| {\left( {\begin{array}{*{20}{c}}
1&1\\
0&1
\end{array}} \right)} \right|}} = \frac{{{x_3}}}{{1 + {x_3}}}
\]

\[\begin{array}{l}
{x_1} = {\left[ {C,A;D,B} \right]_E} = {\left[ {\left( \begin{array}{c}
1\\
0\\
0
\end{array} \right),\left( \begin{array}{c}
1\\
1 + {x_5}\\
{y_5}
\end{array} \right);\left( \begin{array}{c}
0\\
1\\
0
\end{array} \right),\left( \begin{array}{c}
1\\
1 + {x_5} + {z_5}\\
{y_5} + {z_5}{t_5}
\end{array} \right)} \right]_{\left( {\scriptstyle0\atop
{\scriptstyle0\atop
\scriptstyle1}} \right)}} = \end{array}\]
\[ = \frac{{\left| {\left( {\begin{array}{*{20}{c}}
1&0\\
0&1
\end{array}} \right)} \right| \cdot \left| {\left( {\begin{array}{*{20}{c}}
1&1\\
{1 + {x_5}}&{1 + {x_5} + {z_5}}
\end{array}} \right)} \right|}}{{\left| {\left( {\begin{array}{*{20}{c}}
1&1\\
0&{1 + {x_5} + {z_5}}
\end{array}} \right)} \right| \cdot \left| {\left( {\begin{array}{*{20}{c}}
1&0\\
{1 + {x_5}}&1
\end{array}} \right)} \right|}} = \frac{{{z_5}}}{{1 + {x_5} + {z_5}}}\]

\[{y_1} = {\left[ {C,E;D,B} \right]_A} = {\left[ {\left( \begin{array}{c}
1\\
0\\
0
\end{array} \right),\left( \begin{array}{c}
0\\
0\\
1
\end{array} \right);\left( \begin{array}{c}
0\\
1\\
0
\end{array} \right),\left( \begin{array}{c}
1\\
1 + {x_5} + {z_5}\\
{y_5} + {z_5}{t_5}
\end{array} \right)} \right]_{\left( {\scriptstyle1\atop
{\scriptstyle1 + {x_5}\atop
\scriptstyle{y_5}}} \right)}} = \]
\[ = \frac{{\left| {\left( {\begin{array}{*{20}{c}}
1&0&1\\
0&1&{1 + {x_5}}\\
0&0&{{y_5}}
\end{array}} \right)} \right| \cdot \left| {\left( {\begin{array}{*{20}{c}}
0&1&1\\
0&{1 + {x_5} + {z_5}}&{1 + {x_5}}\\
1&{{y_5} + {z_5}{t_5}}&{{y_5}}
\end{array}} \right)} \right|}}{{\left| {\left( {\begin{array}{*{20}{c}}
1&1&1\\
0&{1 + {x_5} + {z_5}}&{1 + {x_5}}\\
0&{{y_5} + {z_5}{t_5}}&{{y_5}}
\end{array}} \right)} \right| \cdot \left| {\left( {\begin{array}{*{20}{c}}
0&0&1\\
0&1&{1 + {x_5}}\\
1&0&{{y_5}}
\end{array}} \right)} \right|}} = \]
\[ = \frac{{{y_5} \cdot \left| {\left( {\begin{array}{*{20}{c}}
0&0&1\\
0&{{z_5}}&{1 + {x_5}}\\
1&{{z_5}{t_5}}&{{y_5}}
\end{array}} \right)} \right|}}{{ - \left| {\left( {\begin{array}{*{20}{c}}
1&0&1\\
0&{{z_5}}&{1 + {x_5}}\\
0&{{z_5}{t_5}}&{{y_5}}
\end{array}} \right)} \right|}} = \frac{{{y_5}{z_5}}}{{{z_5}\left( {{y_5} - {t_5} - {x_5}{t_5}} \right)}} = \frac{{{y_5}}}{{{y_5} - {t_5} - {x_5}{t_5}}}\]

\[{z_1} = {\left[ {E,C;F,B} \right]_A} = {\left[ {\left( \begin{array}{c}
0\\
0\\
1
\end{array} \right),\left( \begin{array}{c}
1\\
0\\
0
\end{array} \right);\left( \begin{array}{l}
1\\
1\\
1
\end{array} \right),\left( \begin{array}{c}
1\\
1 + {x_5} + {z_5}\\
{y_5} + {z_5}{t_5}
\end{array} \right)} \right]_{\left( {\scriptstyle1\atop
{\scriptstyle1 + {x_5}\atop
\scriptstyle{y_5}}} \right)}} = \]
\[ = \frac{{\left| {\left( {\begin{array}{*{20}{c}}
0&1&1\\
0&1&{1 + {x_5}}\\
1&1&{{y_5}}
\end{array}} \right)} \right| \cdot \left| {\left( {\begin{array}{*{20}{c}}
1&1&1\\
0&{1 + {x_5} + {z_5}}&{1 + {x_5}}\\
0&{{y_5} + {z_5}{t_5}}&{{y_5}}
\end{array}} \right)} \right|}}{{\left| {\left( {\begin{array}{*{20}{c}}
0&1&1\\
0&{1 + {x_5} + {z_5}}&{1 + {x_5}}\\
1&{{y_5} + {z_5}{t_5}}&{{y_5}}
\end{array}} \right)} \right| \cdot \left| {\left( {\begin{array}{*{20}{c}}
1&1&1\\
0&1&{1 + {x_5}}\\
0&1&{{y_5}}
\end{array}} \right)} \right|}} = \]
\[ = \frac{{{x_5} \cdot {z_5}\left( {{y_5} - {t_5} - {x_5}{t_5}} \right)}}{{{z_5}\left( {1 + {x_5} - {y_5}} \right)}} 
= \frac{{{x_5}\left( {{y_5} - {t_5} - {x_5}{t_5}} \right)}}{{1 + {x_5} - {y_5}}}\]

\[{t_1} = {\left[ {C,A;F,B} \right]_E} = {\left[ {\left( \begin{array}{c}
1\\
0\\
0
\end{array} \right),\left( \begin{array}{c}
1\\
1 + {x_5}\\
{y_5}
\end{array} \right);\left( \begin{array}{l}
1\\
1\\
1
\end{array} \right),\left( \begin{array}{c}
1\\
1 + {x_5} + {z_5}\\
{y_5} + {z_5}{t_5}
\end{array} \right)} \right]_{\left( {\scriptstyle0\atop
{\scriptstyle0\atop
\scriptstyle1}} \right)}} = \]

\[ = \frac{{\left| {\left( {\begin{array}{*{20}{c}}
1&1\\
0&1
\end{array}} \right)} \right| \cdot \left| {\left( {\begin{array}{*{20}{c}}
1&1\\
{1 + {x_5}}&{1 + {x_5} + {z_5}}
\end{array}} \right)} \right|}}{{\left| {\left( {\begin{array}{*{20}{c}}
1&1\\
0&{1 + {x_5} + {z_5}}
\end{array}} \right)} \right| \cdot \left| {\left( {\begin{array}{*{20}{c}}
1&1\\
{1 + {x_5}}&1
\end{array}} \right)} \right|}} = \frac{{{z_5}}}{{ - {x_5}\left( {1 + {x_5} + {z_5}} \right)}}\]

As these equations produce linearly dependent conditions, we add another equation:

\[\begin{array}{l}
{y_5} = {\left[ {C,E;A,F} \right]_D} = {\left[ {\left( \begin{array}{c}
1 + {x_3}\\
1\\
{y_3}
\end{array} \right),\left( \begin{array}{c}
0\\
0\\
1
\end{array} \right);\left( \begin{array}{c}
1\\
0\\
0
\end{array} \right),\left( \begin{array}{l}
1\\
1\\
1
\end{array} \right)} \right]_{\left( {\scriptstyle1 + {x_3} + {z_3}\atop
{\scriptstyle1\atop
\scriptstyle{y_3} + {z_3}{t_3}}} \right)}} = \\
 = \frac{{\left| {\left( {\begin{array}{*{20}{c}}
{1 + {x_3}}&1&{1 + {x_3} + {z_3}}\\
1&0&1\\
{{y_3}}&0&{{y_3} + {z_3}{t_3}}
\end{array}} \right)} \right| \cdot \left| {\left( {\begin{array}{*{20}{c}}
0&1&{1 + {x_3} + {z_3}}\\
0&1&1\\
1&1&{{y_3} + {z_3}{t_3}}
\end{array}} \right)} \right|}}{{\left| {\left( {\begin{array}{*{20}{c}}
{1 + {x_3}}&1&{1 + {x_3} + {z_3}}\\
1&1&1\\
{{y_3}}&1&{{y_3} + {z_3}{t_3}}
\end{array}} \right)} \right| \cdot \left| {\left( {\begin{array}{*{20}{c}}
0&1&{1 + {x_3} + {z_3}}\\
0&0&1\\
1&0&{{y_3} + {z_3}{t_3}}
\end{array}} \right)} \right|}} = \end{array}\]
\[ = \frac{{{z_3}{t_3}\left( {{x_3} + {z_3}} \right)}}{{\left| {\left( {\begin{array}{*{20}{c}}
{1 + {x_3}}&1&{{z_3}}\\
1&1&0\\
{{y_3}}&1&{{z_3}{t_3}}
\end{array}} \right)} \right|}} = \frac{{{t_3}\left( {{x_3} + {z_3}} \right)}}
{{\left| {\left( {\begin{array}{*{20}{c}}
{{x_3}}&0&1\\
1&1&0\\
{{y_3}}&1&{{t_3}}
\end{array}} \right)} \right|}} = \frac{{{t_3}\left( {{x_3} + {z_3}} \right)}}{{1 - {y_3} + {x_3}{t_3}}}\]

After multiplying by the denominators, and omitting the negligible terms we get 

\[{x_1} \approx  - {z_3}\]
\[{y_1}{y_3} =  - {z_3}{t_3}\]
\[{z_1} \approx  - {y_3}\]
\[{t_1} \approx {x_3}\]
\[{x_1} \approx {z_5}\]
\[{y_1}\left( {{y_5} - {t_5}} \right) \approx {y_5}\]
\[{z_1} \approx 0\]
\[0 \approx {z_5}\]
\[{y_5} \approx {t_3}\left( {{x_3} + {z_3}} \right)\]

It follows that
\[dx_1    =  - dz_3\]
\[y_1dy_3 =  - t_3dz_3\]
\[dz_1    =  - dy_3\]
\[dt_1    = dx_3\]
\[dx_1 = dz_5\]
\[0 = (1-y_1) dy_5+y_1dt_5\]
\[dz_1 = 0\]
\[0 = dz_5\]
\[y_5 \approx {t_3}\left( {dx_3 + dz_3} \right)\]

From the sixth equation we may express $dt_5$ by $dy_5$ and from the ninth equation we may express $dy_5$ from the coordinates 
related to the third chart, so we may omit the both. From the eighth equation we get that $dz_5$ is zero,
hence by the fifth equation $dx_1$ is zero and hence by the first equation $dz_3$ is zero and by the second equation $dy_3$ is zero, 
and by the third equation (or directly by the seventh equation) also $dz_1$ is zero.
So the cotangent space is spanned by $dy_1,dt_1,dx_3, dt_3,dx_5$. But the fourth equation relates $dt_1$ to $dx_3$, 
so the cotangent space is at most four-dimensional.

\subsection{$\mathcal{D}$}
\begin{center}
	\includegraphics{Case_D.png}
\end{center}

We may choose local coordinates related to different charts:
\[
\]

We shall start by explaining that the coordinates of the first chart are related to the second chart,
then that the coordinates of the fourth charts are related to the other charts, in the end we shall relate 
between the coordinates of the second and the third charts.

\[%
} \right)} \right|}} = \frac{{{z_2} - {x_2}}}{{1 - {x_2}}}\]

After multiplying by the denominators and omitting negligible terms, we get the equations
\[\begin{array}{l}
{x_1} \approx {y_2}\\
{y_1} \approx  - {x_2}\\
\left( {1 - {t_2}} \right){z_1} \approx {y_2}\\
{t_1} \approx {z_2} - {x_2}
\end{array}\]

Since $1\neq t_2$, the coordinates of the first chart are locally expressed in terms of the second chart.

\[%
} \right)} \right|}} = \frac{{{z_2} - {x_2}}}{{1 - {x_2}}}\]

After multiplying by the denominators and omitting negligible terms, we get 
\[\begin{array}{l}
{y_3} \approx {y_4}{z_4}\\
{z_3}{z_4} = {t_4}\\
{t_3} = {z_4}\\
{x_4} \approx {z_2} - {x_2}
\end{array}\]

Passing to differentials we get:
\[\begin{array}{l}
d{y_3} = {z_4} \cdot d{y_4}\\
{z_4} \cdot d{z_3} = d{t_4}\\
d{t_3} = d{z_4}\\
d{x_4} = d{z_2} - d{x_2}
\end{array}\]

We see that $dx_4,dy_4, dz_4$ and $dt_4$ are easily expressed, since $z_4\neq 0$.
So we are down to 8 coordinates related to the second and third maps.

\[%
} \right)} \right|}} = \frac{{1 + {y_2} - {t_2}}}{{1 + {y_2}}}\]

After getting rid of the denominators and omitting the negligible terms we get:
\[\begin{array}{l}
{x_3} \approx  - {x_2}\\
{y_3} \approx {y_2}\\
0 \approx {x_2}\\
{t_3} \approx 1 - {t_2}\left( {1 - {y_2}} \right)
\end{array}\]

In terms of differentials:
\[\begin{array}{l}
d{x_3} =  - d{x_2}\\
d{y_3} = d{y_2}\\
0 = d{x_2}\\
d{t_3} =  - d{t_2} + {t_2} \cdot d{y_2}
\end{array}\]

So $dx_2$ and $dx_3$ are zeroes, $dy_3$ and $dt_3$ are expressed from the equations of the second chart, 
so we are down to $dy_2,dz_2,dt_2$ and $dz_3$ and hence the cotangent space is at most 4-dimensional.

\subsection{$\mathcal{E}.1$} 
\begin{center}
	\includegraphics{Case_E1.png}
\end{center}

The local coordinates for related to the charts are:
\[%
\]

We start by relating the first chart to the second:

\[%
} \right)} \right|}} = \frac{{1 - {x_2}{t_2}}}{{\left( {1 - {x_2}} \right){t_2}}}\]

Notice that all the denominators in the last 4 formulas are non-zero, 
when ${x_2},{y_2},{z_2}$ are near zero, and $t_2$ is non-zero.

Hence $dx_1, dy_1. dz_1, dt_1$ can be expressed in terms of coordinates related to the second chart.
So from now on we may disregard the first chart.

Now we connect the coordinates related to the third chart to the coordinates related to the fifth chart.
\[%
} \right)} \right|}} =  - {t_5}\]

Notice that all the denominators in the last 4 formulas are non-zero, 
when $x_5, y_5, z_5, t_5$ are near zero.
Hence $dx_3, dy_3, dz_3, dt_3$ can be expressed in terms of coordinates related to the fifth chart.
So from now on we may disregard the third chart.

So we have only 4 charts to consider: the second, the fourth, the fifth, and the sixth.
We start with relating the fifth chart to the others:

\[%
} \right)} \right|}} = \frac{{{z_5}}}{{{z_5} - 1}}\]

After getting rid of denominators we get:

\[{x_4} \approx {x_5}\]
\[{y_5} \approx {y_2}\]
\[{z_4} \approx {t_5}\]
\[{t_4} \approx  - {z_5}\]

This allows us to express $dz_5, dy_5, dz_5, dt_5$ in terms of the coordinates from second and fourth chart.
Lastly, we connect the coordinates from the second, the fourth, and the sixth charts:

\[%
} \right)} \right|}} = \frac{{{y_6}\left( {{y_6}{z_6} - {t_6}} \right)}}{{{y_6}\left( {1 - {t_6}} \right)}} = \frac{{{y_6}{z_6} - {t_6}}}{{1 - {t_6}}}\]

After getting rid of the denominators, and omitting negligible terms, we get:

\[\begin{array}{c}
{x_2} \approx 0\\
{x_4} \approx {z_4}{y_4}\\
{z_2} = {x_4}\\
{t_2}{y_4} = 1\\
{y_2} \approx  {y_6}\\
0 \approx {x_6} \\
{t_2}{x_6} \approx {t_6}\\
{t_4} \approx {y_6}{z_6} - {t_6}
\end{array}\]

Recall that $t_2, y_4, z_6$ are non-zero, while other coordinates are near zero. 
The differential versions of these are 
\[\begin{array}{c}
d{x_2} = 0\\
d{x_4} = {y_4} \cdot d{z_4}\\
d{z_2} = d{x_4}\\
{t_2}d{y_4} + {y_4}d{t_2} = 0\\
d{y_2} =  d{y_6}\\
0 = d{x_6}\\
{t_2}d{x_6} = d{t_6}\\
d{t_4} = {z_6}d{y_6} - d{t_6}
\end{array}\]

From the first, the sixth, and the seventh equation we see that $dx_2, dx_6, dt_6$ are zeroes. 
From $dz_2,dx_4$ can be expressed can be expressed from the third and the second equations 
by $dz_4$. So cotangent space is spanned by 
\[dy_2,dt_2,dy_4,dz_4,dt_4, dy_6,dz_6,\] 
with relations

\[\begin{array}{c}
{t_2}d{y_4} + {y_4}d{t_2} = 0\\
d{y_2} =  d{y_6}\\
d{t_4} = {z_6}d{y_6}
\end{array}\]

The last two of these three equations allow to express $dt_4$ and $dy_6$ in terms of $dy_2$.
The first equations allows to express $dt_2$ from $dy_4$. So the cotangent space is spanned by
$dy_2,dy_4,dz_4, dz_6$ hence it is at most four-dimensional.

\subsection{$\mathcal{E}.2$} 
\begin{center}
	\includegraphics{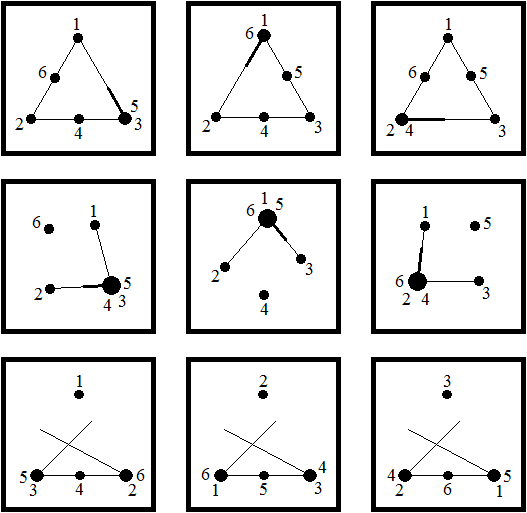}
\end{center}
The points were renumbered to emphasize the symmetry.

We shall not discuss the charts 4-6, as it is obvious all that the first three charts contain all the cross ratios and define the last charts completely.
On the first three charts, we choose the following coordinates:

\[
\]
Notice that since the first three charts are completely degenerate, in the neighborhood of this point 
\[{x_1},{y_1},{z_1},{t_1},{x_2},{y_2},{z_2},{t_2},{x_3},{y_3},{z_3},{t_3} \sim 0.\]

For the last three charts we also choose coordinates (we assume that the extra lines remain the medians of the triangle $P_1P_2P_3$):
\[%
\]

Here ${y_7},{z_7},{t_7},{y_8},{z_8},{t_8},{y_9},{z_9},{t_9} \sim 0$.
The coordinates were chosen in a symmetric way.

From the first chart we may see that 
\[{x_1} = {\left[ {F,C;B,D} \right]_A}\]
\[{y_1} = {\left[ {D,A;B,F} \right]_C}\]
The third coordinate $z_1$ is not a cross ratio, but we can relate it by 
\[{\left[ {C,A;E,F} \right]_D} = \frac{{\left| {\left( {\begin{array}{*{20}{c}}
0&{{z_1}}&0\\
0&{{z_1}{t_1}}&1\\
1&1&1
\end{array}} \right)} \right| \cdot \left| {\left( {\begin{array}{*{20}{c}}
1&1&0\\
0&1&1\\
0&0&1
\end{array}} \right)} \right|}}{{\left| {\left( {\begin{array}{*{20}{c}}
1&{{z_1}}&0\\
0&{{z_1}{t_1}}&1\\
0&1&1
\end{array}} \right)} \right| \cdot \left| {\left( {\begin{array}{*{20}{c}}
0&1&0\\
0&1&1\\
1&0&1
\end{array}} \right)} \right|}} = \frac{{{z_1}}}{{{z_1}{t_1} - 1}}\approx -z_1\]
Since $z_1t_1 \approx 0$ as $z_1,t_1\sim 0$.
Lastly,
\[{t_1} = {\left[ {A,D;E,F} \right]_C}\]

We may relate these formulas to the coordinate of the seventh chart:

\[{x_1} = {\left[ {F,C;B,D} \right]_A} = \frac{{\left| {\left( {%
\]
\[ = \frac{{{y_7} - {t_7}}}{{{y_7} - 1}}\]

Notice also that 
\[{z_7} = {\left[ {C,A;E,F} \right]_B} = {\left[ {\left( \begin{array}{c}
0\\
0\\
1
\end{array} \right),\left( \begin{array}{c}
1\\
{x_3}\\
{y_3}
\end{array} \right);\left( \begin{array}{c}
1\\
0\\
1
\end{array} \right),\left( \begin{array}{c}
1\\
1\\
0
\end{array} \right)} \right]_{\left( {\scriptstyle0\atop
{\scriptstyle1\atop
\scriptstyle0}} \right)}} = \]
\[ = \frac{{\left| {\left( {\begin{array}{*{20}{c}}
0&1\\
1&1
\end{array}} \right)} \right| \cdot \left| {\left( {\begin{array}{*{20}{c}}
1&1\\
{{y_3}}&0
\end{array}} \right)} \right|}}{{\left| {\left( {\begin{array}{*{20}{c}}
0&1\\
1&0
\end{array}} \right)} \right| \cdot \left| {\left( {\begin{array}{*{20}{c}}
1&1\\
{{y_3}}&1
\end{array}} \right)} \right|}} = \frac{{ - {y_3}}}{{1 - {y_3}}}\] 

Multiplying by the denominators and
recalling ${x_1},{y_1},{z_1},{t_1},{y_7},{z_7},{t_7} \sim 0$ we get 
\[\begin{array}{l}
-{x_1} {x_7} \approx {t_7}\\
0 \approx {y_7}\\
0 \approx z_7x_7 \\
-t_1 \approx y_7 - t_7\\
z_7 \approx -y_3.
\end{array}\]

Hence
\[\begin{array}{l}
 - {x_7} \cdot d{x_1} = d{t_7}\\
0 = d{y_7}\\
0 = {x_7} \cdot d{z_7}\\
 - d{t_1} = d{y_7} - d{t_7}\\
dz_7 = -dy_3.
\end{array}\]

Using the second of these four equations, the fourth equation may be rewritten as 
\[d{t_1} = d{t_7}\]

Since our coordinates were chosen symmetrically, we may similarly relate charts 2 and 8, and also charts 3 and 9.
Hence 
\[%
\]

Now we connect charts 1 and 2.

\[{x_1} = {\left[ {F,C;B,D} \right]_A} = \frac{{\left| {\left( {%
\]
\[ =  - \frac{{\left( {1 - {y_2} - {x_2}} \right) \cdot {z_2}\left( {1 - {t_2}} \right)}}
{{1 - {y_2} + {x_2}{z_2}\left( {1 - {t_2}} \right)}}\]

We could also write the equation for $t_1$, but we would not conclude new information from it. 
After multiplying by denominators and recalling $x_i,y_i,z_i,t_i \sim 0$ for $i \in {1,2}$ we see that 
\[\begin{array}{l}
{x_1} =  - {z_2}\\
{y_1} \approx 0\\
{z_1} \approx {z_2}\\
\end{array}\]

Because of cyclic symmetry, similar formulas connect the second chart to the third, and the third to the first.
So together with the former set of formulas and using the symmetry, we get the following set of relations in the cotangent space:

\[%
\]

The cotangent space is spanned by $x_i, y_i, z_i, t_i$, where $i\in \left\{1,2,3,7,8,9\right\}$. 
Anyway, $dy_i=0$ because of the formulas in the second and the seventh lines; it also allows
to simplify the fifth line, which makes the third line redundant.
From the fourth line we may express $dt_1,dt_2,dt_2$, and from the first line we may may express
$dt_7, dt_8,dt_9$. So we may remain with $x_i, z_i$ and relations:

\[%
\]

From the last line we see that $dz_1=dz_2=dz_3$, so we may forget $dz_2$ and $dz_3$ and replace them by $dz_1$.
From the line before the last we conclude that $dx_1=dx_2=dx_3=-dz_1$, so we may forget $dx_1,dx_2,dx_3$ and replace them everywhere by $-dz_1$. From the first line we see that $dz_7=dz_8=dz_9$. 
So the cotangent space is spanned by $dz_1$ and $dx_7,dx_8,dx_9$, and is at most four-dimensional.

\subsection{$\mathcal{E}.3$} 
\begin{center}
	\includegraphics{Case_E3.png}
\end{center}

We choose the coordinates related to the seven charts:

\[%
\]

These charts, except for the last chart, are most degenerate (having no parameters). 
So we have chosen the coordinates so that all of them are zero at this configuration 
and near zero in the neighborhood. In the last chart the only parameter 
that can be non-zero is $x_7$.

Firstly, we shall get rid of the first chart:
\[%
} \right)} \right|}} = \frac{{{t_2} - {x_2}}}{{1 - {x_2}}}\]

Since all denominators are distinct from zero, $dx_1,dy_1,dz_1,dt_1$ can be expressed 
using coordinates related to the second chart.

Secondly, we get rid of the third chart:
\[%
} \right)} \right|}} = \frac{{{t_5}}}{{1 + {t_5}}}\]

Again, all denominators are non-zero, $dx_3,dy_3,dz_3,dt_3$ can be expressed 
using coordinates related to the fifth chart.

Thirdly, we get rid of the sixth chart:
\[%
} \right)} \right|}} = \frac{{{z_5}}}{{{z_5} - 1}}\]

Again, all denominators are non-zero, $dx_6,dy_6,dz_6,dt_6$ can be expressed 
using coordinates related to the fourth and the fifth chart.
Only 4 charts remain: the second, the fourth, the fifth, and the seventh.
Now we shall get rid of the second chart:

\[%
} \right)} \right|}} = \frac{{{x_4} + {z_4}{t_4}}}{{1 + {x_4} + {z_4}{t_4}}}
\approx x_4\]

We see that $dx_2,dy_2,dz_2,dt_2$ are easily expressed in terms of coordinates related to the 
charts 4 and 5. So we remain with 3 charts: the fourth, the fifth and the seventh.
We connect the fifth chart to the fourth:

\[%
} \right)} \right|}} = \frac{{1 - {y_4}}}{{{1 - \left( {{y_4} + {z_4}} \right)} }}\]

After getting rid of the denominators and omitting the negligible terms we get:
\[\begin{array}{c}
{x_5} \approx 0\\
0 \approx  - {y_4}\\
0 \approx {t_4}\\
1 + {t_5} \approx 1 + {z_4}
\end{array}\]

Which means
\[\begin{array}{c}
d{x_5} = 0\\
d{y_4} = 0\\
d{t_4} = 0\\
d{t_5} = d{z_4}
\end{array}\]

To connect to the seventh chart we write:
\[{x_5} = {\left[ {A,E;B,F} \right]_D} = \frac{{\left| {\left( {\begin{array}{*{20}{c}}
0&{{z_7}}&{{t_7}}\\
1&1&{{t_7}}\\
0&{{z_7}}&1
\end{array}} \right)} \right| \cdot \left| {\left( {\begin{array}{*{20}{c}}
1&{{y_7}}&{{t_7}}\\
0&1&{{t_7}}\\
0&{{x_7}}&1
\end{array}} \right)} \right|}}{{\left| {\left( {\begin{array}{*{20}{c}}
0&{{y_7}}&{{t_7}}\\
1&1&{{t_7}}\\
0&{{x_7}}&1
\end{array}} \right)} \right| \cdot \left| {\left( {\begin{array}{*{20}{c}}
1&{{z_7}}&{{t_7}}\\
0&1&{{t_7}}\\
0&{{z_7}}&1
\end{array}} \right)} \right|}} = \]
\[=\frac{{{z_7}\left( {{t_7} - 1} \right)\left( {1 - {x_7}{t_7}} \right)}}{{\left( {{x_7}{t_7} - {y_7}} \right)\left( {1 - {z_7}{t_7}} \right)}}\]


\[{z_5} = {\left[ {D,A;C,F} \right]_E} = \frac{{\left| {\left( {\begin{array}{*{20}{c}}
{{t_7}}&0\\
1&1
\end{array}} \right)} \right| \cdot \left| {\left( {\begin{array}{*{20}{c}}
1&1\\
0&{{x_7}}
\end{array}} \right)} \right|}}{{\left| {\left( {\begin{array}{*{20}{c}}
{{t_7}}&1\\
1&{{x_7}}
\end{array}} \right)} \right| \cdot \left| {\left( {\begin{array}{*{20}{c}}
1&0\\
0&1
\end{array}} \right)} \right|}} = \frac{{{x_7}{t_7}}}{{{x_7}{t_7}-1}}\]

\[1 + {t_5} = {\left[ {D,E;C,F} \right]_A} = \frac{{\left| {\left( {\begin{array}{*{20}{c}}
{{t_7}}&0\\
1&1
\end{array}} \right)} \right| \cdot \left| {\left( {\begin{array}{*{20}{c}}
1&{{y_7}}\\
0&{{x_7}}
\end{array}} \right)} \right|}}{{\left| {\left( {\begin{array}{*{20}{c}}
{{t_7}}&{{y_7}}\\
1&{{x_7}}
\end{array}} \right)} \right| \cdot \left| {\left( {\begin{array}{*{20}{c}}
1&0\\
0&1
\end{array}} \right)} \right|}} = \frac{{{t_7}{x_7}}}{{{t_7}{x_7} - {y_7}}}\]

\[{t_7} = {\left[ {C,E;D,B} \right]_A} = {\left[ {\left( \begin{array}{c}
{z_5}\\
1 + {t_5}\\
1
\end{array} \right),\left( \begin{array}{c}
0\\
1\\
0
\end{array} \right);\left( \begin{array}{c}
0\\
0\\
1
\end{array} \right),\left( \begin{array}{c}
1\\
{x_5}\\
{x_5}{y_5}
\end{array} \right)} \right]_{\left( {\scriptstyle1\atop
{\scriptstyle0\atop
\scriptstyle0}} \right)}} = \]
\[ = \small\frac{{\left| {\left( {\begin{array}{*{20}{c}}
{1 + {t_5}}&0\\
1&1
\end{array}} \right)} \right| \cdot \left| {\left( {\begin{array}{*{20}{c}}
1&{{x_5}}\\
0&{{x_5}{y_5}}
\end{array}} \right)} \right|}}{{\left| {\left( {\begin{array}{*{20}{c}}
{1 + {t_5}}&{{x_5}}\\
1&{{x_5}{y_5}}
\end{array}} \right)} \right| \cdot \left| {\left( {\begin{array}{*{20}{c}}
1&0\\
0&1
\end{array}} \right)} \right|}} = \frac{{\left( {1 + {t_5}} \right) \cdot \left| {\left( {\begin{array}{*{20}{c}}
1&1\\
0&{{y_5}}
\end{array}} \right)} \right|}}{{\left| {\left( {\begin{array}{*{20}{c}}
{1 + {t_5}}&1\\
1&{{y_5}}
\end{array}} \right)} \right|}} = 
\frac{{\left( {1 + {t_5}} \right){y_5}}}{{\left( {1 + {t_5}} \right){y_5} - 1}}\]

From the last 4 formulas we conclude (after multiplying by the denominators and skipping the negligible terms)
\[\begin{array}{l}
0 \approx  - {z_7}\\
{z_5} \approx {x_7}{t_7}\\
{t_7}{x_7} - {y_7} \approx {t_7}{x_7}\\
-t_7 \approx y_5
\end{array}\]

Hence
\[\begin{array}{l}
d{z_7} = 0\\
d{z_5} = {x_7} \cdot d{t_7}\\
d{y_7} = 0 \\
-dt_7=dy_5
\end{array}\]

Previously, we had four more equalities:
\[\begin{array}{c}
d{x_5} = d{y_4} = d{t_4} = 0\\
d{t_5} = d{z_4}
\end{array}\]
So, out of $dx_i, dy_i, dz_i, dt_i$ for $i \in \left\{4,5,7\right\}$, 
five are zeroes, so the cotangent space is spanned by $dx_4, dz_4, dy_5, dz_5,dt_5,dx_7,dt_7$ 
with relations 
\[\begin{array}{l}
d{z_5} = {x_7} \cdot d{t_7}\\
-dt_7=dy_5\\
d{t_5} = d{z_4}
\end{array}\]
which express all the remaining coordinates from the fifth chart from the other charts, 
so the cotangent space is spanned by $dx_4, dz_4,dx_7, dt_7$.

\subsection{$\mathcal{E}.6$}
\begin{center}
	\includegraphics{Case_E6.png}
\end{center}
\[%
\]
The only coordinate which may be far from zero is $t_5$, which governs 
the direction of the line via the double point.

We start by eliminating the coordinates of the first chart:
\[{x_1} = {\left[ {C,E;D,B} \right]_A} = {\left[ {\left( \begin{array}{c}
0\\
0\\
1
\end{array} \right),\left( \begin{array}{c}
{x_2}\\
1 + {y_2}\\
1
\end{array} \right);\left( \begin{array}{c}
0\\
1\\
0
\end{array} \right),\left( \begin{array}{l}
1\\
1\\
1
\end{array} \right)} \right]_{\left( {\scriptstyle1\atop
{\scriptstyle0\atop
\scriptstyle0}} \right)}} = \]

\[ = \frac{{\left| {\left( {\begin{array}{*{20}{c}}
0&1\\
1&0
\end{array}} \right)} \right| \cdot \left| {\left( {\begin{array}{*{20}{c}}
{1 + {y_2}}&1\\
1&1
\end{array}} \right)} \right|}}{{\left| {\left( {\begin{array}{*{20}{c}}
0&1\\
1&1
\end{array}} \right)} \right| \cdot \left| {\left( {\begin{array}{*{20}{c}}
{1 + {y_2}}&1\\
1&0
\end{array}} \right)} \right|}} =  - {y_2}\]

\[{y_1} = {\left[ {E,A;D,B} \right]_C} = {\left[ {\left( \begin{array}{c}
{x_2}\\
1 + {y_2}\\
1
\end{array} \right),\left( \begin{array}{c}
1\\
0\\
0
\end{array} \right);\left( \begin{array}{c}
0\\
1\\
0
\end{array} \right),\left( \begin{array}{l}
1\\
1\\
1
\end{array} \right)} \right]_{\left( {\scriptstyle0\atop
{\scriptstyle0\atop
\scriptstyle1}} \right)}} = \]
\[ = \frac{{\left| {\left( {\begin{array}{*{20}{c}}
{{x_2}}&0\\
{1 + {y_2}}&1
\end{array}} \right)} \right| \cdot \left| {\left( {\begin{array}{*{20}{c}}
1&1\\
0&1
\end{array}} \right)} \right|}}{{\left| {\left( {\begin{array}{*{20}{c}}
{{x_2}}&1\\
{1 + {y_2}}&1
\end{array}} \right)} \right| \cdot \left| {\left( {\begin{array}{*{20}{c}}
1&0\\
0&1
\end{array}} \right)} \right|}} = \frac{{{x_2}}}{{{x_2} - {y_2} - 1}}\]

\[{z_1} = {\left[ {C,E;F,B} \right]_A} = {\left[ {\left( \begin{array}{c}
0\\
0\\
1
\end{array} \right),\left( \begin{array}{c}
{x_2}\\
1 + {y_2}\\
1
\end{array} \right);\left( \begin{array}{c}
{z_2}{t_2}\\
{z_2}\\
1
\end{array} \right),\left( \begin{array}{l}
1\\
1\\
1
\end{array} \right)} \right]_{\left( {\scriptstyle1\atop
{\scriptstyle0\atop
\scriptstyle0}} \right)}} = \]
\[ = \frac{{\left| {\left( {\begin{array}{*{20}{c}}
0&{{z_2}}\\
1&1
\end{array}} \right)} \right| \cdot \left| {\left( {\begin{array}{*{20}{c}}
{1 + {y_2}}&1\\
1&1
\end{array}} \right)} \right|}}{{\left| {\left( {\begin{array}{*{20}{c}}
0&1\\
1&1
\end{array}} \right)} \right| \cdot \left| {\left( {\begin{array}{*{20}{c}}
{1 + {y_2}}&{{z_2}}\\
1&1
\end{array}} \right)} \right|}} = \frac{{{z_2} \cdot {y_2}}}{{1 + {y_2} - {z_2}}}\]

\[{t_1} = {\left[ {E,A;F,B} \right]_C} = {\left[ {\left( \begin{array}{c}
{x_2}\\
1 + {y_2}\\
1
\end{array} \right),\left( \begin{array}{c}
1\\
0\\
0
\end{array} \right);\left( \begin{array}{c}
{z_2}{t_2}\\
{z_2}\\
1
\end{array} \right),\left( \begin{array}{l}
1\\
1\\
1
\end{array} \right)} \right]_{\left( {\scriptstyle0\atop
{\scriptstyle0\atop
\scriptstyle1}} \right)}} = \]
\[ = \frac{{\left| {\left( {\begin{array}{*{20}{c}}
{{x_2}}&{{z_2}{t_2}}\\
{1 + {y_2}}&{{z_2}}
\end{array}} \right)} \right| \cdot \left| {\left( {\begin{array}{*{20}{c}}
1&1\\
0&1
\end{array}} \right)} \right|}}{{\left| {\left( {\begin{array}{*{20}{c}}
{{x_2}}&1\\
{1 + {y_2}}&1
\end{array}} \right)} \right| \cdot \left| {\left( {\begin{array}{*{20}{c}}
1&{{z_2}{t_2}}\\
0&{{z_2}}
\end{array}} \right)} \right|}} = \frac{{{x_2} - \left( {1 + {y_2}} \right){t_2}}}{{{x_2} - 1 - {y_2}}}\] 

All denominators are non-zero, so the coordinates of the first chart are not needed.
Now similarly for the third chart:
\[{x_3} = {\left[ {A,E;B,D} \right]_C} = {\left[ {\left( \begin{array}{c}
1\\
0\\
0
\end{array} \right),\left( \begin{array}{c}
{x_2}\\
1 + {y_2}\\
1
\end{array} \right);\left( \begin{array}{l}
1\\
1\\
1
\end{array} \right),\left( \begin{array}{c}
0\\
1\\
0
\end{array} \right)} \right]_{\left( {\scriptstyle0\atop
{\scriptstyle0\atop
\scriptstyle1}} \right)}} = \]
\[ = \frac{{\left| {\left( {\begin{array}{*{20}{c}}
1&1\\
0&1
\end{array}} \right)} \right| \cdot \left| {\left( {\begin{array}{*{20}{c}}
{{x_2}}&0\\
{1 + {y_2}}&1
\end{array}} \right)} \right|}}{{\left| {\left( {\begin{array}{*{20}{c}}
1&0\\
0&1
\end{array}} \right)} \right| \cdot \left| {\left( {\begin{array}{*{20}{c}}
{{x_2}}&1\\
{1 + {y_2}}&1
\end{array}} \right)} \right|}} = \frac{{{x_2}}}{{{x_2} - 1 - {y_2}}}\]

\[{y_3} = {\left[ {E,C;B,D} \right]_A} = {\left[ {\left( \begin{array}{c}
{x_2}\\
1 + {y_2}\\
1
\end{array} \right),\left( \begin{array}{c}
0\\
0\\
1
\end{array} \right);\left( \begin{array}{c}
1\\
1\\
1
\end{array} \right),\left( \begin{array}{c}
0\\
1\\
0
\end{array} \right)} \right]_{\left( {\scriptstyle1\atop
{\scriptstyle0\atop
\scriptstyle0}} \right)}} = \]
\[ = \frac{{\left| {\left( {\begin{array}{*{20}{c}}
{1 + {y_2}}&1\\
1&1
\end{array}} \right)} \right| \cdot \left| {\left( {\begin{array}{*{20}{c}}
0&1\\
1&0
\end{array}} \right)} \right|}}{{\left| {\left( {\begin{array}{*{20}{c}}
{1 + {y_2}}&1\\
1&0
\end{array}} \right)} \right| \cdot \left| {\left( {\begin{array}{*{20}{c}}
0&1\\
1&1
\end{array}} \right)} \right|}} =  - {y_2}\]

\[{z_3} = {\left[ {C,A;F,D} \right]_E} = {\left[ {\left( \begin{array}{c}
0\\
0\\
1
\end{array} \right),\left( \begin{array}{c}
1\\
1 + {x_6}\\
{y_6}
\end{array} \right);\left( \begin{array}{l}
1\\
1\\
1
\end{array} \right),\left( \begin{array}{c}
0\\
1\\
0
\end{array} \right)} \right]_{\left( {\scriptstyle1\atop
{\scriptstyle0\atop
\scriptstyle0}} \right)}} = \]
\[ = \frac{{\left| {\left( {\begin{array}{*{20}{c}}
0&1\\
1&1
\end{array}} \right)} \right| \cdot \left| {\left( {\begin{array}{*{20}{c}}
{1 + {x_6}}&1\\
{{y_6}}&0
\end{array}} \right)} \right|}}{{\left| {\left( {\begin{array}{*{20}{c}}
0&1\\
1&0
\end{array}} \right)} \right| \cdot \left| {\left( {\begin{array}{*{20}{c}}
{1 + {x_6}}&1\\
{{y_6}}&1
\end{array}} \right)} \right|}} = \frac{{{y_6}}}{{{y_6} - 1 - {x_6}}}\]

\[{t_3} = {\left[ {A,E;F,D} \right]_C} = {\left[ {\left( \begin{array}{c}
1\\
0\\
0
\end{array} \right),\left( \begin{array}{c}
{z_4}\\
1\\
{z_4}{t_4}
\end{array} \right);\left( \begin{array}{l}
1\\
1\\
1
\end{array} \right),\left( \begin{array}{c}
0\\
1\\
0
\end{array} \right)} \right]_{\left( {\scriptstyle0\atop
{\scriptstyle0\atop
\scriptstyle1}} \right)}} = \]
\[ = \frac{{\left| {\left( {\begin{array}{*{20}{c}}
1&1\\
0&1
\end{array}} \right)} \right| \cdot \left| {\left( {\begin{array}{*{20}{c}}
{{z_4}}&0\\
1&1
\end{array}} \right)} \right|}}{{\left| {\left( {\begin{array}{*{20}{c}}
1&0\\
0&1
\end{array}} \right)} \right| \cdot \left| {\left( {\begin{array}{*{20}{c}}
{{z_4}}&1\\
1&1
\end{array}} \right)} \right|}} = \frac{{{z_4}}}{{{z_4} - 1}}\]
Again, all denominators are non-zero, so the third chart might also be omitted.
Now, the fourth chart:
\[{x_4} = {\left[ {A,D;B,F} \right]_C} = {\left[ {\left( \begin{array}{c}
1\\
0\\
0
\end{array} \right),\left( \begin{array}{c}
0\\
1\\
0
\end{array} \right);\left( \begin{array}{l}
1\\
1\\
1
\end{array} \right),\left( \begin{array}{c}
{z_2}{t_2}\\
{z_2}\\
1
\end{array} \right)} \right]_{\left( {\scriptstyle0\atop
{\scriptstyle0\atop
\scriptstyle1}} \right)}} = \]
\[ = \frac{{\left| {\left( {\begin{array}{*{20}{c}}
1&1\\
0&1
\end{array}} \right)} \right| \cdot \left| {\left( {\begin{array}{*{20}{c}}
0&{{z_2}{t_2}}\\
1&{{z_2}}
\end{array}} \right)} \right|}}{{\left| {\left( {\begin{array}{*{20}{c}}
1&{{z_2}{t_2}}\\
0&{{z_2}}
\end{array}} \right)} \right| \cdot \left| {\left( {\begin{array}{*{20}{c}}
0&1\\
1&1
\end{array}} \right)} \right|}} =  - \frac{{\left| {\left( {\begin{array}{*{20}{c}}
0&{{t_2}}\\
1&1
\end{array}} \right)} \right|}}{{\left| {\left( {\begin{array}{*{20}{c}}
1&{{t_2}}\\
0&1
\end{array}} \right)} \right|}} = {t_2}\]

\[{y_4} = {\left[ {D,C;B,F} \right]_A} = {\left[ {\left( \begin{array}{c}
0\\
1\\
0
\end{array} \right),\left( \begin{array}{c}
0\\
0\\
1
\end{array} \right);\left( \begin{array}{l}
1\\
1\\
1
\end{array} \right),\left( \begin{array}{c}
{z_2}{t_2}\\
{z_2}\\
1
\end{array} \right)} \right]_{\left( {\scriptstyle1\atop
{\scriptstyle0\atop
\scriptstyle0}} \right)}} = \]
\[ = \frac{{\left| {\left( {\begin{array}{*{20}{c}}
1&1\\
0&1
\end{array}} \right)} \right| \cdot \left| {\left( {\begin{array}{*{20}{c}}
0&{{z_2}}\\
1&1
\end{array}} \right)} \right|}}{{\left| {\left( {\begin{array}{*{20}{c}}
1&{{z_2}}\\
0&1
\end{array}} \right)} \right| \cdot \left| {\left( {\begin{array}{*{20}{c}}
0&1\\
1&1
\end{array}} \right)} \right|}} = {z_2}\]

\[{z_4} = {\left[ {D,A;E,F} \right]_C} = {\left[ {\left( \begin{array}{c}
0\\
1\\
0
\end{array} \right),\left( \begin{array}{c}
1\\
1 + {x_6}\\
{y_6}
\end{array} \right);\left( \begin{array}{c}
1\\
0\\
0
\end{array} \right),\left( \begin{array}{l}
1\\
1\\
1
\end{array} \right)} \right]_{\left( {\scriptstyle0\atop
{\scriptstyle0\atop
\scriptstyle1}} \right)}} = \]
\[ = \frac{{\left| {\left( {\begin{array}{*{20}{c}}
0&1\\
1&0
\end{array}} \right)} \right| \cdot \left| {\left( {\begin{array}{*{20}{c}}
1&1\\
{1 + {x_6}}&1
\end{array}} \right)} \right|}}{{\left| {\left( {\begin{array}{*{20}{c}}
0&1\\
1&1
\end{array}} \right)} \right| \cdot \left| {\left( {\begin{array}{*{20}{c}}
1&1\\
{1 + {x_6}}&0
\end{array}} \right)} \right|}} = \frac{{{x_6}}}{{1 + {x_6}}}\]

\[{t_4} = {\left[ {A,C;E,F} \right]_D} = {\left[ {\left( \begin{array}{c}
1\\
1 + {x_6}\\
{y_6}
\end{array} \right),\left( \begin{array}{c}
0\\
0\\
1
\end{array} \right);\left( \begin{array}{c}
1\\
0\\
0
\end{array} \right),\left( \begin{array}{l}
1\\
1\\
1
\end{array} \right)} \right]_{\left( {\scriptstyle0\atop
{\scriptstyle1\atop
\scriptstyle0}} \right)}} = \]
\[ = \frac{{\left| {\left( {\begin{array}{*{20}{c}}
1&1\\
{{y_6}}&0
\end{array}} \right)} \right| \cdot \left| {\left( {\begin{array}{*{20}{c}}
0&1\\
1&1
\end{array}} \right)} \right|}}{{\left| {\left( {\begin{array}{*{20}{c}}
1&1\\
{{y_6}}&1
\end{array}} \right)} \right| \cdot \left| {\left( {\begin{array}{*{20}{c}}
0&1\\
1&0
\end{array}} \right)} \right|}} = \frac{{ - {y_6}}}{{1 - {y_6}}}\]
And here again, all denominators are non-zero.

We remain with just 3 charts: the second, the fifth, and the sixth.
\[{t_2} = {\left[ {D,A;F,B} \right]_C} = {\left[ {\left( \begin{array}{c}
1\\
1 + {x_5}\\
{y_5}
\end{array} \right),\left( \begin{array}{c}
1\\
0\\
0
\end{array} \right);\left( \begin{array}{l}
1\\
1\\
1
\end{array} \right),\left( \begin{array}{c}
0\\
1\\
0
\end{array} \right)} \right]_{\left( {\scriptstyle0\atop
{\scriptstyle0\atop
\scriptstyle1}} \right)}} = \]
\[ = \frac{{\left| {\left( {\begin{array}{*{20}{c}}
1&1\\
{1 + {x_5}}&1
\end{array}} \right)} \right| \cdot \left| {\left( {\begin{array}{*{20}{c}}
1&0\\
0&1
\end{array}} \right)} \right|}}{{\left| {\left( {\begin{array}{*{20}{c}}
1&0\\
{1 + {x_5}}&1
\end{array}} \right)} \right| \cdot \left| {\left( {\begin{array}{*{20}{c}}
1&1\\
0&1
\end{array}} \right)} \right|}} =  - {x_5}\]
\[{z_2} = {\left[ {C,D;F,B} \right]_A} = {\left[ {\left( \begin{array}{c}
0\\
0\\
1
\end{array} \right),\left( \begin{array}{c}
1\\
1 + {x_5}\\
{y_5}
\end{array} \right);\left( \begin{array}{l}
1\\
1\\
1
\end{array} \right),\left( \begin{array}{c}
0\\
1\\
0
\end{array} \right)} \right]_{\left( {\scriptstyle1\atop
{\scriptstyle0\atop
\scriptstyle0}} \right)}} = \]

\[ = \frac{{\left| {\left( {\begin{array}{*{20}{c}}
0&1\\
1&1
\end{array}} \right)} \right| \cdot \left| {\left( {\begin{array}{*{20}{c}}
{1 + {x_5}}&1\\
{{y_5}}&0
\end{array}} \right)} \right|}}{{\left| {\left( {\begin{array}{*{20}{c}}
0&1\\
1&0
\end{array}} \right)} \right| \cdot \left| {\left( {\begin{array}{*{20}{c}}
{1 + {x_5}}&1\\
{{y_5}}&1
\end{array}} \right)} \right|}} = \frac{{{y_5}}}{{{y_5} - 1 - {x_5}}} \approx  - {y_5}\]

On the other hand, 
\[{z_2} = {\left[ {C,D;F,B} \right]_A} = {\left[\small {\left( \begin{array}{c}
0\\
0\\
1
\end{array} \right),\left( \begin{array}{c}
0\\
1\\
0
\end{array} \right);\left( \begin{array}{l}
1\\
1\\
1
\end{array} \right),\left( \begin{array}{c}
1\\
1 + {x_6} + {z_6}\\
{y_6} + {z_6}{t_6}
\end{array} \right)} \right]_{\left( {\scriptstyle1\atop
{\scriptstyle1 + {x_6}\atop
\scriptstyle{y_6}}} \right)}} = \]
\[ = \frac{{\left| {\left( {\begin{array}{*{20}{c}}
0&1&1\\
0&1&{1 + {x_6}}\\
1&1&{{y_6}}
\end{array}} \right)} \right| \cdot \left| {\left( {\begin{array}{*{20}{c}}
0&1&1\\
1&{1 + {x_6} + {z_6}}&{1 + {x_6}}\\
0&{{y_6} + {z_6}{t_6}}&{{y_6}}
\end{array}} \right)} \right|}}{{\left| {\left( {\begin{array}{*{20}{c}}
0&1&1\\
0&{1 + {x_6} + {z_6}}&{1 + {x_6}}\\
1&{{y_6} + {z_6}{t_6}}&{{y_6}}
\end{array}} \right)} \right| \cdot \left| {\left( {\begin{array}{*{20}{c}}
0&1&1\\
1&1&{1 + {x_6}}\\
0&1&{{y_6}}
\end{array}} \right)} \right|}} = \]
\[ = \frac{{\left| {\left( {\begin{array}{*{20}{c}}
1&1\\
1&{1 + {x_6}}
\end{array}} \right)} \right| \cdot \left| {\left( {\begin{array}{*{20}{c}}
0&1\\
{{z_6}{t_6}}&{{y_6}}
\end{array}} \right)} \right|}}{{\left| {\left( {\begin{array}{*{20}{c}}
0&1\\
{{z_6}}&{1 + {x_6}}
\end{array}} \right)} \right| \cdot \left| {\left( {\begin{array}{*{20}{c}}
1&1\\
1&{{y_6}}
\end{array}} \right)} \right|}} = \frac{{{x_6}{t_6}}}{{{y_6} - 1}} \approx 0\]

\[t_2 = {\left[ {D,A;F,B} \right]_C} = {\left[\footnotesize {\left( \begin{array}{c}
0\\
1\\
0
\end{array} \right),\left( \begin{array}{c}
1\\
1 + {x_6}\\
{y_6}
\end{array} \right);\left( \begin{array}{l}
1\\
1\\
1
\end{array} \right),\left( \begin{array}{c}
1\\
1 + {x_6} + {z_6}\\
{y_6} + {z_6}{t_6}
\end{array} \right)} \right]_{\left( {\scriptstyle0\atop
{\scriptstyle0\atop
\scriptstyle1}} \right)}} = \]
\[ = \frac{{\left| {\left( {\begin{array}{*{20}{c}}
0&1\\
1&1
\end{array}} \right)} \right| \cdot \left| {\left( {\begin{array}{*{20}{c}}
1&1\\
{1 + {x_6}}&{1 + {x_6} + {z_6}}
\end{array}} \right)} \right|}}{{\left| {\left( {\begin{array}{*{20}{c}}
0&1\\
1&{1 + {x_6} + {z_6}}
\end{array}} \right)} \right| \cdot \left| {\left( {\begin{array}{*{20}{c}}
1&1\\
{1 + {x_6}}&1
\end{array}} \right)} \right|}} = \frac{{\left| {\left( {\begin{array}{*{20}{c}}
1&0\\
{1 + {x_6}}&{{z_6}}
\end{array}} \right)} \right|}}{{\left( {\begin{array}{*{20}{c}}
1&1\\
{{x_6}}&0
\end{array}} \right)}} =  - \frac{{{z_6}}}{{{x_6}}}\]
Now we shall compute ${\left[ {A,E;B,D} \right]_C}$ in all the remaining charts.
\[{\left[ {\left( \begin{array}{c}
1\\
1 + {x_6}\\
{y_6}
\end{array} \right),\left( \begin{array}{c}
1\\
0\\
0
\end{array} \right);\left( \begin{array}{c}
1\\
1 + {x_6} + {z_6}\\
{y_6} + {z_6}{t_6}
\end{array} \right),\left( \begin{array}{c}
0\\
1\\
0
\end{array} \right)} \right]_{\left( {\scriptstyle0\atop
{\scriptstyle0\atop
\scriptstyle1}} \right)}} = \]
\[ = \frac{{\left| {\left( {\begin{array}{*{20}{c}}
1&1\\
{1 + {x_6}}&{1 + {x_6} + {z_6}}
\end{array}} \right)} \right| \cdot \left| {\left( {\begin{array}{*{20}{c}}
1&0\\
0&1
\end{array}} \right)} \right|}}{{\left| {\left( {\begin{array}{*{20}{c}}
1&0\\
{1 + {x_6}}&1
\end{array}} \right)} \right| \cdot \left| {\left( {\begin{array}{*{20}{c}}
1&1\\
0&{1 + {x_6} + {z_6}}
\end{array}} \right)} \right|}} = \frac{{{z_6}}}{{1 + {x_6} + {z_6}}} \approx {z_6}\]

\[{\left[ {A,E;B,D} \right]_C} = {\left[ {\left( \begin{array}{c}
1\\
0\\
0
\end{array} \right),\left( \begin{array}{c}
{x_2}\\
1 + {y_2}\\
1
\end{array} \right);\left( \begin{array}{l}
1\\
1\\
1
\end{array} \right),\left( \begin{array}{c}
0\\
1\\
0
\end{array} \right)} \right]_{\left( {\scriptstyle0\atop
{\scriptstyle0\atop
\scriptstyle1}} \right)}} = \]
\[ = \frac{{\left| {\left( {\begin{array}{*{20}{c}}
1&1\\
0&1
\end{array}} \right)} \right| \cdot \left| {\left( {\begin{array}{*{20}{c}}
{{x_2}}&0\\
{1 + {y_2}}&1
\end{array}} \right)} \right|}}{{\left| {\left( {\begin{array}{*{20}{c}}
1&0\\
0&1
\end{array}} \right)} \right| \cdot \left| {\left( {\begin{array}{*{20}{c}}
{{x_2}}&1\\
{1 + {y_2}}&1
\end{array}} \right)} \right|}} = \frac{{{x_2}}}{{{x_2} - 1 - {y_2}}} \approx  - {x_2}\]

\[ {\left[ {A,E;B,D} \right]_C} = {\left[ {\left( \begin{array}{c}
1\\
0\\
0
\end{array} \right),\left( \begin{array}{c}
1\\
1 + {x_5} + {z_5}\\
{y_5} + {z_5}{t_5}
\end{array} \right);\left( \begin{array}{c}
0\\
1\\
0
\end{array} \right),\left( \begin{array}{c}
1\\
1 + {x_5}\\
{y_5}
\end{array} \right)} \right]_{\left( {\scriptstyle0\atop
{\scriptstyle0\atop
\scriptstyle1}} \right)}} = \]
\[ = \frac{{\left| {\left( {\begin{array}{*{20}{c}}
1&0\\
0&1
\end{array}} \right)} \right| \cdot \left| {\left( {\begin{array}{*{20}{c}}
1&1\\
{1 + {x_5} + {z_5}}&{1 + {x_5}}
\end{array}} \right)} \right|}}{{\left| {\left( {\begin{array}{*{20}{c}}
1&1\\
0&{1 + {x_5}}
\end{array}} \right)} \right| \cdot \left| {\left( {\begin{array}{*{20}{c}}
1&0\\
{1 + {x_5} + {z_5}}&1
\end{array}} \right)} \right|}} = \small\frac{{\left| {\left(\small {\begin{array}{*{20}{c}}
0&1\\
{{z_5}}&{1 + {x_5}}
\end{array}} \right)} \right|}}{{1 + {x_5}}} =\]
\[= \frac{{-{z_5}}}{{1 + {x_5}}} \approx {-z_5}\]
Therefore ${z_6} \approx  - {x_2} \approx -{z_5}$.
\[1 + {y_2} = {\left[ {C,D;E,B} \right]_A} = {\footnotesize\left[ {\left( \begin{array}{c}
0\\
0\\
1
\end{array} \right),\left( \begin{array}{c}
0\\
1\\
0
\end{array} \right);\left( \begin{array}{c}
1\\
0\\
0
\end{array} \right),\left( \begin{array}{c}
1\\
1 + {x_6} + {z_6}\\
{y_6} + {z_6}{t_6}
\end{array} \right)} \right]_{\left( {\scriptstyle1\atop
{\scriptstyle1 + {x_6}\atop
\scriptstyle{y_6}}} \right)}} = \]
\[ = \frac{{\left| {\left( {\begin{array}{*{20}{c}}
0&1&1\\
0&0&{1 + {x_6}}\\
1&0&{{y_6}}
\end{array}} \right)} \right| \cdot \left| {\left( {\begin{array}{*{20}{c}}
0&1&1\\
1&{1 + {x_6} + {z_6}}&{1 + {x_6}}\\
0&{{y_6} + {z_6}{t_6}}&{{y_6}}
\end{array}} \right)} \right|}}{{\left| {\left( {\begin{array}{*{20}{c}}
0&1&1\\
0&{1 + {x_6} + {z_6}}&{1 + {x_6}}\\
1&{{y_6} + {z_6}{t_6}}&{{y_6}}
\end{array}} \right)} \right| \cdot \left| {\left( {\begin{array}{*{20}{c}}
0&1&1\\
1&0&{1 + {x_6}}\\
0&0&{{y_6}}
\end{array}} \right)} \right|}} = \]
\[ = \frac{{\left( {1 + {x_6}} \right) \cdot \left| {\left( {\begin{array}{*{20}{c}}
1&1\\
{{y_6} + {z_6}{t_6}}&{{y_6}}
\end{array}} \right)} \right|}}{{\left| {\left( {\begin{array}{*{20}{c}}
1&1\\
{1 + {x_6} + {z_6}}&{1 + {x_6}}
\end{array}} \right)} \right| \cdot {y_6}}} = \frac{{\left( {1 + {x_6}} \right) \cdot \left| {\left( {\begin{array}{*{20}{c}}
0&1\\
{{z_6}{t_6}}&{{y_6}}
\end{array}} \right)} \right|}}{{\left| {\left( {\begin{array}{*{20}{c}}
0&1\\
{{z_6}}&{1 + {x_6}}
\end{array}} \right)} \right| \cdot {y_6}}} = \frac{{\left( {1 + {x_6}} \right){t_6}}}{{{y_6}}}\]
In other words $y_6\approx t_6$.

Now we compute ${\left[ {B,E;C,D} \right]_A}$ in the coordinates of the second and the sixth charts:
\[{\left[ {B,E;C,D} \right]_A} = {\left[ {\left( \begin{array}{l}
1\\
1\\
1
\end{array} \right),\left( \begin{array}{c}
{x_2}\\
1 + {y_2}\\
1
\end{array} \right);\left( \begin{array}{c}
0\\
0\\
1
\end{array} \right),\left( \begin{array}{c}
0\\
1\\
0
\end{array} \right)} \right]_{\left( {\scriptstyle1\atop
{\scriptstyle0\atop
\scriptstyle0}} \right)}} = \]
\[ = \frac{{\left| {\left( {\begin{array}{*{20}{c}}
1&0\\
1&1
\end{array}} \right)} \right| \cdot \left| {\left( {\begin{array}{*{20}{c}}
{{x_2}}&1\\
{1 + {y_2}}&0
\end{array}} \right)} \right|}}{{\left| {\left( {\begin{array}{*{20}{c}}
1&1\\
1&0
\end{array}} \right)} \right| \cdot \left| {\left( {\begin{array}{*{20}{c}}
{{x_2}}&0\\
{1 + {y_2}}&1
\end{array}} \right)} \right|}} = \frac{{1 + {y_2}}}{{{x_2}}}\]

\[{\left[ {B,E;C,D} \right]_A} = {\left[ {\left( \begin{array}{c}
1\\
1 + {x_6} + {z_6}\\
{y_6} + {z_6}{t_6}
\end{array} \right),\left( \begin{array}{c}
1\\
0\\
0
\end{array} \right);\left( \begin{array}{c}
0\\
0\\
1
\end{array} \right),\left( \begin{array}{c}
0\\
1\\
0
\end{array} \right)} \right]_{\left( {\scriptstyle1\atop
{\scriptstyle1 + {x_6}\atop
\scriptstyle{y_6}}} \right)}} = \]
\[ = \frac{{\left| {\left( {\begin{array}{*{20}{c}}
1&0&1\\
{1 + {x_6} + {z_6}}&0&{1 + {x_6}}\\
{{y_6} + {z_6}{t_6}}&1&{{y_6}}
\end{array}} \right)} \right| \cdot \left| {\left( {\begin{array}{*{20}{c}}
1&0&1\\
0&1&{1 + {x_6}}\\
0&0&{{y_6}}
\end{array}} \right)} \right|}}{{\left| {\left( {\begin{array}{*{20}{c}}
1&0&1\\
{1 + {x_6} + {z_6}}&1&{1 + {x_6}}\\
{{y_6} + {z_6}{t_6}}&0&{{y_6}}
\end{array}} \right)} \right| \cdot \left| {\left( {\begin{array}{*{20}{c}}
1&0&1\\
0&0&{1 + {x_6}}\\
0&1&{{y_6}}
\end{array}} \right)} \right|}} = \]
\[ = \frac{{\left| {\left( {\begin{array}{*{20}{c}}
0&0&1\\
{{z_6}}&0&{1 + {x_6}}\\
{{z_6}{t_6}}&1&{{y_6}}
\end{array}} \right)} \right| \cdot {y_6}}}{{ - \left| {\left( {\begin{array}{*{20}{c}}
0&0&1\\
{{z_6}}&1&{1 + {x_6}}\\
{{z_6}{t_6}}&0&{{y_6}}
\end{array}} \right)} \right| \cdot \left( {1 + {x_6}} \right)}} = \frac{{\left| {\left( {\begin{array}{*{20}{c}}
0&1\\
{{z_6}}&{1 + {x_6}}
\end{array}} \right)} \right| \cdot {y_6}}}{{\left| {\left( {\begin{array}{*{20}{c}}
0&1\\
{{z_6}{t_6}}&{{y_6}}
\end{array}} \right)} \right| \cdot \left( {1 + {x_6}} \right)}} = \]
\[ = \frac{{\left| {\left( {\begin{array}{*{20}{c}}
0&1\\
1&{1 + {x_6}}
\end{array}} \right)} \right| \cdot {y_6}}}{{\left| {\left( {\begin{array}{*{20}{c}}
0&1\\
{{t_6}}&{{y_6}}
\end{array}} \right)} \right| \cdot \left( {1 + {x_6}} \right)}} = 
\frac{{{y_6}}}{{{t_6}\left( {1 + {x_6}} \right)}}\]

Therefore:
\[\frac{{{y_6}}}{{{t_6}\left( {1 + {x_6}} \right)}} = \frac{{1 + {y_2}}}{{{x_2}}}\]
\[0 \approx {y_6}{x_2} = {t_6}\left( {1 + {x_6}} \right)\left( {1 + {y_2}} \right) \approx {t_6}\]

To summarize: 
\[\begin{array}{c}
{t_2} =  - {x_5}\\
 - {y_5} \approx {z_2} \approx 0\\
0 \approx  - {x_6}{t_2} \approx {z_6}\\
{z_6} \approx  - {x_2} \approx - {z_5}\\
{y_6} \approx {t_6}\\
{t_6} \approx 0
\end{array}\]
After simplification: 
\[\begin{array}{c}
{t_2} =  - {x_5}\\
{y_5} \approx {z_2} \approx 0 \approx {z_6} \approx  {x_2} \approx {z_5}  \\
{y_6} \approx {t_6} \approx 0
\end{array}\]
Hence
\[dt_2=-dx_5\]
\[0= dy_5 = dz_2 = dz_6 = dx_2= dz_5=dy_6=dt_6\]
Therefore $dx_2,dz_2, dy_5,dz_5,dy_6,dz_6,dt_6$ are zeroes, and $dt_2$ is expressed in terms of 
$dx_5$.
So the cotangent space is spanned by $dy_2,dx_5, dt_5,dx_6$ and is at most four-dimensional.

\subsection{$\mathcal{F}.1$}
We shall place the charts in different order, than what came up during the classification, to emphasize the symmetry of this configuration:
\begin{center}
	\includegraphics{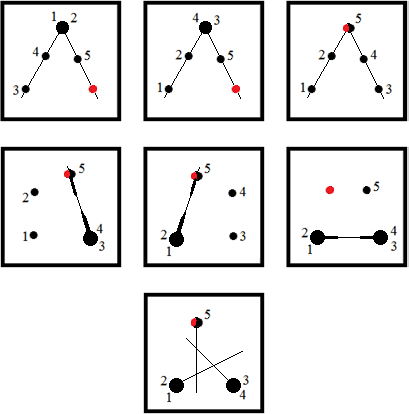}
\end{center}

We choose the coordinates:
\[\small
%
\]

We shall show that the fourth chart is redundant; 
because of the symmetry, so are the fifth and the sixth chart, 
so we won't introduce coordinates for those charts.

For the seventh chart all coordinates are zero, except for $t_7$. 
If $t_7$ becomes zero (or symmetrically, $\infty$) we get an even more degenerate 
collection of charts which we shall discuss later. 
So, for now, all coordinates except $t_7$ are zeroes, while $t_7\neq 0$.

We have
\[{x_4} = {\left[ {C,E;D,B} \right]_A} = {\left[ {\left( \begin{array}{c}
0\\
0\\
1
\end{array} \right),\left( \begin{array}{c}
1\\
1 + {x_3}\\
{y_3}
\end{array} \right);\left( \begin{array}{l}
1\\
1\\
1
\end{array} \right),\left( \begin{array}{c}
0\\
1\\
0
\end{array} \right)} \right]_{\left( {\scriptstyle1\atop
{\scriptstyle0\atop
\scriptstyle0}} \right)}} = \]
\[ = \frac{{\left| {\left( {\begin{array}{*{20}{c}}
0&1\\
1&1
\end{array}} \right)} \right| \cdot \left| {\left( {\begin{array}{*{20}{c}}
{1 + {x_3}}&1\\
{{y_3}}&0
\end{array}} \right)} \right|}}{{\left| {\left( {\begin{array}{*{20}{c}}
0&1\\
1&0
\end{array}} \right)} \right| \cdot \left| {\left( {\begin{array}{*{20}{c}}
{1 + {x_3}}&1\\
{{y_3}}&1
\end{array}} \right)} \right|}} =  - \frac{{{y_3}}}{{1 + {x_3} - {y_3}}}\]

\[{y_4} = {\left[ {E,A;D,B} \right]_C} = {\left[ {\left( \begin{array}{c}
1\\
1 + {x_3}\\
{y_3}
\end{array} \right),\left( \begin{array}{c}
1\\
0\\
0
\end{array} \right);\left( \begin{array}{l}
1\\
1\\
1
\end{array} \right),\left( \begin{array}{c}
0\\
1\\
0
\end{array} \right)} \right]_{\left( {\scriptstyle0\atop
{\scriptstyle0\atop
\scriptstyle1}} \right)}} = \]
\[ = \frac{{\left| {\left( {\begin{array}{*{20}{c}}
1&1\\
{1 + {x_3}}&1
\end{array}} \right)} \right| \cdot \left| {\left( {\begin{array}{*{20}{c}}
1&0\\
0&1
\end{array}} \right)} \right|}}{{\left| {\left( {\begin{array}{*{20}{c}}
1&0\\
{1 + {x_3}}&1
\end{array}} \right)} \right| \cdot \left| {\left( {\begin{array}{*{20}{c}}
1&1\\
0&1
\end{array}} \right)} \right|}} =  - {x_3}\]

\[{z_4} = {\left[ {E,C;F,B} \right]_A} = {\left[ {\left( \begin{array}{c}
1\\
0\\
0
\end{array} \right),\left( \begin{array}{c}
1\\
1 + {x_2}\\
{y_2}
\end{array} \right);\left( \begin{array}{c}
0\\
1\\
0
\end{array} \right),\left( \begin{array}{l}
1\\
1\\
1
\end{array} \right)} \right]_{\left( {\scriptstyle0\atop
{\scriptstyle0\atop
\scriptstyle1}} \right)}} = \]
\[ = \frac{{\left| {\left( {\begin{array}{*{20}{c}}
1&0\\
0&1
\end{array}} \right)} \right| \cdot \left| {\left( {\begin{array}{*{20}{c}}
1&1\\
{1 + {x_2}}&1
\end{array}} \right)} \right|}}{{\left| {\left( {\begin{array}{*{20}{c}}
1&1\\
0&1
\end{array}} \right)} \right| \cdot \left| {\left( {\begin{array}{*{20}{c}}
1&0\\
{1 + {x_2}}&1
\end{array}} \right)} \right|}} =  - {x_2}\]

\[{t_4} = {\left[ {C,A;F,B} \right]_E} = {\left[ {\left( \begin{array}{c}
1\\
1 + {x_2}\\
{y_2}
\end{array} \right),\left( \begin{array}{c}
0\\
0\\
1
\end{array} \right);\left( \begin{array}{c}
0\\
1\\
0
\end{array} \right),\left( \begin{array}{l}
1\\
1\\
1
\end{array} \right)} \right]_{\left( {\scriptstyle1\atop
{\scriptstyle0\atop
\scriptstyle0}} \right)}} \]
\[ = \frac{{\left| {\left( {\begin{array}{*{20}{c}}
{1 + {x_2}}&1\\
{{y_2}}&0
\end{array}} \right)} \right| \cdot \left| {\left( {\begin{array}{*{20}{c}}
0&1\\
1&1
\end{array}} \right)} \right|}}{{\left| {\left( {\begin{array}{*{20}{c}}
{1 + {x_2}}&1\\
{{y_2}}&1
\end{array}} \right)} \right| \cdot \left| {\left( {\begin{array}{*{20}{c}}
0&1\\
1&0
\end{array}} \right)} \right|}} = \frac{{ - {y_2}}}{{1 + {x_2} - {y_2}}}\]

So all coordinates related to the fourth chart are expressed explicitly from the other charts, 
hence we may disregard the fourth chart, and for the similar reason the fifth and the sixth chart.
Now we shall relate charts 1-3 to each other. We may use the symmetry, so for each formula we write 
there are two more symmetric formulas.

We have already noticed that $ - {x_3} = \left[ {E,A;D,B} \right]_C$. Hence
\[ - {x_3} = {\left[ {E,A;D,B} \right]_C} = \small\frac{{\left| {\left( {\begin{array}{*{20}{c}}
1&1&1\\
0&{1 + {x_2} + {z_2}}&{1 + {x_2}}\\
0&{{y_2} + {z_2}{t_2}}&{{y_2}}
\end{array}} \right)} \right| \cdot \left| {\left( {\begin{array}{*{20}{c}}
0&1&1\\
0&1&{1 + {x_2}}\\
1&1&{{y_2}}
\end{array}} \right)} \right|}}{{\left| {\left( {\begin{array}{*{20}{c}}
1&1&1\\
0&1&{1 + {x_2}}\\
0&1&{{y_2}}
\end{array}} \right)} \right| \cdot \left| {\left( {\begin{array}{*{20}{c}}
0&1&1\\
0&{1 + {x_2} + {z_2}}&{1 + {x_2}}\\
1&{{y_2} + {z_2}{t_2}}&{{y_2}}
\end{array}} \right)} \right|}} = \]
\[ = \frac{{\left| {\left( {\begin{array}{*{20}{c}}
{{z_2}}&{1 + {x_2}}\\
{{z_2}{t_2}}&{{y_2}}
\end{array}} \right)} \right| \cdot \left| {\left( {\begin{array}{*{20}{c}}
1&1\\
1&{1 + {x_2}}
\end{array}} \right)} \right|}}{{\left| {\left( {\begin{array}{*{20}{c}}
1&{1 + {x_2}}\\
1&{{y_2}}
\end{array}} \right)} \right| \cdot \left| {\left( {\begin{array}{*{20}{c}}
0&1\\
{{z_2}}&{1 + {x_2}}
\end{array}} \right)} \right|}} = \]
\[= \frac{{\left| {\left( {\begin{array}{*{20}{c}}
1&{1 + {x_2}}\\
{{t_2}}&{{y_2}}
\end{array}} \right)} \right| \cdot {x_2}}}{{\left( {{y_2} - {x_2} - 1} \right) \cdot \left| {\left( {\begin{array}{*{20}{c}}
0&1\\
1&{1 + {x_2}}
\end{array}} \right)} \right|}} = 
 \frac{{{x_2}\left( {{y_2} - {t_2} - {t_2}{x_2}} \right)}}{{1 - {x_2} - {y_2}}}\]

Also
\[{y_3} = {\left[ {A,C;E,D} \right]_B} = \frac{{\left| {\left( {\begin{array}{*{20}{c}}
0&1&1\\
0&0&1\\
1&0&1
\end{array}} \right)} \right| \cdot \left| {\left( {\begin{array}{*{20}{c}}
1&1&1\\
{1 + {x_2}}&{1 + {x_2} + {z_2}}&1\\
{{y_2}}&{{y_2} + {z_2}{t_2}}&1
\end{array}} \right)} \right|}}{{\left| {\left( {\begin{array}{*{20}{c}}
0&1&1\\
0&{1 + {x_2} + {z_2}}&1\\
1&{{y_2} + {z_2}{t_2}}&1
\end{array}} \right)} \right| \cdot \left| {\left( {\begin{array}{*{20}{c}}
1&1&1\\
{1 + {x_2}}&0&1\\
{{y_2}}&0&1
\end{array}} \right)} \right|}} = \]
\[ = \frac{{\left| {\left( {\begin{array}{*{20}{c}}
1&0&1\\
{1 + {x_2}}&{{z_2}}&1\\
{{y_2}}&{{z_2}{t_2}}&1
\end{array}} \right)} \right|}}{{\left| {\left( {\begin{array}{*{20}{c}}
0&1\\
{{x_2} + {z_2}}&1
\end{array}} \right)} \right| \cdot \left( {{y_2} - {x_2} - 1} \right)}} = \frac{{{z_2} - {y_2}{z_2} + {x_2}{z_2}{t_2}}}{{\left( {{x_2} + {z_2}} \right) \cdot \left( {1 + {x_2} - {y_2}} \right)}}\]

Now we shall express \[{\left[ {A,B;E,D} \right]_F}\] in two different ways:

\[{\left[ {A,B;E,D} \right]_F} = \small\frac{{\left| {\left( {\begin{array}{*{20}{c}}
1&1&1\\
0&{1 + {x_3}}&{1 + {x_3} + {z_3}}\\
0&{{y_3}}&{{y_3} + {z_3}{t_3}}
\end{array}} \right)} \right| \cdot \left| {\left( {\begin{array}{*{20}{c}}
0&1&1\\
1&1&{1 + {x_3} + {z_3}}\\
0&1&{{y_3} + {z_3}{t_3}}
\end{array}} \right)} \right|}}{{\left| {\left( {\begin{array}{*{20}{c}}
1&1&1\\
0&1&{1 + {x_3} + {z_3}}\\
0&1&{{y_3} + {z_3}{t_3}}
\end{array}} \right)} \right| \cdot \left| {\left( {\begin{array}{*{20}{c}}
0&1&1\\
1&{1 + {x_3}}&{1 + {x_3} + {z_3}}\\
0&{{y_3}}&{{y_3} + {z_3}{t_3}}
\end{array}} \right)} \right|}} = \]
\[ = \frac{{\left| {\left( {\begin{array}{*{20}{c}}
{1 + {x_3}}&{{z_3}}\\
{{y_3}}&{{z_3}{t_3}}
\end{array}} \right)} \right| \cdot \left( {1 - {y_3} - {z_3}{t_3}} \right)}}
{{\left( {{y_3} + {z_3}{t_3} - \left( {1 + {x_3} + {z_3}} \right)} \right) 
\cdot \left| {\left( {\begin{array}{*{20}{c}}
0&1&0\\
1&{1 + {x_3}}&{{z_3}}\\
0&{{y_3}}&{{z_3}{t_3}}
\end{array}} \right)} \right|}} = \]
\[= \frac{{\left| {\left( {\begin{array}{*{20}{c}}
{1 + {x_3}}&1\\
{{y_3}}&{{t_3}}
\end{array}} \right)} \right| \cdot \left( {1 - {y_3} - {z_3}{t_3}} \right)}}
{{\left( {{y_3} + {z_3}{t_3} - \left( {1 + {x_3} + {z_3}} \right)} \right) 
\cdot \left| {\left( {\begin{array}{*{20}{c}}
0&1&0\\
1&{1 + {x_3}}&1\\
0&{{y_3}}&{{t_3}}
\end{array}} \right)} \right|}} = \]
\[=\frac{{\left( {{t_3} + {t_3}{x_3} - {y_3}} \right)
\left( {1 - {y_3} - {z_3}{t_3}} \right)}}
{{\left( {1 + {x_3} + {z_3} - {y_3} - {z_3}{t_3}} \right) \cdot {t_3}}}\]

\[{\left[ {A,B;E,D} \right]_F} = \frac{{\left| {\left( {\begin{array}{*{20}{c}}
0&1\\
1&0
\end{array}} \right)} \right| \cdot \left| {\left( {\begin{array}{*{20}{c}}
1&1\\
1&{{y_2} + {z_2}{t_2}}
\end{array}} \right)} \right|}}{{\left| {\left( {\begin{array}{*{20}{c}}
0&1\\
1&{{y_2} + {z_2}{t_2}}
\end{array}} \right)} \right| \cdot \left| {\left( {\begin{array}{*{20}{c}}
1&1\\
1&0
\end{array}} \right)} \right|}} = 1 - {y_2} - {z_2}{t_2}\]

Therefore 
\[\frac{{\left( {{t_3} + {t_3}{x_3} - {y_3}} \right)\left( {1 - {y_3} - {z_3}{t_3}} \right)}}{{\left( {1 + {x_3} + {z_3} - {y_3} - {z_3}{t_3}} \right) \cdot {t_3}}} = 1 - {y_2} - {z_2}{t_2}\]

We may multiply by denominators and skip the terms of second order, so we get
\[ - {x_3} \approx -t_2 x_2\]
\[0 \approx {z_2}\]
\[{t_3} + {t_3}{x_3} - {y_3} - {t_3}{y_3} - {z_3}t_3^2 \approx {t_3}
\left( {1 - {y_2} - {z_2}{t_2} + {x_3} + {z_3} - {y_3} - {z_3}{t_3}} \right)\]
We may simplify these equations; in particular we may use that $z_2\approx 0$, 
and by symmetry $z_3\approx 0 \approx z_1$, which allows to simplify the third equation.
\[ {x_3} \approx t_2 x_2 \]
\[ 0 \approx {z_2} \]
\[ {y_3} \approx {t_3} {y_2} \]

After simplifying and using symmetry, we get 
\[\begin{array}{c}
\begin{array}{*{20}{c}}
\begin{array}{c}
{x_3} \approx {t_2}{x_2}\\
{y_3} \approx {t_3}{y_2}
\end{array}&\vline& \begin{array}{c}
{x_1} \approx {t_3}{x_3}\\
{y_1} \approx {t_1}{y_3}
\end{array}&\vline& \begin{array}{l}
{x_2} \approx {t_1}{x_1}\\
{y_2} \approx {t_2}{y_1}
\end{array}
\end{array}\\
{z_1},{z_2},{z_3} \approx 0
\end{array}\]

To get more equations we shall compute $\left[ {A,C;F,D} \right]_E$ in 3 different ways:
\[{\left[ {A,C;F,D} \right]_E} = {\left[ {\left( \begin{array}{c}
1\\
1 + {x_1}\\
{y_1}
\end{array} \right),\left( \begin{array}{c}
1\\
0\\
0
\end{array} \right);\left( \begin{array}{c}
1\\
1\\
1
\end{array} \right),\left( \begin{array}{c}
0\\
1\\
0
\end{array} \right)} \right]_{\left( {\scriptstyle0\atop
{\scriptstyle0\atop
\scriptstyle1}} \right)}} = \]
\[ = \frac{{\left| {\left( {\begin{array}{*{20}{c}}
1&1\\
{1 + {x_1}}&1
\end{array}} \right)} \right| \cdot \left| {\left( {\begin{array}{*{20}{c}}
1&0\\
0&1
\end{array}} \right)} \right|}}{{\left| {\left( {\begin{array}{*{20}{c}}
1&0\\
{1 + {x_1}}&1
\end{array}} \right)} \right| \cdot \left| {\left( {\begin{array}{*{20}{c}}
1&1\\
0&1
\end{array}} \right)} \right|}} =  - {x_1}\]

\[ {\left[ {A,C;F,D} \right]_E} = {\left[\small {\left( \begin{array}{c}
0\\
0\\
1
\end{array} \right),\left( \begin{array}{c}
1\\
1 + {x_2}\\
{y_2}
\end{array} \right);\left( \begin{array}{c}
0\\
1\\
0
\end{array} \right),\left( \begin{array}{c}
1\\
1 + {x_2} + {z_2}\\
{y_2} + {z_2}{t_2}
\end{array} \right)} \right]_{\left( {\scriptstyle1\atop
{\scriptstyle0\atop
\scriptstyle0}} \right)}} = \]
\[ = \frac{{\left| {\left( {\begin{array}{*{20}{c}}
0&1\\
1&0
\end{array}} \right)} \right| \cdot \left| {\left( {\begin{array}{*{20}{c}}
{1 + {x_2}}&{1 + {x_2} + {z_2}}\\
{{y_2}}&{{y_2} + {z_2}{t_2}}
\end{array}} \right)} \right|}}{{\left| {\left( {\begin{array}{*{20}{c}}
0&{1 + {x_2} + {z_2}}\\
1&{{y_2} + {z_2}{t_2}}
\end{array}} \right)} \right| \cdot \left| {\left( {\begin{array}{*{20}{c}}
1&1\\
{1 + {x_2}}&0
\end{array}} \right)} \right|}} = \frac{{ - \left| {\left( {\begin{array}{*{20}{c}}
{1 + {x_2}}&{{z_2}}\\
{{y_2}}&{{z_2}{t_2}}
\end{array}} \right)} \right|}}{{\left( {1 + {x_2} + {z_2}} \right) \cdot \left( {1 + {x_2}} \right)}} = \]
\[ = \frac{{{y_2}{z_2} - {z_2}{t_2}\left( {1 + {x_2}} \right)}}{{\left( {1 + {x_2} + {z_2}} \right) \cdot \left( {1 + {x_2}} \right)}}\]

\[ {\left[ {A,C;F,D} \right]_E} = {\left[ {\left( \begin{array}{c}
1\\
0\\
0
\end{array} \right),\left( \begin{array}{c}
0\\
1\\
0
\end{array} \right);\left( \begin{array}{c}
{z_7}{t_7}\\
{z_7}\\
1
\end{array} \right),\left( \begin{array}{c}
{y_7}\\
1\\
{y_7}
\end{array} \right)} \right]_{\left( {\scriptstyle0\atop
{\scriptstyle0\atop
\scriptstyle1}} \right)}} = \]
\[ = \frac{{\left| {\left( {\begin{array}{*{20}{c}}
1&{{z_7}{t_7}}\\
0&{{z_7}}
\end{array}} \right)} \right| \cdot \left| {\left( {\begin{array}{*{20}{c}}
0&{{y_7}}\\
1&1
\end{array}} \right)} \right|}}{{\left| {\left( {\begin{array}{*{20}{c}}
1&{{y_7}}\\
0&1
\end{array}} \right)} \right| \cdot \left| {\left( {\begin{array}{*{20}{c}}
0&{{z_7}{t_7}}\\
1&{{z_7}}
\end{array}} \right)} \right|}} = \frac{{{y_7}}}{{{t_7}}}\]

Therefore 
\[ - {x_1} = {\left[ {A,C;F,D} \right]_E} = \frac{{{y_2}{z_2} - {z_2}{t_2}\left( {1 + {x_2}} \right)}}{{\left( {1 + {x_2} + {z_2}} \right) \cdot \left( {1 + {x_2}} \right)}} \approx 0,\]
as we have seen that $z_2\approx 0$, and 
\[ - {x_1} = {\left[ {A,C;F,D} \right]_E} = \frac{{{y_7}}}{{{t_7}}}.\]
Hence $x_1\approx 0$ and by symmetry $x_2,x_3 \approx 0$. 
Also $y_7=-t_7x_1\approx 0$.

Now we write more formulas related to the seventh chart.
We have already seen that ${x_3} = {\left[ {E,A;D,B} \right]_C}$. Therefore 
\[ - {x_3} = {\left[ {E,A;D,B} \right]_C} = \frac{{\left| {\left( {\begin{array}{*{20}{c}}
0&{{y_7}}\\
1&{{y_7}}
\end{array}} \right)} \right| \cdot \left| {\left( {\begin{array}{*{20}{c}}
1&1\\
0&{{x_7}}
\end{array}} \right)} \right|}}{{\left| {\left( {\begin{array}{*{20}{c}}
0&1\\
1&{{x_7}}
\end{array}} \right)} \right| \cdot \left| {\left( {\begin{array}{*{20}{c}}
1&{{y_7}}\\
0&{{y_7}}
\end{array}} \right)} \right|}} = \frac{{\left| {\left( {\begin{array}{*{20}{c}}
0&1\\
1&1
\end{array}} \right)} \right| \cdot {x_7}}}{{ - 1 \cdot \left| {\left( {\begin{array}{*{20}{c}}
1&1\\
0&1
\end{array}} \right)} \right|}} = {x_7}\]

We also see that
\[{y_3} = {\left[ {A,C;E,D} \right]_B} = \frac{{\left| {\left( {\begin{array}{*{20}{c}}
1&0&1\\
0&0&{{x_7}}\\
0&1&{{x_7}}
\end{array}} \right)} \right| \cdot \left| {\left( {\begin{array}{*{20}{c}}
0&{{y_7}}&1\\
1&1&{{x_7}}\\
0&{{y_7}}&{{x_7}}
\end{array}} \right)} \right|}}{{\left| {\left( {\begin{array}{*{20}{c}}
1&{{y_7}}&1\\
0&1&{{x_7}}\\
0&{{y_7}}&{{x_7}}
\end{array}} \right)} \right| \cdot \left| {\left( {\begin{array}{*{20}{c}}
0&0&1\\
1&0&{{x_7}}\\
0&1&{{x_7}}
\end{array}} \right)} \right|}} = \]
\[ = \frac{{\left| {\left( {\begin{array}{*{20}{c}}
0&{{x_7}}\\
1&{{x_7}}
\end{array}} \right)} \right| \cdot \left| {\left( {\begin{array}{*{20}{c}}
0&{{y_7}}&1\\
1&1&{{x_7}}\\
0&{{y_7}}&{{x_7}}
\end{array}} \right)} \right|}}{{\left| {\left( {\begin{array}{*{20}{c}}
1&{{x_7}}\\
{{y_7}}&{{x_7}}
\end{array}} \right)} \right| \cdot 1}} = \frac{{\left| {\left( {\begin{array}{*{20}{c}}
0&1\\
1&1
\end{array}} \right)} \right| \cdot \left| {\left( {\begin{array}{*{20}{c}}
0&{{y_7}}&1\\
1&1&{{x_7}}\\
0&{{y_7}}&{{x_7}}
\end{array}} \right)} \right|}}{{\left| {\left( {\begin{array}{*{20}{c}}
1&1\\
{{y_7}}&1
\end{array}} \right)} \right| \cdot 1}} = \]
\[ = \frac{{ - \left( {{y_7} - {x_7}{y_7}} \right)}}{{1 - {y_7}}}\]

\[{y_1} = {\left[ {C,E;A,F} \right]_D} = \frac{{\left| {\left( {\begin{array}{*{20}{c}}
0&1&{{y_7}}\\
1&0&1\\
0&0&{{y_7}}
\end{array}} \right)} \right| \cdot \left| {\left( {\begin{array}{*{20}{c}}
0&{{z_7}{t_7}}&{{y_7}}\\
0&{{z_7}}&1\\
1&1&{{y_7}}
\end{array}} \right)} \right|}}{{\left| {\left( {\begin{array}{*{20}{c}}
0&{{z_7}{t_7}}&{{y_7}}\\
1&{{z_7}}&1\\
0&1&{{y_7}}
\end{array}} \right)} \right| \cdot \left| {\left( {\begin{array}{*{20}{c}}
0&1&{{y_7}}\\
0&0&1\\
1&0&{{y_7}}
\end{array}} \right)} \right|}} = \]
\[ = \frac{{ - {y_7} \cdot {z_7}\left( {{t_7} - {y_7}} \right)}}{{{y_7}
\left( {1 - {z_7}{t_7}} \right)}} = \frac{-{{z_7}\left( {{t_7} - {y_7}} \right)}}{{1 - {z_7}{t_7}}}\]

By multiplying by denominators and dismissing the terms of second order we get 
\[-x_3 = x_7\]
\[ y_3 \approx -y_7 \]
\[ y_1 \approx t_7z_7\]
Recall that 
\[z_1,z_2,z_3\approx 0\]
\[x_1,x_2,x_3\approx 0\]
\[y_7\approx 0\]
Hence 
\[x_7 \approx -x_3 \approx 0\]
\[y_3 \approx -y_7 \approx 0\]
\[{y_1} \approx {t_1}{y_3} \approx 0\]
\[{y_2} \approx {t_2}{y_1} \approx 0\]
And since $t_7\neq 0$, also 
\[ z_7= \frac {y_1}{t_7} \approx 0.\]
To summarize: $dx_i,dy_i,dz_i \approx 0$ for $i\in \{1,2,3,7\}$.
Therefore the cotangent space is spanned by 
$dt_1,dt_2,dt_3$ and $dt_7$.

\subsection{$\mathcal{F}.1'$}

The last subcase we verify is a degeneration of the previous subcase, 
when the triple ratio $\left\{ {P_1,P_5,P_3;P_2,P_6,P_4} \right\}$ tends to zero; 
if the triple ratio tends to infinity, it is symmetrical.

We choose the coordinates symmetrically, consistently with the previous subcase,
in order to use the formulas from the previous subcase.
\begin{center}
	\includegraphics{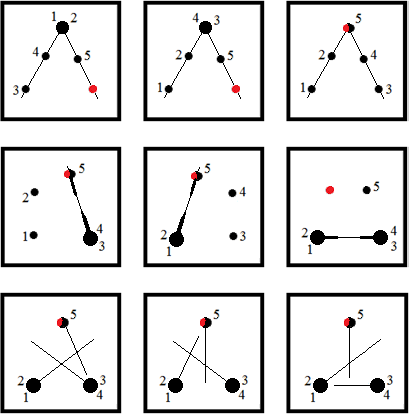}
\end{center}

\[\begin{array}{*{20}{c}}
\begin{array}{l}
{A_1} = \left( {1:1 + {x_1}:{y_1}} \right)\\
{B_1} = \left( {1:1 + {x_1} + {z_1}} \right.\\
\,\,\,\,\,\,\,\,\,\,\,\,\,\,\,\,\,\,\,\,\left. {:{y_1} + {z_1}{t_1}} \right)\\
{C_1} = \left( {1:0:0} \right)\\
{D_1} = \left( {0:1:0} \right)\\
{E_1} = \left( {0:0:1} \right)\\
{F_1} = \left( {1:1:1} \right)\\
{x_1},{y_1},{z_1} \sim 0
\end{array}&\vline& \begin{array}{l}
{A_2} = \left( {0:0:1} \right)\\
{B_2} = \left( {1:1:1} \right)\\
{C_2} = \left( {1:1 + {x_2}:{y_2}} \right)\\
{D_2} = \left( {1:1 + {x_2} + {z_2}} \right.\\
\,\,\,\,\,\,\,\,\,\,\,\,\,\,\,\,\,\,\,\,\left. {:{y_2} + {z_2}{t_2}} \right)\\
{E_2} = \left( {1:0:0} \right)\\
{F_2} = \left( {0:1:0} \right)\\
{x_2},{y_2},{z_2} \sim 0
\end{array}&\vline& \begin{array}{l}
{A_3} = \left( {1:0:0} \right)\\
{B_3} = \left( {0:1:0} \right)\\
{C_3} = \left( {0:0:1} \right)\\
{D_3} = \left( {1:1:1} \right)\\
{E_3} = \left( {1:1 + {x_3}:{y_3}} \right)\\
{F_3} = \left( {1:1 + {x_3} + {z_3}} \right.\\
\,\,\,\,\,\,\,\,\,\,\,\,\,\,\,\,\,\,\,\,\left. {:{y_3} + {z_3}{t_3}} \right)\\
{x_3},{y_3},{z_3} \sim 0
\end{array}
\end{array}\]

\[\begin{array}{*{20}{c}}
\begin{array}{l}
{A_4} = \left( {1:0:0} \right)\\
{B_4} = \left( {1:1:1} \right)\\
{C_4} = \left( {0:1:0} \right)\\
{D_4} = \left( {{x_4}{y_4}:1:{x_4}} \right)\\
{E_4} = \left( {0:0:1} \right)\\
{F_4} = \left( {{z_4}{t_4}:{z_4}:1} \right)\\
x_4,y_4,z_4,t_4\sim 0
\end{array}&\vline&  * &\vline&  * 
\end{array}\]

\[\begin{array}{*{20}{c}}
\begin{array}{l}
{A_7} = \left( {1:0:0} \right)\\
{B_7} = \left( {1:{x_7}:{x_7}} \right)\\
{C_7} = \left( {0:1:0} \right)\\
{D_7} = \left( {{y_7}:1:{y_7}} \right)\\
{E_7} = \left( {0:0:1} \right)\\
{F_7} = \left( {{z_7}{t_7}:{z_7}:1} \right)\\
x_7,y_7,z_7,t_7\sim 0
\end{array}&\vline& \begin{array}{l}
{A_8} = \left( {1:0:0} \right)\\
{B_8} = \left( {1:{z_8}{t_8}:{z_8}} \right)\\
{C_8} = \left( {0:1:0} \right)\\
{D_8} = \left( {{x_8}:1:{x_8}} \right)\\
{E_8} = \left( {0:0:1} \right)\\
{F_8} = \left( {{y_8}:{y_8}:1} \right)\\
x_8,y_8,z_8,t_8\sim 0
\end{array}&\vline& \begin{array}{l}
{A_9} = \left( {1:0:0} \right)\\
{B_9} = \left( {1:{y_9}:{y_9}} \right)\\
{C_9} = \left( {0:1:0} \right)\\
{D_9} = \left( {{z_9}:1:{z_9}{t_9}} \right)\\
{E_9} = \left( {0:0:1} \right)\\
{F_9} = \left( {{x_9}:{x_9}:1} \right)\\
x_9,y_9,z_9,t_9\sim 0
\end{array}
\end{array}\]

As in previous case, the coordinates from the charts 4-6 are expressed from 
the coordinates from the charts 1-3.

We can write again most of the formulas we had in the previous case for the same reasons, such as:
\[\begin{array}{c}
\begin{array}{*{20}{c}}
\begin{array}{c}
{x_3} \approx {t_2}{x_2}\\
{y_3} \approx {t_3}{y_2}
\end{array}&\vline& \begin{array}{c}
{x_1} \approx {t_3}{x_3}\\
{y_1} \approx {t_1}{y_3}
\end{array}&\vline& \begin{array}{l}
{x_2} \approx {t_1}{x_1}\\
{y_2} \approx {t_2}{y_1}
\end{array}
\end{array}\\
{z_1},{z_2},{z_3} \approx 0
\end{array}\]

As before we may also develop the formulas $x_1\approx 0$ and $y_7=-t_7x_1 \approx 0$, and now both 
formulas have 3 symmetric versions. We still have $x_7=-x_3$ and $y_7=-y_3$, so we may 
conclude immediately that $x_7\approx 0$ and $y_3\approx 0$, with two more symmetric versions 
for both last formulas
have 3 symmetric versions

\[\begin{array}{*{20}{c}}
\begin{array}{c}
{x_1} \approx 0\\
{y_7} \approx 0\\
x_7 \approx 0\\
y_1 \approx 0
\end{array}&\vline& \begin{array}{c}
{x_2} \approx 0\\
{y_8} \approx 0\\
x_8 \approx 0\\
y_2 \approx 0
\end{array}&\vline& \begin{array}{l}
{x_3} \approx 0\\
{y_9} \approx 0\\
x_9 \approx 0\\
y_3 \approx 0
\end{array}
\end{array}\]

The treatment of $z_7$ in the previous case relied on division by $t_7$ which is zero in this case,
so it can't be simply repeated. However  
\[{z_7} = {\left[ {E,C;F,B} \right]_A} = \frac{{\left| {\left( {\begin{array}{*{20}{c}}
0&{{x_9}}\\
1&1
\end{array}} \right)} \right| \cdot \left| {\left( {\begin{array}{*{20}{c}}
1&{{y_9}}\\
0&{{y_9}}
\end{array}} \right)} \right|}}{{\left| {\left( {\begin{array}{*{20}{c}}
0&{{y_9}}\\
1&{{y_9}}
\end{array}} \right)} \right| \cdot \left| {\left( {\begin{array}{*{20}{c}}
1&{{x_9}}\\
0&1
\end{array}} \right)} \right|}} = \frac{{ - {x_9} \cdot \left| {\left( {\begin{array}{*{20}{c}}
1&1\\
0&1
\end{array}} \right)} \right|}}{{\left| {\left( {\begin{array}{*{20}{c}}
0&1\\
1&1
\end{array}} \right)} \right| \cdot 1}} =  - {x_9},\]
and there are two more symmetric versions of this, hence also $z_i\approx 0$ for $i\in \left\{7,8,9\right\}$.
Therefore $dx_i=dy_i=dz_i=0$ for $i\in \left\{1,2,3,7,8,9\right\}$.
So, the cotangent space is spanned by $dt_1,dt_2,dt_3,dt_7,dt_8,dt_9$.

We may express the triple ratio $\left\{P_1,P_5,P_3;P_2,P_6,P_4\right\}$ by 
coordinates of the seventh chart:
\[\left\{ {A,E,C;B,F,D} \right\} = \frac{{\left| {\left( {\begin{array}{*{20}{c}}
0&{{z_7}{t_7}}&0\\
1&{{z_7}}&0\\
0&1&1
\end{array}} \right)} \right| \cdot \left| {\left( {\begin{array}{*{20}{c}}
0&1&1\\
0&{{x_7}}&0\\
1&{{x_7}}&0
\end{array}} \right)} \right| \cdot \left| {\left( {\begin{array}{*{20}{c}}
1&{{y_7}}&0\\
0&1&1\\
0&{{y_7}}&0
\end{array}} \right)} \right|}}
{{\left| {\left( {\begin{array}{*{20}{c}}
{{z_7}{t_7}}&1&0\\
{{z_7}}&0&0\\
1&0&1
\end{array}} \right)} \right| \cdot \left| {\left( {\begin{array}{*{20}{c}}
1&0&1\\
{{x_7}}&1&0\\
{{x_7}}&0&0
\end{array}} \right)} \right| \cdot \left| {\left( {\begin{array}{*{20}{c}}
{{y_7}}&0&0\\
1&0&1\\
{{y_7}}&1&0
\end{array}} \right)} \right|}} = \]
\[ = \frac{{ - {z_7}{t_7} \cdot {x_7} \cdot {y_7}}}{{{-z_7} \cdot {x_7} \cdot {y_7}}} =   {t_7}\]

We could do it similarly using the eighth or the ninth chart, and conclude that $t_7=t_8=t_9$.
Hence the cotangent space is spanned by $dt_1,dt_2, dt_3$ and $dt_7$, hence it is at most four-dimensional.

To summarize: our moduli space is smooth at all points.


\section{Bibliography.}
\begin{thebibliography}{99}

\bibitem[AC]{AC} E. Arbarello, M. Cornalba. Calculating cohomology groups of moduli spaces of curves via algebraic geometry.
Publ. math. de l'I.H.E.S., 88 (1998), 97-127.

\bibitem[FP]{FP} W. Fulton, R. Pandharipande, Notes on stable maps and quantum cohomology
Algebraic geometry - Santa Cruz (1995), 45-96.

\bibitem[DM]{DM}  P. Deligne, D. Mumford. The irreducibility of the space of curves of given genus. 
Publ. math. de l'I.H.E.S. 36 (1969), 75-109.

\bibitem [GKZ] {GKZ} I. M. Gelfand, M. Kapranov, and A. Zelevinsky. {\it Discriminants, resultants, and multidimensional determinants.} Springer Science \& Business Media (2008).

\bibitem[GH]{GH} P. Griffiths, J. Harris. {\it Principles of Algebraic geometry.}
 John Wiley \& Sons, (2014).

\bibitem[HM]{HM}  J. Harris, I. Morrison. Moduli of curves (1998), Grad. texts in math., 187.

\bibitem [H]{H} R. Hartshorne. {\it Algebraic geometry}. Vol. 52. Springer Science \& Business Media (2013).

\bibitem[Kap]{Kap} M. M. Kapranov. Chow quotients of Grassmannians I. in I. M. Gelfand Seminar,
Adv. Soviet Math. 16, Part 2, Amer. Math. Soc., Providence, 1993, 29 – 110.
MR 1237834 Number 2 (1992).

\bibitem[Keel]{Keel}  S. Keel. Intersection theory of moduli space of stable N-pointed curves of genus zero. Trans. of the Amer. Math. Soc., volume 330, Number 2 (1992).

\bibitem[KT]{KT} S. Keel, J. Tevelev. Geometry of Chow quotients of Grassmannians. Duke Math. J. 134(2), 259–311 (2006) 

\bibitem[Kn]{Kn}  F. F. Knudsen. Projectivity of the moduli space of stable curves, II: the stacks $\mathcal{M}_{g,n}$. Math. Scand. 52 (1983),
1225-1265.









\end{thebibliography}
\end{document}